\documentclass[10pt,a4paper]{article}
\usepackage{amstext}
\usepackage{array}
\usepackage{amsfonts}
\usepackage{amssymb}
\usepackage{graphicx}
\usepackage{epstopdf}
\usepackage{algorithm}     
\usepackage{algpseudocode} 
\usepackage{url}
\usepackage{caption}
\usepackage{subfig}
\usepackage{csquotes}
\usepackage{tikz}
\usetikzlibrary{arrows}
\usepackage{mathtools}
\usepackage{bbm}
\usepackage{amsthm}
\usepackage{pgfplots}
\usepackage{multirow}
\usepackage{framed}

\usepgfplotslibrary{external}
\DeclareRobustCommand{\tikzcaption}[1]{\tikzset{external/export next=false}#1}
\DeclareRobustCommand{\tikzref}[1]{\tikzcaption{\resizebox{!}{\refsize}{\ref{#1}}}}

\newcommand{\vertiii}[1]{{\left\vert\kern-0.25ex\left\vert\kern-0.25ex\left\vert #1 
   \right\vert\kern-0.25ex\right\vert\kern-0.25ex\right\vert}}

\newtheorem{Definition}{Definition}[section]
\newtheorem{theorem}{Theorem}[section]
\newtheorem{theorem1}{Theorem}[section]
\newtheorem{theorem2}{Theorem}[section]
\newtheorem{proposition}[theorem]{Proposition}
\newtheorem{lemma}[theorem1]{Lemma}
\newtheorem{remark}[theorem2]{Remark}

\usepackage{cleveref}

\usepackage{geometry}
\title{An optimization-based registration approach to    geometry reduction}
\author{Tommaso Taddei}
\date{}
\usepackage{fancyhdr}  
  
\fancypagestyle{plain}{
\fancyhead{}
\fancyhead[C]{\hfill Submitted to International Journal of Numerical Methods in Engineering, November 2022}
}

\begin{document}

\maketitle

\begin{abstract}
We develop and assess an optimization-based approach to parametric  geometry reduction.
Given a family of parametric domains, we aim to determine a parametric diffeomorphism $\Phi$ that maps a fixed reference domain $\Omega$ into each element of the family, for different values of the parameter; the ultimate goal of our study is to determine an effective tool for parametric projection-based model order reduction of partial differential equations in parametric geometries.
For practical problems in engineering, explicit parameterizations of the geometry are likely unavailable: for this reason, our approach   takes  as  inputs a reference mesh of $\Omega$ and a point cloud $\{y_i^{\rm raw}\}_{i=1}^Q$ that belongs to the boundary of the target domain $V$ and returns a bijection $\Phi$ that approximately maps $\Omega$ in $V$.
We propose a two-step procedure:
given the point clouds $\{ x_j  \}_{j=1}^N\subset \partial \Omega$ and 
$\{y_i^{\rm raw}\}_{i=1}^Q \subset \partial V$, we first resort to a point-set registration algorithm to determine the  displacements  $\{ v_j  \}_{j=1}^N$ such that the deformed point cloud
$\{y_j:= x_j+v_j  \}_{j=1}^N$ approximates $\partial V$; then, we solve a nonlinear non-convex optimization problem to build a mapping $\Phi$ that is bijective from  $\Omega$ in $\mathbb{R}^d$ and (approximately) satisfies $\Phi(x_j) = y_j$ for $j=1,\ldots,N$.
We present a rigorous mathematical analysis to justify our approach; we further present thorough numerical experiments to show the effectiveness of the proposed method.
 \end{abstract}

\section{Introduction}
\label{sec:intro}
A broad range of problems in science and engineering requires parametric studies to assess the influence of geometry on the solution to a given  
 partial differential equation (PDE). Given the family of 
parametric  domains $\{ \Omega_{\mu} : \mu \in \mathcal{P} \}$ in $\mathbb{R}^d$, where $\mu$ denotes a vector of parameters in the compact set $\mathcal{P}$ in $\mathbb{R}^p$, we seek a reference domain $\Omega$ and a parameterized mapping $\Phi:\Omega\times \mathcal{P} \to \mathbb{R}^d$ such that $\Phi_{\mu}(\Omega)$ approximates $\Omega_{\mu}$ --- in a sense to be defined --- and is bijective in $\Omega$ for all $\mu \in \mathcal{P}$. The aim of this work is to develop and analyze an optimization-based
geometry registration and reduction (GRR) framework for the construction of the mapping $\Phi$ and for its rapid evaluation for new configurations.
In this work,
we refer to the problem of determining a mapping from the reference  domain $\Omega$ to the target domain $\Omega_{\mu}$ as \emph{geometry registration}; 
on the other hand, we refer to the process of determining a low-rank representation of the parametric  mapping $\Phi$ as \emph{geometry reduction}.

The ultimate goal of our study is parametric (projection-based) model order reduction (pMOR, \cite{hesthaven2016certified,quarteroni2015reduced}) for PDEs in parametric geometries.
Given the parametric PDE model $\mathfrak{L}_{\mu}( \mathbbm{u}_{\mu}   ) = 0$ defined over  the domain $\Omega_{\mu}$, pMOR methods rely on  a parametric mapping $\Phi$ to recast the problem over a parameter-independent configuration $\Omega$ ; then, pMOR methods resort to the projection of the PDE over  suitable parameter-independent low-dimensional trial and test spaces to approximate 
the mapped field $\mathbbm{u}_{\mu}\circ \Phi_{\mu}$. As discussed in Remark \ref{remark:MtD_vs_DtM}, the formulation in the reference domain might be performed implicitly (discretize-then-map) or explicitly (map-then-discretize): 
for finite element (FE) discretizations, 
both approaches rely on the introduction of  a mesh $\mathcal{T}^{\rm hf} = \left( \{ x_i^{\rm hf} 
 \}_{i=1}^{N_{\rm hf}}, \texttt{T} \right)$  with nodes 
 $ \{x_i^{\rm hf}  \}_{i=1}^{N_{\rm hf}}$ and connectivity matrix $\texttt{T}$.
Despite the many contributions to the field 
\cite{lassila2014model,lovgren2006reduced,manzoni2017efficient,manzoni2012model,rozza2008reduced},
development of rapid and reliable geometry reduction techniques for large deformations is still a challenging task.
 
Effective GRR techniques should fulfill the following requirements.
\begin{itemize}
\item
\emph{Hidden parameterization.}
In many applications (e.g., modeling of 
biological processes in patient-specific geometries), the parameterization is unknown and the only available information is a set of points $\{ y_i^{\rm raw}   \}_{i=1}^Q$ that belong to the boundary of $\Omega_{\mu}$. GRR algorithms should thus take as input the mesh $\mathcal{T}_{\rm hf}$ and the point cloud $\{ y_i^{\rm raw}   \}_{i=1}^Q \subset \partial \Omega_{\mu}$ and return the deformed mesh $\Phi_{\mu}(  \mathcal{T}_{\rm hf}  )$
with deformed nodes
 $ \{\Phi_{\mu}(  x_i^{\rm hf})  \}_{i=1}^{N_{\rm hf}}$ and fixed connectivity matrix $\texttt{T}$.
\item
\emph{Real-time computations.}
In many scenarios, the computation of the deformed mesh should be performed extremely rapidly, possibly with little computing resources. 
To achieve this goal,  effective reduction strategies based on the low-rank approximation of the sought  mapping should be developed.
 \item
 \emph{Control of mesh quality and geometric error.}
GRR methods should explicitly ensure that the geometric error --- measured by the Hausdorff distance between $\Phi_{\mu}(\Omega)$ and $\Omega_{\mu}$ (see  section \ref{sec:analysis}) ---
 and the quality of the deformed mesh meet user-defined tolerances. 
 \end{itemize}
 
 In this work, 
we propose a two-step procedure:
given the point clouds $\{ x_j  \}_{j=1}^N\subset \partial \Omega$ and 
$\{y_i^{\rm raw}\}_{i=1}^Q \subset \partial \Omega_{\mu}$, we first resort to a point-set registration
(PSR,\cite{horaud2010rigid,ma2015robust,ma2018nonrigid,maiseli2017recent,zitova2003image}) algorithm
 to determine the displacements $\{ v_j  \}_{j=1}^N$ such that the deformed point cloud
$\{y_j:= x_j+v_j  \}_{j=1}^N$ approximates $\partial \Omega_{\mu}$; then, we seek a mapping $\Phi_{\mu}$ that is bijective in $\Omega$ and (approximately) satisfies $\Phi_{\mu}(x_j) = y_j$ for $j=1,\ldots,N$.
To accomplish the first task, we rely on the coherent point drift 
 (CPD, \cite{myronenko2006non,myronenko2010point})
method: CPD is a well-established non-rigid PSR technique that is broadly used for image processing and pattern recognition applications  and is also the point of departure of several more recent methods.
On the other hand, following   \cite{taddei2020registration}, 
 we recast the problem of finding the mapping $\Phi$ as an optimization problem:
  the optimization framework allows us to directly control the geometry error and also the quality of the deformed mesh; furthermore, the spectacular advances in non-convex optimization enable the solution to highly-nonlinear and non-convex problems in a reasonable time frame.
Finally, we resort to proper orthogonal decomposition
(POD, \cite{volkwein2011model})  to identify a low-rank approximation space for the mapping and ultimately speed up the registration process. 
 
We present a rigorous mathematical analysis of the geometry error that provides a rigorous foundation for the optimization approach.
In more detail, we show that for smooth domains, under mild assumptions on $\Phi$ the Hausdorff distance between $\Omega_{\mu}$ and $\Phi_{\mu}(\Omega)$ is controlled  by $\max_{j=1,\ldots,N} \|  \Phi_{\mu}(x_j) - y_j \|_2$ (see Propositions \ref{th:quasi_hausdorff_bound} and \ref{th:geo_theorem}); we also show that the same result does not hold for Lipschitz domains with corners.
Moreover, we establish a connection between elasticity-based mesh morphing and registration methods (cf. Remark \ref{remark:link_elasticity}). 
Finally, we investigate the performance of several variants of the optimization statement considered in 
\cite{taddei2020registration}  for a number of two-dimensional model problems.
  
This work is also linked to optimization-based mesh morphing, 
which rely on the solution to a (convex) optimization problem to determine the deformed mesh  (see the review  
\cite{staten2011comparison}); here, we emphasize the application to parametric problems and the use of a global --- as opposed to local, compactly-supported --- basis for the mapping.
Furthermore, our approach exploits the connection 
between GRR and the problem of point-set  registration   in bounded domains, which was previously considered in  \cite{iollo2022mapping} for MOR applications.

We further remark that several authors have proposed mesh moving strategies for unsteady PDEs that rely on multiple incremental solutions to a suitable elasticity problem with mesh-dependent properties
(see  \cite{tezduyar1992new,tonon2021linear}  and references therein). In the  parametric setting, we can interpret this strategy as feed-forward maps of the form
$\Phi = \Phi_N\circ\ldots \circ \Phi_1$ where the $i$-th map $\Phi_i$ maps $\Omega_i$ into $\Omega_{i+1}$ with $\Omega_1=\Omega$ and $\Omega_{N+1}=\Omega_{\mu}$. In this work, we consider instead linear approximations (cf. \eqref{eq:Noperator_square}) that are determined by solving a single highly-nonlinear optimization problem: this choice provides a natural framework for linear-subspace data compression methods (such as POD) and is thus well-suited for MOR applications. 

The outline of the paper is as follows.
In  section \ref{sec:method}, we present the methodology: first, we consider the problem 
of point-set registration in bounded domains; second, we present the extension to geometry registration; third, we comment on dimensionality reduction.
In section \ref{sec:analysis}, we present the analysis of the geometry error;
in section \ref{sec:numerics}, we present the results of the numerical experiments;
finally, in section \ref{sec:conclusions}, we provide a summary of the contributions and of the numerical  results . Several appendices   conclude the paper.

\section{Optimization-based registration}
\label{sec:method}

In  section \ref{sec:notation}, we introduce relevant notation and preliminary definitions that are useful for the subsequent discussion. Then, in  section \ref{sec:opt_statement_PBR}, we discuss the problem of point-set registration in bounded domains:
we provide a complete discussion for tensorized domains and we comment on the extension to more general domains. 
In  section \ref{sec:geometry_reduction}, we discuss the problem of geometry registration and reduction. In the remainder, the spatial dimension $d$ is either $d=2$ or $d=3$.

\subsection{Notation and preliminary definitions}
\label{sec:notation}
We denote  by $x$ a generic point of the domain $\Omega$, and by $\mathbf{n}(x)$ the outward normal to $\Omega$ at $x\in \partial \Omega$. We define the identity map 
$\texttt{id}:\Omega\to \Omega$ such that $\texttt{id}(x)=x$; then, given the mapping $\Phi:\Omega\to \mathbb{R}^d$, we denote by $u = \Phi - \texttt{id}$ the displacement field and by $J(\Phi) = {\rm det} (\nabla \Phi )$ the Jacobian determinant.
We denote by $\mathcal{T}_{\rm hf}$ a finite element (FE) mesh of $\Omega$ with $N_{\rm e}$ elements $\{ \texttt{D}_k  \}_{k=1}^{N_{\rm e}}$; we define the reference element 
$\widehat{\texttt{D}} = \{ \widetilde{x}\in (0,1)^d: \sum_{i=1}^d (\widetilde{x})_i<  1  \}$  and the shape functions
$\{ \ell_i \}_{i=1}^{n_{\rm lp}}$  of the polynomial space $\mathbb{P}_{\texttt{p}}(\widehat{\texttt{D}})$ associated with the nodes  $\{  \widetilde{x}_i  \}_{i=1}^{n_{\rm lp}}$. Then, we represent the elemental mappings $\Psi_k^{\rm hf}$
from  $\widehat{\texttt{D}}$ to $\texttt{D}_k$ for $k=1,\ldots,N_{\rm e}$ as
  \begin{equation}
\label{eq:psi_mapping}
 {\Psi}_k^{\rm hf}( \widetilde{x}  )
 =
 \sum_{i=1}^{n_{\rm lp}} \;
 {x}_{i,k}^{\rm hf}  \; \ell_i(\widetilde{x}),
\end{equation}
where $\{ {x}_{i,k}^{\rm hf}   := {\Psi}_k^{\rm hf}( \widetilde{x}_i): \, i=1,\ldots,n_{\rm lp}, k=1,\ldots,N_{\rm e}  \}$ are the nodes of the mesh.
To simplify the presentation, we here consider the same basis for both FE fields and FE elemental mappings (isoparametric FE); the discussion below can be trivially extended to sub- and sup-parametric discretizations.
We say that the mapping $\Phi$ is bijective if $J(\Phi)$ is strictly positive in $\Omega$; on the other hand, we say that $\Phi$ is discretely bijective with respect to the mesh $\mathcal{T}_{\rm hf}$ if the elemental mappings
\begin{equation}
\label{eq:psi_mapping_def}
 {\Psi}_{\Phi,k}^{\rm hf}( \widetilde{x}  )
 =
 \sum_{i=1}^{n_{\rm lp}} \;
 \Phi \left( {x}_{i,k}^{\rm hf} \right) \; \ell_i(\widetilde{x})
\end{equation}
are all bijective. 
Note that  the function $\Phi:\Omega\to \mathbb{R}^d$ is bijective with respect to $\mathcal{T}_{\rm hf}$ if
$\Phi_{\rm hf}:\Omega\to \mathbb{R}^d$ given by
\begin{equation}
\label{eq:mapping_DtM}
\Phi_{\rm hf}(x)
=
\Psi_{\Phi,k}^{\rm hf} \left( 
\left(
\Psi_{k}^{\rm hf} 
\right)^{-1}
(x)
\right),
\quad
\forall \, x\in \texttt{D}_k, \;\;
k=1,\ldots,N_{\rm e},
\end{equation}
is a bijection.
Bijectivity and discrete bijectivity are two independent conditions, which are both important for model reduction  of parametric systems in parameterized geometries as discussed in the following remark.

\begin{remark}
\label{remark:MtD_vs_DtM}
To clarify the distinction between bijectivity and discrete bijectivity, consider the problem of approximating the solution to the following Laplace equation in the family parameterized domains $\{\Omega_{\mu}  : \mu\in \mathcal{P} \}$:
\begin{equation}
\label{eq:laplace_eq_remark}
-\Delta U_{\mu} = f_{\mu} \;\; {\rm in } \;\; \Omega_{\mu},
\;\;
U_{\mu} |_{\partial \Omega_{\mu}}=0,
\end{equation}
for some $f_{\mu}\in H^{-1}(\Omega_{\mu})$. 
Given the reference domain $\Omega$, we denote by $\Phi:\Omega\times \mathcal{P} \to \mathbb{R}^d$ a mapping such that $\Phi_{\mu}(\Omega) = \Omega_{\mu}$ for all $\mu\in \mathcal{P}$,  by $\mathcal{T}_{\rm hf}$ a FE mesh of $\Omega$, and by $\mathcal{X}_{\rm hf} \subset H_0^1(\Omega)$ the corresponding FE space. To  devise a FE discretization of  \eqref{eq:laplace_eq_remark}, we first observe that the mapped field
$\widetilde{U}_{\mu} := U_{\mu}\circ \Phi_{\mu}$ satisfies
\begin{equation}
\label{eq:MtD_paradigm}
\int_{\Omega} J( \Phi_{\mu}) \left(
\nabla \Phi_{\mu}^{-1} (\nabla \Phi_{\mu})^{-\star} \nabla 
\widetilde{U}_{\mu} \cdot \nabla v - f_{\mu}\circ \Phi_{\mu} \, v
\right) \, dx = 0 \quad
\forall \, v\in H_0^1(\Omega),
\end{equation}
where  $\nabla \Phi_{\mu}^{\star}$ denotes the   transpose  of  $\nabla \Phi_{\mu}$;
then, we proceed to discretize \eqref{eq:MtD_paradigm} by projecting the equations over the FE space $\mathcal{X}_{\rm hf}$. Alternatively, we might consider the formulation: find $U_{\mu}^{\rm hf} \in \mathcal{X}_{\rm hf,\mu}$ such that
\begin{equation}
\label{eq:DtM_paradigm}
\sum_{k=1}^{N_{\rm e}} \;\;
\int_{ \texttt{D}_{\mu,k}   }   \left(
 \nabla 
U_{\mu}^{\rm hf} \cdot \nabla v - f_{\mu} \, v
\right) \, dx = 0 \quad
\forall \, v\in \mathcal{X}_{\rm hf,\mu},
\end{equation}
where $\mathcal{X}_{\rm hf,\mu}$ is the FE space associated with the deformed FE mesh
$\Phi_{\mu}(  \mathcal{T}_{\rm hf}    )$, while 
$ \texttt{D}_{\mu,k} = \{  \Psi_{\mu,k}^{\rm hf}(\widetilde{x}) : \widetilde{x} \in \widehat{\texttt{D}} \}$ denotes the $k$-th deformed element of the mesh.
Clearly, \eqref{eq:MtD_paradigm}  and  \eqref{eq:DtM_paradigm} are equivalent  if and only if
$\Phi =\Phi_{\rm hf}$ (see \eqref{eq:mapping_DtM}).
Both approaches \eqref{eq:MtD_paradigm}  and  \eqref{eq:DtM_paradigm} are broadly used for the treatment of parameterized geometries in pMOR:
following  \cite{taddei2021discretize}, we refer to
\eqref{eq:MtD_paradigm} as to map-then-discretize approach and to  \eqref{eq:DtM_paradigm} as to discretize-then-map approach. Note that \eqref{eq:MtD_paradigm} requires the bijectivity of the mapping $\Phi$, while  \eqref{eq:DtM_paradigm}  requires the discrete bijectivity with respect to the mesh  $\mathcal{T}_{\rm hf}$.
\end{remark}

\subsection{Point-set registration in bounded domains}
\label{sec:opt_statement_PBR}

\subsubsection{Optimization statement}
Given the domain $\Omega\subset \mathbb{R}^d$, and the point clouds $\{x_i \}_{i=1}^N$  and $\{y_i \}_{i=1}^N$, we consider the problem of finding a mapping $\Phi = \texttt{id} + u$ such that $\Phi(x_i) = y_i$ for $i=1,\ldots,N$ and $\Phi$ is bijective from $\Omega$ in itself. Towards this end, we define the operator 
$B: [C(\Omega)]^d \to \mathbb{R}^{d\cdot N}$ and the vector $\mathbf{z}$ such that 
\begin{equation}
\label{eq:interpolation_operators}
B\Phi 
\, =\,
\left[
\begin{array}{l}
\left(   \Phi(x_1) \right)_1 \\
\vdots \\
\left(   \Phi(x_N) \right)_1 \\
\vdots \\
\left(   \Phi(x_N) \right)_d \\
\end{array}
\right],
\quad
\mathbf{z}
=
\left[
\begin{array}{l}
\left(   y_1 \right)_1 \\
\vdots \\
\left(   y_1 \right)_1 \\
\vdots \\
\left(   y_N \right)_d \\
\end{array}
\right].
\end{equation}
Then, we can state the constrained optimization statement as follows:
given 
(i) the set of vector-valued functions $\mathcal{X}$, 
(ii) the objective function $\mathfrak{f}^{\rm obj}:\mathcal{X} \to \mathbb{R}$,
 we seek $\Phi\in \mathcal{X}$ to minimize
\begin{equation}
\label{eq:constrained_registration}
\min_{\Phi\in \mathcal{W}} \;
\mathfrak{f}^{\rm obj}(\Phi)  \quad
{\rm s.t.} \;\;
B\Phi  = \mathbf{z}.
\end{equation}
Clearly, the performance of registration relies on the choice of $\mathcal{W}$, $\mathfrak{f}^{\rm obj}$ and also on the optimization algorithm that is employed to determine local minima of \eqref{eq:constrained_registration} . We address these points in the next sections.

In many contexts, the hard constraint 
$B\Phi  = \mathbf{z}$ is either inappropriate (e.g., if the point clouds are subject to noise) or
might lead to an ill-posed statement (e.g., if two points are coincident). Exploiting classical theory in inverse problems. we consider two regularization strategies for \eqref{eq:constrained_registration}:
\begin{equation}
\label{eq:tykhonov_registration}
\min_{\Phi\in \mathcal{W}} \;
\xi \, \mathfrak{f}^{\rm obj}(\Phi)     + \frac{1}{2}  \|  B\Phi  - \mathbf{z} \|_2^2,
\end{equation}
and
\begin{equation}
\label{eq:morozov_registration}
\min_{\Phi\in \mathcal{W}} \;
 \mathfrak{f}^{\rm obj}(\Phi)   
  \quad
{\rm s.t.} \;\;  \|  B\Phi  - \mathbf{z} \|_{\infty} \leq \delta.
\end{equation}
Following taxonomy from the inverse problem literature 
\cite{ivanov2013theory}, we refer to \eqref{eq:tykhonov_registration} as to Tykhonov-regularized registration statement and to \eqref{eq:morozov_registration} as to Morozov-regularized registration statement. 
We observe that both approaches rely on the introduction of an hyper-parameter --- $\xi$ and $\delta$, respectively.

Morozov regularization allows to directly control the geometry error through the choice of the parameter $\delta$; if $\mathcal{W}$ is a linear space,
the choice of the norm $\| \cdot  \|_{\infty}$ in  \eqref{eq:morozov_registration} 
leads to  an   optimization statement with linear inequality constraints, which can be effectively tackled using interior point methods.
Tykhonov regularization does not enable direct control of the geometry error and requires a careful selection of the parameter $\xi$; however, the resulting optimization problem is unconstrained and, in our experience,  can be solved more efficiently than \eqref{eq:morozov_registration} using quasi-Newton methods, particularly for the choice of the squared 2-norm in \eqref{eq:tykhonov_registration}. 

\subsubsection{Choice of the objective function and of the search space $\mathcal{W}$ for tensorized domains}
\label{sec:rectangular_domain}

The choice of the space $\mathcal{W}$ and of the objective function $\mathfrak{f}^{\rm obj}$
 should ensure that 
(i) $B$ is a bounded functional over $\mathcal{W}$, 
(ii) any minimizer $\Phi$ of the problem \eqref{eq:tykhonov_registration} or \eqref{eq:morozov_registration} 
should be  a bijection in $\Omega$ and/or should be discretely bijective with respect to the target mesh $\mathcal{T}_{\rm hf}$, and 
(iii) the space $\mathcal{W}$ 
is well-suited to fit a broad range of bijections in $\Omega$ (\emph{high expressive power}).
In this section, we discuss the choice of  $\mathcal{W}$  and  $\mathfrak{f}^{\rm obj}$ for tensorized domains. Without loss of generality, we consider $\Omega=(0,1)^d$.

If $\Omega$ is a tensorized domain, we can prove that $\Phi=\texttt{id} + u$ is bijective in $\Omega$ if $u\cdot \mathbf{n} |_{\partial \Omega} = 0$ and $J(\Phi)$ is strictly positive in $\overline{\Omega}$ 
(cf. Proposition 2.3 in \cite{taddei2020registration}). We can thus choose the affine space 
$\mathcal{W}$ as 
\begin{equation}
\label{eq:inf_dimensional_calW}
\mathcal{W} = \left\{
\texttt{id} +\varphi : 
\varphi \in [H^2(\Omega)]^d, 
\varphi  \cdot \mathbf{n} |_{\partial \Omega} = 0
\right\}.
\end{equation}
Note that $\mathcal{W}$ is a linear affine space of continuous functions for $d\leq 3$; more precisely, 
in two and three dimensions $H^2$ is the Sobolev space of minimal regularity that is contained in the space of continuous functions \cite{adams2003sobolev}.
For practical computations, we replace the infinite-dimensional space $H^2(\Omega)$ in \eqref{eq:inf_dimensional_calW} with a finite-dimensional tensorized polynomial space of dimension $(n_{\rm lp}+1)^d$, that is
$\mathcal{W}_{\rm hf}  =
\texttt{id} + \mathcal{U}_{\rm hf}$ with 
\begin{equation}
\label{eq:polynomial_space}
\mathcal{U}_{\rm hf} 
= \left\{
 \varphi   \in [\mathbb{Q}_{n_{\rm lp}} ]^d \; : \;
\varphi  \cdot \mathbf{n} |_{\partial \Omega} = 0
\right\},
\quad
\mathbb{Q}_{n_{\rm lp}}
={\rm span} \left\{
\varphi(x) = \prod_{i=1}^d
 \ell_i((x)_i) \; : \;
  \ell_1,\ldots,  \ell_d\in \mathbb{P}_{n_{\rm lp}} (\mathbb{R})
\right\},
\end{equation}
where $\mathbb{P}_{n_{\rm lp}}(\mathbb{R})$ denotes the space of one-dimensional polynomials of degree up to $n_{\rm lp}$. We observe that $\mathcal{W}_{\rm hf}$ is dense in $\mathcal{W}$ in the limit $n_{\rm lp}\to \infty$; we further observe that the space 
$\mathcal{U}_{\rm hf}$ is of dimension $M=d (n_{\rm lp}+1)^d - 2 d (n_{\rm lp}+1)^{d-1}$. 
Given the orthonormal basis of $\mathcal{U}_{\rm hf}$  $\{\psi_m\}_{m=1}^M$, we can introduce the operator
$\texttt{N}: \mathbb{R}^M \to \mathcal{W}_{\rm hf}$ such that
\begin{equation}
\label{eq:Noperator_square}
\texttt{N}(x; \mathbf{a})
\, : = \,
x + \sum_{m=1}^M \left( \mathbf{a} \right)_m \psi_m(x);
\end{equation}
note that $\texttt{N}$ is a linear affine operator and 
 $\texttt{N}(\cdot; \mathbf{0}) = \texttt{id}$.

We propose two different objective functions: the former is designed to enforce bijectivity of $\Phi$, while the latter is designed to enforce discrete bijectivity.
In both cases, we add the regularization term
\begin{equation}
\label{eq:penalty_term} 
\mathfrak{P}(\Phi)
=
\frac{1}{2 |\Omega|}   
| \Phi  |_{H^2(\Omega)}^2
=
\frac{1}{2 |\Omega|}
\int_{\Omega} 
\sum_{i,j,k=1}^d \, 
\left(
\partial_{i,j} (\Phi)_k \right)^2 \, dx ,
\end{equation}
which is designed to promote smoothness: the regularization $\mathfrak{P}$ \eqref{eq:penalty_term} ensures that the $H^2$ norm of any global minimum of \eqref{eq:tykhonov_registration} or \eqref{eq:morozov_registration} is uniformly bounded in $H^2$ for $M\to \infty$.
Given $\epsilon,C_{\rm exp}>0$, we define the objective function $\mathfrak{f}^{\rm obj}  = 
\mathfrak{f}_{\rm jac} + \mathfrak{P}$ such that
\begin{equation}
\label{eq:exp_jac}
\mathfrak{f}_{\rm jac}(\Phi)
=
\frac{1}{|\Omega|}
\int_{\Omega} \; {\rm exp} \left(
\frac{\epsilon - J(\Phi) } {C_{\rm exp}}
\right) \, dx.
\end{equation}  
The parameter $\epsilon$ provides a weak lower bound for the value of the Jacobian determinant of the mapping $\Phi$: provided that
$e^{\epsilon/C_{\rm exp}} \gg 1$, the objective function
\eqref{eq:exp_jac}  forces
the optimizer to ensure the condition
$J(\Phi) \gtrsim \epsilon$ everywhere in $\Omega$;
the $H^2$ penalization enforces smoothness of $\Phi$ and is also instrumental to ensure bijectivity (see
   \cite[section 2.2]{taddei2020registration}).
Given the constant $\kappa_{\rm msh}$, following   \cite{knupp2001algebraic} (see also   \cite{zahr2018optimization}), we 
define the objective function $\mathfrak{f}^{\rm obj} = 
\mathfrak{f}_{\rm msh} + \mathfrak{P}$ such that
\begin{equation}
\label{eq:exp_mesh}
\mathfrak{f}_{\rm msh}({\Phi}   )
=
\frac{1}{|\Omega|}
\sum_{k=1}^{N_{\rm e}} \; \; 
|  \texttt{D}_k  |
\int_{  \widehat{\texttt{D}}  } 
{\rm exp} \left(
{\kappa_{\rm msh}
\, - \, 
 q_{\Phi,k} }
\right)
 \; d\widetilde{x},
\quad
 q_{\Phi,k} : = 
 \frac{1}{d^2}
\left(
\frac{ \|  \nabla \Psi_{\Phi,k}^{\rm hf} \|_{\rm F}^2  }{( {\rm det} ( \nabla  \Psi_{\Phi,k}^{\rm hf}   )    )_+^{2/d}}\right)^2,
\end{equation}
where 
$\Psi_{\Phi,k}^{\rm hf} $ is defined in \eqref{eq:psi_mapping_def},
$\| \cdot \|_{\rm F}$ is the Frobenius norm and  
$(\cdot)_+ = \max( \cdot, 0 )$.
The ratio $ q_{\Phi,k} $ measures the degree of anisotropy of the mesh: it can thus be interpreted as a measure of the quality of the deformed mesh. 
Note that if $\Psi_{\Phi,k}$ is perfectly isotropic, 
$ q_{\Phi,k} = 1$.
Similarly to \eqref{eq:exp_jac}, $\mathfrak{f}^{\rm obj} = 
\mathfrak{f}_{\rm msh}  + \mathfrak{P}$ strongly penalizes deformations for which 
$ q_{\Phi,k} $ exceeds $\kappa_{\rm msh}$ for some element $\texttt{D}_k$ of the mesh, and promotes smoothness of the mapping in regions where 
$q_{\Phi,k} < \kappa_{\rm msh}$.

We observe that \eqref{eq:exp_jac} and to \eqref{eq:exp_mesh} depend on several hyper-parameters --- $\epsilon,C_{\rm exp}$ and $\kappa_{\rm msh}$.
In our numerical results, we investigate the sensitivity with respect to $\epsilon$ and 
$\kappa_{\rm msh}$; on the other hand, we set $C_{\rm exp}=0.025\epsilon$.
We  further observe that neither \eqref{eq:exp_jac} nor \eqref{eq:exp_mesh} are convex: therefore, we do not expect that  problems \eqref{eq:tykhonov_registration} and \eqref{eq:morozov_registration} admit a unique solution;
we provide extensive numerical investigations of the registration approach based on the two choices of the objective function in section \ref{sec:numerics}.
In the remainder, we refer to \eqref{eq:exp_jac} and to \eqref{eq:exp_mesh} as to ``exp-jac'' objective and ``exp-mesh'' objective, respectively.
Finally, we observe that 
\eqref{eq:penalty_term}, 
\eqref{eq:exp_jac} and   \eqref{eq:exp_mesh}
--- and thus also  the objective function 
$\mathfrak{f}^{\rm obj}$ --- 
are invariant under translations and rotations, that is
$$
\mathfrak{P}(  \mathbf{b} + \mathbf{R}  {\Phi}   )
=
\mathfrak{P}(  {\Phi}   ),
\quad
\mathfrak{f}_{\rm jac}(  \mathbf{b} + \mathbf{R}  {\Phi}   )
=
\mathfrak{f}_{\rm jac}(  {\Phi}   ),
\quad
\mathfrak{f}_{\rm msh}(  \mathbf{b} + \mathbf{R}  {\Phi}   )
=
\mathfrak{f}_{\rm msh}(  {\Phi}   ),
$$
for any $\mathbf{b}\in \mathbb{R}^d$ and any rotation $\mathbf{R}\in \mathbb{R}^{d\times d}$.

\begin{remark}
\label{remark:link_elasticity}
We might attempt to establish a direct link between optimization-based registration and elasticity-based mesh morphing.
We might indeed formally consider the statement
\begin{equation}
\label{eq:elasticity_based_a}
\min_{\Phi\in \mathcal{W}} \mathfrak{f}^{\rm el}(\Phi)
\quad {\rm s.t.} 
\quad
B\Phi = \mathbf{z},
\end{equation}
where $\mathfrak{f}^{\rm el}(\Phi) = \int_{\Omega} \psi^{\rm se}(x; \Phi) \, dx$ is the integral over $\Omega$ of the strain energy $\psi^{\rm se}$. To provide concrete references that are considered in the numerical experiments, we might consider the linear isotropic model:
\begin{equation}
\label{eq:isotropic_linear_strain}
\psi^{\rm se}(\cdot; \Phi= \texttt{id} + u)
\, =\,
\lambda_1 (\nabla\cdot u)^2 + 2 \lambda_2 \| 
\frac{1}{2} \left( 
\nabla u + (\nabla u)^{\star}  \right) \|_{\rm F}^2,
\end{equation}
or the nonlinear Neohookean model (for $d=2$):
\begin{equation}
\label{eq:neohookean_strain}
\psi^{\rm se}(\cdot; \Phi)
\, =\,
\frac{1}{2} \lambda_2 
\| \nabla  \Phi \|_{\rm F}^2
\, - \,
\lambda_2 {\rm log} ( J(\Phi)    )
\, + \,
\lambda_1 \left(
{\rm log} ( J(\Phi)   )
\right)^2,
\end{equation}
where $\lambda_1,\lambda_2$ are the Lamè constants.

We observe that the models \eqref{eq:isotropic_linear_strain} and \eqref{eq:neohookean_strain} are naturally defined in $H^1(\Omega)$, which is not contained in the space of continuum functions in $\Omega$ for $d\geq 1$. On the other hand,  the objective functions \eqref{eq:exp_jac} and \eqref{eq:exp_mesh} involve derivatives of the deformation gradient $\nabla \Phi$: therefore, they cannot be interpreted as strain energies of a suitable hyper-elastic material.
\end{remark}

\subsubsection{Generalization to arbitrary domains}
\label{sec:approx_class}
The extension of the registration approach to non-tensorized domains is particularly challenging due to the difficulty to ensure the condition $\Phi(\Omega)=\Omega$.
From an algorithmic standpoint, we should identify a new operator $\texttt{N}$ (cf. \eqref{eq:Noperator_square}) and a new  penalization term $\mathfrak{P}$ (cf. \eqref{eq:penalty_term}). 
A potential strategy based on a spectral element approximation to deal with this case is proposed in 
\cite{taddei2021registration}. In this paper, we do not address this issue.

\subsubsection{Generalization to unsorted point clouds}
\label{sec:unsorted_point_clouds}

In many applications, the reference and target point clouds 
$\{ x_i \}_{i=1}^N$ and
$\{ y_i^{\rm raw} \}_{i=1}^Q$ might be of different cardinality (i.e., $N\neq Q$) and/or might not be properly sorted.
We are hence in need of a point-set registration (PSR) algorithm that aligns the two point clouds, that is, it finds $\{ v_i \}_{i=1}^N$ such that 
$\{ y_i:=x_i+v_i \}_{i=1}^N$  approximates --- in a sense to be defined --- the target point cloud
$\{ y_i^{\rm raw}  \}_{i=1}^Q$.
If we define the matrices
$\mathbf{X} = [x_1,\ldots,x_N]^{\star}  \in \mathbb{R}^{N\times d}$,
$\mathbf{Y}^{\rm raw}  = [y_1^{\rm raw} ,\ldots,y_Q^{\rm raw} ]^{\star}  \in \mathbb{R}^{Q\times d}$ and
$\mathbf{Y} = [x_1+v_1,\ldots,y_N+v_N]^{\star}  \in \mathbb{R}^{N\times d}$, we can formalize  an abstract PSR algorithm as
$$
\mathbf{Y}
=\texttt{PSR} \left( \mathbf{X}, \mathbf{Y}^{\rm raw}  \right).
$$
The literature on PSR is extremely vast  and a thorough review of the subject is beyond the scope of the present work. In this paper, we rely on the coherent point drift 
(CPD, \cite{myronenko2006non,myronenko2010point})
algorithm: CPD is a well-established non-rigid PSR technique that is broadly used  and is also the point of departure of several more recent methods.

The CPD method seeks mappings of the form
\begin{equation}
\label{eq:CPD_functional}
\Phi^{\rm cpd}(x) =
x + \sum_{i=1}^N w_i \phi\left( \|x-x_i\|_2 \right),
\quad
{\rm with} \;\;
\phi(r) = {\rm exp} \left(  -\frac{1}{2\beta^2} r^2  \right).
\end{equation}
where  
$\mathbf{W} = [w_1,\ldots,w_N]^{\star}\in \mathbb{R}^{N\times d}$  is  a matrix of  coefficients 
that is chosen through an iterative procedure;
therefore, the CPD method seeks mappings in the affine space
\begin{equation}
\label{eq:CPD_affine_space}
\mathcal{W}_{\rm hf}^{\rm cpd} =\texttt{id}
\, + \, \mathcal{U}_{\rm hf}^{\rm cpd},
\quad
{\rm with} \;\;
\mathcal{U}_{\rm hf}^{\rm cpd} ={\rm span} 
\left\{
\phi(\|  \cdot - x_i \|_2) \mathbf{e}_j: 
i=1,\ldots,N, j=1,\ldots,d
\right\},
\end{equation}
where $\mathbf{e}_1,\ldots,\mathbf{e}_d$ are the elements of the canonical basis of $\mathbb{R}^d$.
Since the CPD methodology is at this stage well-established, we postpone the detailed description to Appendix \ref{sec:CPD_appendix}. Here, we anticipate that the method can naturally cope with low-dimensional search spaces $\mathcal{U}_{M}^{\rm cpd} \subset H^2(\mathbb{R}^d) $: we discuss in 
section  \ref{sec:dimensionality_reduction} the importance of this observation.

\subsection{Geometry registration and reduction}
\label{sec:geometry_reduction}

\subsubsection{Optimization statement}
\label{sec:optimization_georeg}
We adapt the framework of section \ref{sec:opt_statement_PBR} to the problem of geometry registration: given the reference domain $\Omega$ and the point cloud $\{   y_i^{\rm raw}\}_{i=1}^M$ that belongs to the boundary of the domain $V$, we seek a bijective mapping $\Phi$ from $\Omega$ to $V$. This task is propedeutic to the problem of geometry reduction, which is addressed in the next section. 
In the remainder, we denote by $\{   x_i \}_{i=1}^N$ the discretization of the boundary $\partial \Omega$ that is used for registration; we also denote by  $\mathcal{T}_{\rm hf}$ a FE mesh in $\Omega$.

We  introduce the hyper-rectangle $\Omega_{\rm box}$ that strictly contains $\Omega$ 
and the point cloud $\{y_i^{\rm raw} \}_i$; 
 we define the affine space $\mathcal{W}_{\rm hf} = \texttt{id} + \mathcal{U}_{\rm hf}$ with $\mathcal{U}_{\rm hf} = [\mathbb{Q}_{n_{\rm lp}}]^d$;
given $\Phi\in \mathcal{W}_{\rm hf}$, we further define the counterparts of the quantities in section \ref{sec:rectangular_domain}
$\mathfrak{P}(\Phi) = \frac{1}{2|\Omega_{\rm box}|} | \Phi  |_{H^2(\Omega_{\rm box})}^2$
(cf. \eqref{eq:penalty_term}), 
$\mathfrak{f}_{\rm jac}(\Phi) = \frac{1}{|\Omega_{\rm box}| }  \int_{\Omega_{\rm box}}  {\rm exp} \left(  \frac{\epsilon - J(\Phi)}{C_{\rm exp}} \right)  \, dx$ (cf. \eqref{eq:exp_jac}), and 
$\mathfrak{f}_{\rm msh}(\Phi)$ as in  \eqref{eq:exp_mesh}. Then, we consider the following two-step registration procedure: first, we resort to PSR to determine the deformed points
$\{ y_i = x_i + v_i   \}_{i=1}^N$ (cf. section \ref{sec:unsorted_point_clouds}); second, we solve a registration problem in $\Omega_{\rm box}$ with inputs
$\{   x_i \}_{i=1}^N$ and $\{   y_i \}_{i=1}^N$. Towards this end, we consider the Tykhonov-regularized registration statement
\begin{equation}
\label{eq:mesh_morphing_optimization_tikhonov}
\min_{\Phi\in \mathcal{W}_{\rm hf}} \xi \mathfrak{f}^{\rm obj}(\Phi) + 
\frac{1}{2}
\|  B \Phi - \mathbf{z} \|_2^2;
\end{equation}
and the Morozov-regularized registration statement
\begin{equation}
\label{eq:mesh_morphing_optimization_morozov}
\min_{\Phi\in \mathcal{W}_{\rm hf}}  \mathfrak{f}^{\rm obj}(\Phi)  
\quad
{\rm s.t.} \;\;
\|  B \Phi - \mathbf{z} \|_{\infty}\leq \delta.
\end{equation}
The operator $B$ and the vector $\mathbf{z}$ are defined as in \eqref{eq:interpolation_operators};
as in section \ref{sec:rectangular_domain}, 
the objective $\mathfrak{f}^{\rm obj}$ is given by $\mathfrak{f}^{\rm obj} = \mathfrak{f}_{\rm jac}+\mathfrak{P}$ (``exp-jac'') or by 
$\mathfrak{f}^{\rm obj} = \mathfrak{f}_{\rm msh}+\mathfrak{P}$ (``exp-mesh'').

In view of the numerical assessment, we further introduce the formulation 
\begin{equation}
\label{eq:inverted_formulation}
\min_{\Phi\in \mathcal{W}_{\rm hf}}  
\frac{1}{2} \|  B \Phi - \mathbf{z} \|_2^2 \, + \,
\xi  \mathfrak{P}(\Phi)  
\quad
{\rm s.t.} \;\;
\mathfrak{f}_{\star}(\Phi)  
-\delta_{\rm con}
\leq  0,
\end{equation}
which  was considered in  \cite{taddei2020registration}. Statement \eqref{eq:inverted_formulation} reads as a nonlinear constrained optimization problem with quadratic objective function and nonlinear non-convex inequality constraint. We envision that the advantage of \eqref{eq:inverted_formulation} is that it allows to control more explicitly the  minimum of the Jacobian determinant (for $\star = {\rm jac}$) or the 
mesh quality (for $\star = {\rm msh}$); on the other hand, it does not explicitly control the geometry error (unlike \eqref{eq:mesh_morphing_optimization_morozov}) and it involves a nonlinear constraint. We compare performance of \eqref{eq:inverted_formulation} 
with \eqref{eq:mesh_morphing_optimization_tikhonov} and  \eqref{eq:mesh_morphing_optimization_morozov} 
in the numerical experiments.

We here consider a polynomial space $\mathcal{W}_{\rm hf}$ in \eqref{eq:mesh_morphing_optimization_tikhonov}, \eqref{eq:mesh_morphing_optimization_morozov},
\eqref{eq:inverted_formulation}. This choice is motivated by the need to compute second-order derivatives for the regularization term and  by the strong approximation properties of polynomials in moderate dimensions. A thorough assessment of other approximation classes is  beyond the scope of the present work.

We remark that the goal of geometry registration is to find a mapping $\Phi$  such that $\Phi(\Omega)$ is close (in the sense of Hausdorff) to the domain $V$. Statements \eqref{eq:mesh_morphing_optimization_tikhonov} and \eqref{eq:mesh_morphing_optimization_morozov}
\eqref{eq:inverted_formulation}
  control the error
$B \Phi - \mathbf{z}$ in a convenient norm; in section \ref{sec:analysis}, we discuss under what conditions the control of the difference $B \Phi - \mathbf{z}$  ensures appropriate reconstruction of $V$.

\subsubsection{Dimensionality reduction}
\label{sec:dimensionality_reduction}
Given the results of $n_{\rm train}$ geometry registration problems $\{ \Phi_i = \texttt{id} + u_i  \}_{i=1}^{n_{\rm train}}$ for different configurations, we might apply dimensionality reduction techniques to identify a low-dimensional representation of the mapping. The objective of dimensionality reduction is twofold:
first, by reducing the number of unknowns, we might significantly speed up computations;
second, dimensionality reduction techniques might also act as low-pass filters to avoid overfitting, especially for noisy datasets.

In this work, we rely on 
proper orthogonal decomposition
(POD, \cite{volkwein2011model}) based on the $H^2(\Omega_{\rm box})$ inner product to determine a $M$-dimensional reduced space $\mathcal{U}_M = {\rm span} \{ \psi_m \}_{m=1}^M \subset \mathcal{U}_{\rm hf}$. 
We first assemble the Gramian matrix $\mathbf{C} \in \mathbb{R}^{n_{\rm train} \times n_{\rm train}}$ such that
$\mathbf{C}_{i,j} = ( u_i ,  u_j)_{H^2(\Omega_{\rm box})}$ for $i,j=1,\ldots,n_{\rm train}$; then, we compute the eigenpairs of $\mathbf{C}$ $\{ (\lambda_i, \boldsymbol{\psi}_i)  \}_i$ such that $\lambda_1\geq \ldots \geq \lambda_{n_{\rm train}} \geq 0$. Finally, we define the POD space
$$
\psi_m = \frac{1}{\sqrt{\lambda_m}} \sum_{i=1}^{n_{\rm train}} \, u_i \left(  \boldsymbol{\psi}_m \right)_i,
\quad
m=1,\ldots,M.
$$
The size $M$ of the reduced space is chosen 
according to the criterion
\begin{equation}
\label{eq:POD_cardinality_selection}
M := \min \left\{
m: \,  {\rm EC}_m(\boldsymbol{\lambda})  \geq   1 - tol_{\rm pod} 
\right\},
\quad {\rm where} \;
{\rm EC}_m(\boldsymbol{\lambda}) = 
\left(  \sum_{j=1}^{n_{\rm train}} \lambda_j \right)^{-1} 
\sum_{i=1}^m \lambda_i.
\end{equation} 
 In the remainder, ${\rm EC}_m(\boldsymbol{\lambda}) $ is referred to as the relative energy content of the POD space $\mathcal{U}_M$.
 
We observe that the cost of solving the optimization problems 
\eqref{eq:mesh_morphing_optimization_tikhonov}, \eqref{eq:mesh_morphing_optimization_morozov}, \eqref{eq:inverted_formulation}  with reduced search space $\texttt{id} + \mathcal{U}_M$ might still be significant due to the need to compute the  function \eqref{eq:exp_mesh} over all elements of the mesh. To address this issue, we might extend reduced quadrature techniques developed in the pMOR framework\cite{farhat2021computational,yano2021model} to the registration framework. This extension is   beyond the scope of the present paper.
  
As discussed above, we first resort to PSR to determine the sorted deformed points $\{y_i\}_{i=1}^N$ and then to one of the statements   
  \eqref{eq:mesh_morphing_optimization_tikhonov}, \eqref{eq:mesh_morphing_optimization_morozov}, \eqref{eq:inverted_formulation} to obtain the bijective  mapping. So far, we proposed to apply POD to reduce the cost of the second step of the procedure. Even if the latter  dominates by far computational costs, since the problem of interest is inherently non-convex, the difference in the ansatz spaces employed by CPD (cf. \eqref{eq:CPD_affine_space}) and by registration might lead to inaccurate performance.
  To investigate this issue, in the numerical results we consider two distinct strategies:  \textbf{(i)} apply CPD with $\mathcal{W}_{\rm hf}^{\rm cpd}$;
  \textbf{(ii)} apply CPD with 
\begin{equation}
\label{eq:CPD_reduced_space}
 \mathcal{W}_{M}^{\rm cpd}
 =
 \texttt{id} +
 {\rm span}
 \left\{
 \psi_m 
 \right\}_{m=1}^M.
\end{equation}
The second strategy requires to  slightly   modify the original CPD formulation that are described in Appendix \ref{sec:CPD_appendix}.
 
\section{Analysis}
\label{sec:analysis}

We denote by $U$ and $V$ two Lipschitz domains that are isomorphic to the unit ball and are compactly embedded in the hyper-cube $\Omega_{\rm box}$. We denote by ${\rm dist}_{\rm H}(U,V)$ the Hausdorff distance between $U$ and $V$ such that
\begin{equation}
\label{eq:hausdorff_distance}
{\rm dist}_{\rm H}(U,V) = {\rm max} \left\{
\sup_{x\in U} {\rm dist}(x,V), \;\;\
\sup_{x\in V} {\rm dist}(x,U)
\right\},
\end{equation}
where ${\rm dist}(x,\omega) = \inf_{y\in \omega} \| x - y \|_2$ for any $x\in \mathbb{R}^d$ and any measurable set $\omega\subset \mathbb{R}^d$.
Given the target domain $V$ and the reference domain $\Omega$, the goal of geometry registration algorithms is to determine a bijective mapping $\Phi$ from $\Omega$ to $\mathbb{R}^d$ such  that the Hausdorff distance  ${\rm dist}_{\rm H}(\Phi(\Omega),V)$
is below a given tolerance.

In order to bridge the gap between the target \eqref{eq:hausdorff_distance} and the computational procedure introduced in this paper, we introduce the non-symmetric boundary distance:
\begin{equation}
\label{eq:quasi_hausdorff_distance}
{\rm dist}_{\rm bnd}(U; V) = 
\sup_{x\in \partial U} {\rm dist}(x,\partial V).
\end{equation}
Next Proposition clarifies the link between \eqref{eq:quasi_hausdorff_distance} and  the geometry constraint
in \eqref{eq:mesh_morphing_optimization_tikhonov} and
\eqref{eq:mesh_morphing_optimization_morozov}; the proof is contained in  Appendix \ref{sec:proofs}.

\begin{proposition}
\label{th:quasi_hausdorff_bound}
Let $\{ x_i \}_{i=1}^N \subset \partial U$ be an $\epsilon$-cover of $\partial \Omega$, that is $\sup_{x\in \partial \Omega}  {\rm dist} \left(x, \{ x_i \}_i \right) = \epsilon$. Assume that $\Phi:\Omega\to \mathbb{R}^d$ is a Lipschitz bijective map with Lipschitz constant $K$ and let $\{ y_i \}_{i=1}^N \subset \partial V$. Then, we have 
\begin{equation}
\label{eq:quasi_hausdorff_bound}
{\rm dist}_{\rm bnd}(\Phi(U); V)
\leq
\max_{i=1,\ldots,N} \|  \Phi(x_i) - y_i\|_2 + K \epsilon. 
\end{equation}
Note that 
$\max_{i=1,\ldots,N} \|  \Phi(x_i) - y_i\|_2 \leq
\sqrt{d} 
\max_{i=1,\ldots,N} \|  \Phi(x_i) - y_i\|_{\infty}
=
\sqrt{d} 
\| B \Phi - \mathbf{z}  \|_{\infty}
$, which is the constraint of the optimization statement \eqref{eq:mesh_morphing_optimization_morozov}.
\end{proposition}

The objective of this section is to establish conditions for which \eqref{eq:quasi_hausdorff_distance} is equal to \eqref{eq:hausdorff_distance}.
Towards this end, we first present a number of definitions and preliminary results (cf.   section \ref{sec:justification_pre}),  then we present the main proposition for smooth domains  (cf.  section \ref{sec:justification_main}) and 
we 
show through the vehicle of a two-dimensional  example that \eqref{eq:hausdorff_distance} and \eqref{eq:quasi_hausdorff_distance} cannot be equal for Lipschitz domains with corners
(cf.   section \ref{sec:Lipschitz_analysis}).
Proofs of the Lemmas below are contained in  Appendix \ref{sec:proofs}.

\subsection{Preliminary results and definitions}
\label{sec:justification_pre}

Given $\delta>0$, we define the $\delta$-neighborhood of $\partial U$ as
\begin{equation}
\label{eq:neighU}
{\rm Neigh}_{\delta}(\partial U) : =
\left\{
x \in \Omega_{\rm box} :
{\rm dist} \left(
x, \partial U
\right) < \delta
\right\},
\end{equation}
and the tubular neighborhood of  $\partial U$  as
\begin{equation}
\label{eq:tubular_neigh}
{\rm Neigh}_{\delta}^{\rm t}( \partial U  ) := \{ y \in \mathbb{R}^d: \, 
y=x+ t \mathbf{n}({x}),
\:
|t| < \delta , \;  x \in \partial U \}.
\end{equation}
Next Lemma summarizes important properties of ${\rm Neigh}_{\delta}( \partial U  ), {\rm Neigh}_{\delta}^{\rm t}( \partial U  )$ that will be used below  --- we remark that the second statement of Lemma \ref{th:neigh_prel} is an immediate consequence of the Weyl's tube formula\cite{gray2003tubes}. Given the set $A$, we denote by $|A|_{(d-1)}$ (resp., $|A|_{(d)}$) the $(d-1)$-dimensional (resp., $d$-dimensional) measure of  $A$ --- to fix ideas, if $d=2$,  $|A|_{(1)}$ is the length of the curve $A$ embedded in $\mathbb{R}^2$, while  $|A|_{(2)}$ is the area of the region $A$. If not specified otherwise, we shall assume that $\delta$ is small enough so that ${\rm Neigh}_{\delta}(\partial U)$ is compactly embedded in $\Omega_{\rm box}$.

\begin{lemma}
\label{th:neigh_prel}
Let $U \subset \Omega_{\rm box}$ be  a $C^1$ domain isomorphic to the unit ball. Then, 
(i) ${\rm Neigh}_{\delta}(\partial U) = {\rm Neigh}_{\delta}^{\rm t}( \partial U  )$, and
(ii) 
\begin{equation}
\label{eq:bound_neigh}
| {\rm Neigh}_{\delta}(\partial U) |_{(d)} \;  \leq 
\left\{
\begin{array}{ll}
\displaystyle{2\delta | \partial U  |_{(d-1)}
} & {\rm if} \;\; d=2, \\[3mm]
\displaystyle{2\delta | \partial U  |_{(d-1)}
\,+\,
\frac{8 \pi^2}{3} \delta^3
} & {\rm if} \;\; d=3. \\
\end{array}
\right.
\end{equation} 
\end{lemma}

We introduce a special class of domains.
\begin{Definition}
\label{def:delta_regularity}
We say that $U$ is $\delta$-regular if for any given $x \in \partial U$ any  set $V$  diffeomorphic to the unit ball  such that
$\partial V \subset {\rm Neigh}_{\delta}(\partial U) \setminus \mathcal{B}_{\delta}(x)$ is contained in ${\rm Neigh}_{\delta}(\partial U)$.
\end{Definition}

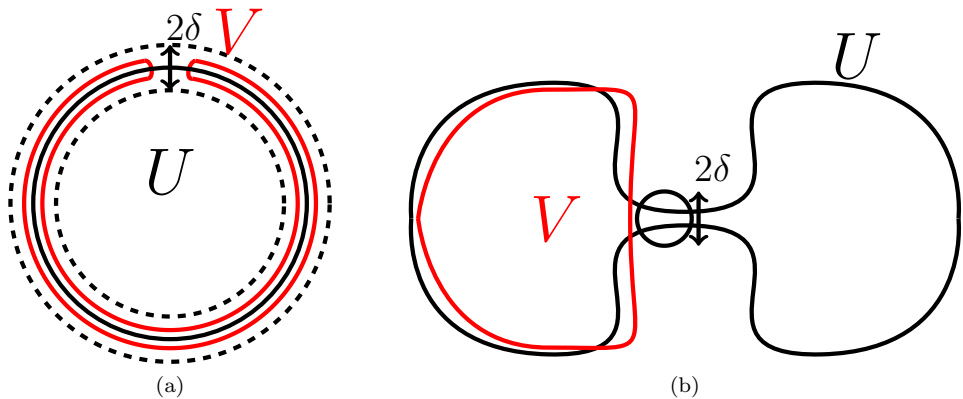
\begin{figure}[h!]
\centering

 \subfloat[] {
\begin{tikzpicture}[scale=0.6]
  \draw [red,ultra thick,domain=100:440,samples=300] plot ({3.2*cos(\x)}, {3.2*sin(\x)});

 \draw [red,ultra thick,domain=100:440,samples=300] plot ({2.8*cos(\x)}, {2.8*sin(\x)});

\draw[<->,ultra thick] (0,2.5) --(0,3.5);
\coordinate [label={right:  {\Large{$2 \delta$}}}] (E) at (-0.3,3.9) ;

\coordinate (G) at (-0.48621489746,2.75746170843);
\coordinate (R) at (-0.55567416853,  3.15138480964);

\coordinate (G2) at (0.48621489746,2.75746170843);
\coordinate (R2) at (0.55567416853,  3.15138480964);

  \draw [black,ultra thick,dashed,domain=0:360,samples=300] plot ({2.5*cos(\x)}, {2.5*sin(\x)});
 
   \draw [black,ultra thick,dashed,domain=0:360,samples=300] plot ({3.5*cos(\x)}, {3.5*sin(\x)});

\draw [red,ultra thick] (G) to[out=10,in=10] (R);
\draw [red,ultra thick] (G2) to[out=190,in=190] (R2) ;

\coordinate [label={above:  {\Huge {${\color{red} V}$}}}] (E) at (1.5,3.1) ;
\coordinate [label={above:  {\Huge {${U}$}}}] (E) at (0,0) ;
  \draw [ultra thick] (0,0) circle [radius=3]; 
   
\end{tikzpicture}
}
~~~~
\subfloat[] {
\begin{tikzpicture}[scale=0.9]

\draw [ultra thick,black] plot [smooth, tension=2] coordinates { (0,0)   (2,2)  (4,0.1)  (6,2) (8,0)};

%

\draw [ultra thick,black] plot [smooth, tension=2] coordinates { (0,0)   (2,-2)  (4,-0.1)  (6,-2) (8,0)};

 \draw [ultra thick,red] plot [smooth, tension=2] coordinates { (0.1,0)   (2,1.9)  (3.2,0) (2,-1.9) (0.1,0)   }; 

\coordinate [label={right:  {\Huge {${U}$}}}] (E) at (6,2.4) ;
\coordinate [label={right:  {\Huge {${\color{red}  V}$}}}] (E) at (1.6,0) ;

\draw [black,ultra thick,domain=0:360,samples=300] plot ({3.7+0.4*cos(\x)}, {0.4*sin(\x)});

\draw[<->,ultra thick] (4.2,-0.4) --(4.2,0.4);
\coordinate [label={right:  {\Large{$2\delta$}}}] (E) at (4,0.75) ;

\end{tikzpicture}
 }

 \caption{geometric interpretation of   Definition  \ref{def:delta_regularity}.
Dashed lines in (a)    denote the boundary of ${\rm Neigh}_{\delta}(\partial U)$.
}
\label{fig:polyB}
\end{figure}

Figure \ref{fig:polyB} provides an interpretation of   Definition \ref{def:delta_regularity}:
for a circle $U= \mathcal{B}_1(0)$ (cf.   Figure \ref{fig:polyB}(a)), we find by inspection that any $V$ such that $\partial V \subset {\rm Neigh}_{\delta}(\partial U) \setminus \mathcal{B}_{\delta}({x})$ for some ${x} \in \partial U$ must be contained in the annulus $\mathcal{B}_{1+\delta}({0}) \setminus  \mathcal{B}_{1-\delta}({0})$, provided that $0<\delta<1$.
For any fixed $\delta>0$, we might construct  smooth domains with bounded  curvature that are not  $\delta$-regular  (cf.  example in   Figure \ref{fig:polyB}(b)).
The example in Figure \ref{fig:polyB}(b) also shows that $\delta$-regularity is a global property of $\partial U$ in the sense that does not uniquely depend on local properties (e.g., curvature, Lipschitz constant) of the boundary.

Next three Lemmas offer further interpretations and properties of 
$\delta$-regular domains: Lemma \ref{th:equivalent_condition_deltareg} provides a sufficient condition for $\delta$-regularity;
Lemma \ref{th:lemma_smooth} shows that smooth domains are $\delta$-regular for sufficiently small values of $\delta$; 
Lemma \ref{th:geo_lemma} shows an interesting property of $\delta$-regular domains that is exploited below.
As currently stated, Lemma \ref{th:lemma_smooth} does not provide an explicit estimate of $\delta_0$; in  Appendix \ref{sec:proofs} we provide a more detailed expression for $\delta_0$. We anticipate that $\delta_0$ depends on the maximum principal curvature over $\partial U$ and also on a global property of the boundary.

\begin{lemma}
\label{th:equivalent_condition_deltareg}
Let $U$ be diffeomorphic to the unit ball.
Assume that the hyper-surfaces 
$\{x\pm \delta \mathbf{n}(x): x\in \partial U \}$ do not have self intersections. Then, $U$ is $\delta$-regular. 
\end{lemma}

\begin{lemma}
\label{th:lemma_smooth}
Let $U$ be a $C^2$ domain isomorphic to the unit ball. Then, there exists $\delta_0>0$ such that  $U$ is $\delta$-regular for any $\delta<\delta_0$. 
\end{lemma}

\begin{lemma}
\label{th:geo_lemma}
Let $U,V$ be Lipschitz domains   that  are isomorphic to the unit ball. Assume that
(i) $V$ is $\delta$-regular for some $\delta>0$;
(ii) $|U|_{(d)} > | {\rm Neigh}_{\delta}(\partial V) |_{(d)}$;
(iii) ${\rm dist}_{\rm bnd}(U, V) = 
\max_{{x} \in \partial U} {\rm dist} ({x}, \partial V ) \leq \delta$.
Then, 
${\rm dist}_{\rm H}(\partial U, \partial V) \leq \delta$ and 
${\rm dist}_{\rm H}(U, V) \leq \delta$.
\end{lemma}

\subsection{Equivalence of  \eqref{eq:hausdorff_distance} and \eqref{eq:quasi_hausdorff_distance}  for smooth domains}
\label{sec:justification_main}
Proposition \ref{th:geo_theorem} contains the result that is  relevant for our discussion.

\begin{proposition}
\label{th:geo_theorem}
Let  ${\Phi}: \Omega_{\rm box} \to \mathbb{R}^d$ be  a bijection in $\Omega_{\rm box}$ and let $\Omega,V \subset \Omega_{\rm box}$ be  
diffeomorphic to the unit ball.
Define $\epsilon: = \min_{{x} \in \overline{\Omega}}  {\rm det} (\nabla {\Phi})$, 
$\delta :={\rm dist}_{\rm bnd}(\Phi(\Omega); V )$. Assume that 
(i) $\epsilon | \Omega   |_{(d)} >   | {\rm Neigh}_{\delta}(\partial V) |_{(d)}$  and
(ii) $V$ is $\delta$-regular.
Then, we have 
${\rm dist}_{\rm H}( {\Phi}(\Omega), V ) \leq \delta$, ${\rm dist}_{\rm H}( \partial {\Phi}( \Omega), \partial V) \leq \delta$.
\end{proposition}

\begin{proof}
We define $U= {\Phi}(\Omega)$. Clearly, $U$ is isomorphic to the unit ball; furthermore, using the change-of-variable formula and the definition of $\epsilon$, we obtain
$$
|U|_{(d)}
=
\int_{\Omega} \;
 {\rm det} (\nabla {\Phi}) \, d{x}
\; \geq \;
\epsilon |\Omega|_{(d)}
>  | {\rm Neigh}_{\delta}(\partial V) |_{(d)}.
$$
Using the definition \eqref{eq:quasi_hausdorff_distance}  and 
 the identity $\partial \Phi(\Omega) = \Phi(\partial \Omega)$, 
 which is valid for any  bijection  $\Phi$, we also find
$$
\max_{{x} \in \partial U} {\rm dist} ({x}, \partial V  ) 
=
\max_{x \in \partial \Omega} {\rm dist} (  
{\Phi}({x}), \partial V  ) 
\leq \delta.
$$
Therefore,
$U,V$ satisfy the hypotheses of 
Lemma \ref{th:geo_lemma}: we conclude that  
${\rm dist}_{\rm H}(  U,   V ),  {\rm dist}_{\rm H}( \partial U,  \partial V ) \leq \delta$.
\end{proof}

Proposition  \ref{th:geo_theorem}  --- together with     Proposition  \ref{th:quasi_hausdorff_bound} ---
provides a rigorous justification of the geometry registration strategy proposed in   section \ref{sec:method}.
We observe that the value $\epsilon$ in 
   Proposition   \ref{th:geo_theorem} is weakly controlled by the exp-jac objective function \eqref{eq:exp_jac}. In practice, we expect $\epsilon$ to be significantly larger than the geometric error $\delta$ at the boundary: therefore, the condition $\epsilon|\Omega|_{(d)} \geq  | {\rm Neigh}_{\delta}(\partial V) |_{(d)}$ should  in practice be easy to enforce, even for slender domains with small $\frac{|V|_{(d)}}{|\partial V   |_{(d-1)}}$. 

On the other hand, we cannot establish a rigorous  connection between the exp-mesh objective function \eqref{eq:exp_mesh}
and the quantity $\epsilon$ in 
Proposition   \ref{th:geo_theorem}.
In this respect, we recall  that 
in the discretize-then-map framework (cf.  Remark \ref{remark:MtD_vs_DtM})
the mapping to be considered in   Proposition  \ref{th:geo_theorem} for the mesh $\mathcal{T}_{\rm hf}$ is 
$\Phi_{\rm hf}$  (cf. \eqref{eq:mapping_DtM}):
for purely isotropic P1 FE discretizations, it is easy to verify that
$\nabla \Phi_{\rm hf}|_{\texttt{D}_k} = 
\frac{| \texttt{D}_{\Phi,k}   |}{|\texttt{D}_k|}
\mathbbm{1}$, while the ratio $q_{\Phi,k}$ in \eqref{eq:exp_mesh} is equal to one.
We conclude that the solution to the registration  problem with objective given by \eqref{eq:exp_mesh} might lead to extremely small elements 
(see Figure \ref{fig:vis_mesh_three} in the numerical example of section \ref{sec:three_point_registration}),
which might prevent the application of   
Proposition \ref{th:geo_theorem}. If large deformations are expected in the proximity of the boundary, it might thus be necessary to combine  \eqref{eq:exp_jac} and \eqref{eq:exp_mesh}.

\subsection{Strict inequality for  Lipschitz domains}
\label{sec:Lipschitz_analysis}
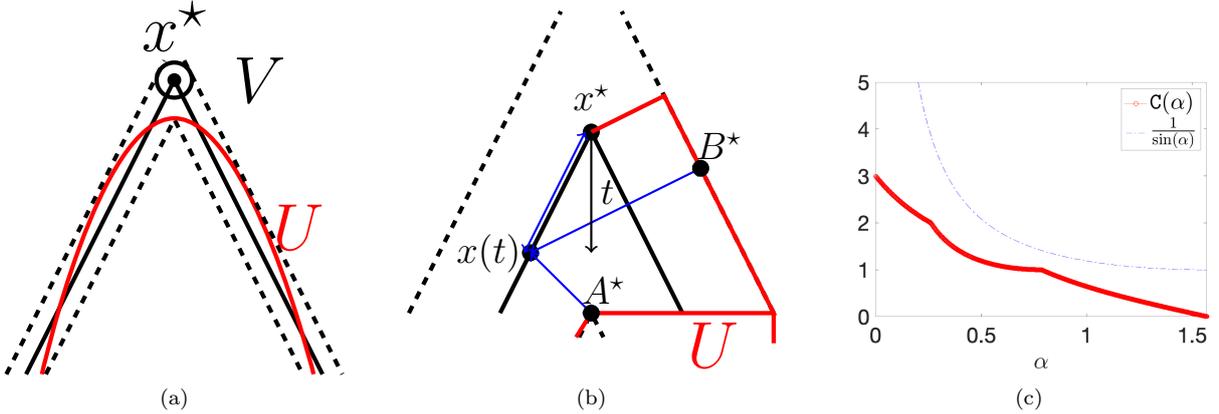
\begin{figure}[h!]
\centering
 \subfloat[] {
\begin{tikzpicture}[scale=1.3]

\draw [black,ultra thick] (0,0) -- (1.5,3) -- (3,0);

\draw [black,dashed,ultra thick] (0.2,0) -- (1.5,2.6) -- (2.8,0);

\draw [ultra thick,red] plot [smooth, tension=1] coordinates { (0.16,0) (1.5,2.61)  (2.91,0)};

\draw [black,ultra thick,domain=0:360,samples=300] plot ({1.5+0.1789*cos(\x)}, {3+0.1789*sin(\x)});

\coordinate [label={right:  {\Huge {${V}$}}}] (E) at (2,3) ;
\coordinate [label={right:  {\Huge {${\color{red}  U}$}}}] (E) at (2.4,1.5) ;

\fill (1.5,3)  circle[radius=2pt];
\coordinate [label={above:  {\Huge {${x}^{\star}$}}}] (E) at (1.5,3.2) ;

\draw [black,dashed,ultra thick] (-0.2,0) -- (1.4,3.2);
\draw [black,dashed,ultra thick]   (1.6,3.2) -- (3.2,0);
 
 \end{tikzpicture}
} 
 ~~~~ 
  \subfloat[] {
\begin{tikzpicture}[scale=0.8]

\fill (1.5,3)  circle[radius=4pt];
\coordinate [label={above:  {\Large {${x}^{\star}$}}}] (E) at (1.5,3.2) ;

\draw [black,ultra thick] (0,0) -- (1.5,3) -- (3,0);
\draw [black,dashed,ultra thick] (-1.5,0) -- (1,5) ;
\draw [black,dashed,ultra thick] (4.5,0) -- (2,5) ;
\draw [black,dashed,ultra thick] (1.3,-0.4) -- (1.5,0) -- (1.7,-0.4) ;

\draw [red,ultra thick] (1.5,3) -- (2.7,3.6) -- (4.5,0) -- (1.5,0) --  (1.25,-0.4) ;

\draw [ultra thick,red] (4.5,0) --  (4.5,-0.5);

\fill (0.5,1)  circle[radius=4pt];
\coordinate [label={left:  {\Large {${x}(t)$}}}] (E) at (0.5,1) ;
\coordinate [label={below:  {\Huge {${\color{red}  U}$}}}] (E) at (3.5,0) ;

\draw [<->,blue,thick] (0.5,1) -- (3.3,2.4);
\draw [<->,blue,thick] (0.5,1) -- (1.5,0);
\draw [<->,blue,thick] (0.4,1) -- (1.4,3);
\draw [<->,black,thick] (1.5,3) -- (1.5,1);
\coordinate [label={right:  {\Large {$t$}}}] (E) at (1.5,2) ;

\coordinate [label={above:  {\Large {${A}^{\star}$}}}] (E) at (1.7,0) ;
\fill (1.5,0)  circle[radius=4pt];

\coordinate [label={above:  {\Large {${B}^{\star}$}}}] (E) at (3.6,2.4) ;
\fill (3.3,2.4)  circle[radius=4pt];

 \end{tikzpicture}
  }
 ~~
 \subfloat[ ] 
{  \includegraphics[width=0.33\textwidth]
 {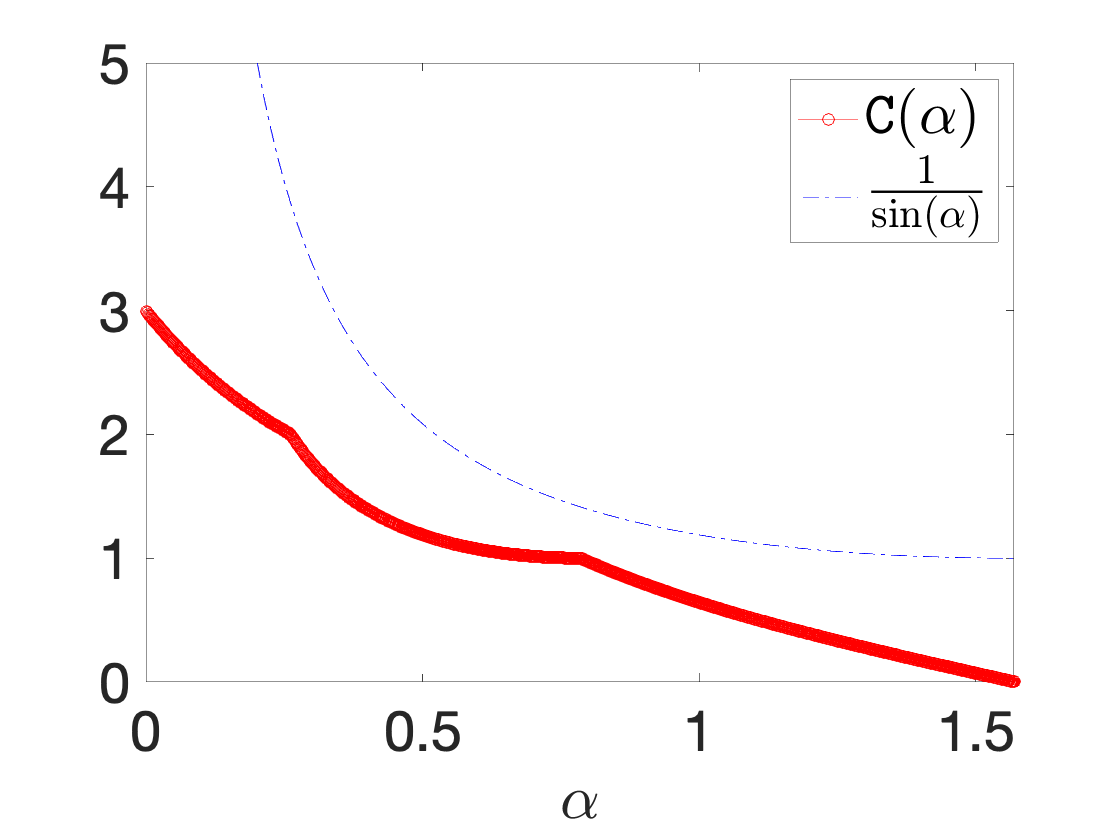}}

 \caption{analysis for Lipschitz domains. 
 (a) example of domains $U,V$ isomorphic to $\mathcal{B}_1(0)$ with $|U|_{(2)} > 2\delta | \partial V|_{(1)}$ and 
 $\partial U \subset {\rm Neigh}_{\delta}(\partial V)$
 such that 
${\rm dist}_{\rm H}(\partial U, \partial V) = \frac{\delta}{\sin(\alpha)}$.
  (b) domain $U$ satisfying (i)-(iv) such that ${\rm dist}_{\rm H}(\partial U, \partial V) = 
\texttt{C}(\alpha) \delta$.
 (c) behavior of $\frac{1}{\sin(\alpha)}$ and $\texttt{C}(\alpha)$ with respect to $\alpha \in [0,\pi/2]$. 
}
\label{fig:polyB_lipschitz}
\end{figure}

The   analysis of the previous section cannot be readily extended to Lipschitz domains with corners. To investigate the problem, we shall   consider the domain depicted in   Figure \ref{fig:polyB_lipschitz}: we are here interested in the neighborhood of the vertex 
${x}^{\star}$;   we assume that $V$ is of class $C^{2}$ elsewhere.
We denote by $2 \alpha \in [0,\pi]$ the angle associated with the corner at ${x}^{\star}$; the case 
 $2 \alpha \in [\pi,2 \pi]$ is analogous.

We first consider the domains $U$  and $V$ in   Figure \ref{fig:polyB_lipschitz}(a): the dashed lines denote the boundary of ${\rm Neigh}_{\delta}(\partial V)$. 
Clearly,  $U$ belongs to ${\rm Neigh}_{\delta}(\partial U)$;  the point ${x}^{\star}$ satisfies ${\rm dist}({x}^{\star}, \partial  U) = \frac{\delta}{\sin(\alpha)} =: r_{\delta}> \delta$. 
The example shows that, for a domain $V$ with a corner of angle $2 \alpha$ we might construct a domain $U$ such that
$U$ is isomorphic to $\mathcal{B}_1({0})$,
$|U|_{(d)} > 2 \delta |\partial V |_{(d-1)} $, 
$\partial U \subset {\rm Neigh}_{\delta}(\partial V)$, and
$\max_{x \in \partial V} {\rm dist}(x, U) = r_{\delta}> \delta$.

We seek an upper bound for 
$\max_{{x} \in \partial V} {\rm dist}({x}, U) $ among all domains $U$ satisfying
\textbf{(i)}
$U$ is isomorphic to the unit ball $\mathcal{B}_1({0})$,
\textbf{(ii)}
$|U|_{(d)} > 2 \delta |\partial V |_{(d-1)}$,
\textbf{(iii)}
$\partial U \subset {\rm Neigh}_{\delta}(\partial V)$, and 
\textbf{(iv)}
${x}^{\star} \in \partial U$.
Exploiting the same argument as in   Lemma \ref{th:lemma_smooth} (cf. Appendix \ref{sec:proofs}), we find that there exists $\delta_0>0$ such that the curves 
${\gamma}_{\delta}^{\pm} =
{\gamma} \pm \delta \mathbf{n}$ do not have self-intersections in $\partial V \setminus \mathcal{B}_{r_{\delta}}( {x}^{\star}   )$ for all $\delta \leq \delta_0$. Therefore, we have
$$
\max_{{x} \in \partial V  \setminus \mathcal{B}_{r_{\delta}}( {x}^{\star}   )} {\rm dist}({x}, \partial U) \leq \delta,
$$
for any $U$ that satisfies (i)-(iv).
On the other hand, in the neighborhood of ${x}^{\star}$ , we can show that the domain $U$ depicted in   Figure \ref{fig:polyB_lipschitz}(b) maximizes 
$\max_{{x} \in \partial V \cap \mathcal{B}_{r_{\delta}}( {x}^{\star}   )} {\rm dist}({x}, \partial U)$. Note that   $U$ is the limit of domains isomorphic to the unit ball; therefore, computation of  the distance of $\partial V$ from $U$ provides an upper bound for 
$\max_{{x} \in \partial V \setminus \mathcal{B}_{r_{\delta}}( {x}^{\star}   )} {\rm dist}({x}, \partial U)$ among Lipschitz domains $U$ that  satisfy  (i)-(iv).
 Exploiting well-known trigonometric identities, we find that the point
 ${x}(t)$ in  Figure \ref{fig:polyB_lipschitz}(b) satisfies
 $$
 \| {x}(t) - {A}^{\star} \|_2
 =
 \sqrt{
(r_{\delta} - t)^2 + (\tan(\alpha) t )^2
 }
 =
 \delta
  \sqrt{
 (1/\sin(\alpha)  - t')^2 + (\tan(\alpha) t' )^2
 },
 \quad
 t'=\frac{t}{\delta}, 
 $$
 and 
 $$
 \| {x}(t) - {x}^{\star} \|_2
 =
 \delta \frac{t'}{\cos(\alpha)},
\qquad
 \| {x}(t) - {B}^{\star} \|_2
 =
\left(
1+2 \sin(\alpha) t'
\right) \delta.
 $$
In conclusion, we obtain
\begin{subequations}
\label{eq:crazy_identity_lipschitz}
\begin{equation}
\max_{{x} \in \partial U} \;
{\rm dist} ({x}, \partial V)
\; =  
\max_{  t \in \left[ 0, \frac{\delta}{\sin(\alpha)} - \delta \sin(\alpha)  \right]  } \;
{\rm dist} ({x}(t), \partial U)
 = 
 \texttt{C}(\alpha) \; \delta,
\end{equation}
where
\begin{equation}
\texttt{C}(\alpha) = 
\max_{t' \in \left[ 0, \frac{1}{\sin(\alpha)} -   \sin(\alpha)  \right]  } 
\min \left\{
  \sqrt{
 (1/\sin(\alpha)  - t')^2 + (\tan(\alpha) t' )^2
 },
\;
 \frac{t'}{\cos(\alpha)},  \;
 1+2 \sin(\alpha) t'
\right\}.
\end{equation}
\end{subequations}

In   Figure \ref{fig:polyB_lipschitz}(c), we show the behavior of $\texttt{C}(\alpha)$ with respect to $\alpha$ and we compare it with $\frac{1}{\sin(\alpha)}$. Note that 
$\texttt{C}(\alpha) \leq  \min \{3, \frac{1}{\sin(\alpha)} \}$ for all $\alpha \in [0,\pi/2]$; $\texttt{C}$ is monotonic decreasing with $\alpha$; and 
$\texttt{C}(\alpha=\pi/4) = 1$.
The  assumption 
${x}^{\star} \in \partial U$ leads to a reduction of the worst-case scenario for $\max_{{x} \in \partial V} \;
{\rm dist} ({x}, \partial U)$ (and thus for ${\rm dist}_{\rm H}(\partial V, \partial U)$) but does not guarantee that 
$\max_{{x} \in \partial V} \;
{\rm dist} ({x}, \partial U) = \delta$ asymptotically. We note, however, that, while the domain  $U$ depicted in  Figure \ref{fig:polyB_lipschitz}(a) is regular, the domain  in  Figure \ref{fig:polyB_lipschitz}(b) is not Lipschitz and thus the upper bound in \eqref{eq:crazy_identity_lipschitz} cannot be achieved. 

The analysis of this section shows that the approximation of slender bodies and/or of bodies with cusps and corners is challenging for the  optimization-based registration procedure discussed in this paper.
However, if we are able to to ensure appropriate approximation --- via interpolation --- of the target shape at corners and cusps, we can recover near-optimal bounds for the geometry error.

\section{Numerical results}
\label{sec:numerics}

We present several numerical experiments for a number of two-dimensional model problems. In order to assess the performance of a given mapping $\Phi$, we report the behavior of the minimum value $q_{\rm min}$ of the radius ratio\footnote{The radius ratio is the ratio between the radius of the circle inscribed in the triangle and the radius of the circle circumscribed around the triangle.} over the elements of the deformed target mesh $\mathcal{T}_{\rm hf}$.
We also report the minimum value of the Jacobian determinant over the domain $\overline{\Omega}$, $J_{\rm min} =\min_{x \in \overline{\Omega}} J(\Phi)(x)$.  All numerical simulations are performed in Matlab 2020b on a commodity laptop.

\subsection{Three-point registration}
\label{sec:three_point_registration}
Given $\Omega=(0,1)^2$, we consider the problem of deforming the points
$x_1=(1/2,1/2)$, $x_2=(1/4,1/4)$, $x_3=(3/4,1/4)$ into the points
$y_1=(1/4,3/4)$, $y_2=(1/16,1/16)$, $y_3=(1/2,1/4)$.
The aim of this test is to illustrate the impact of the objective function on performance.
Towards this end, we consider the Morozov-regularized formulation \eqref{eq:morozov_registration} with objective function equal to $\mathfrak{f}^{\rm obj} = \mathfrak{P}$ (dubbed ``H2'' below),
$\mathfrak{f}^{\rm obj} = \mathfrak{f}_{\rm jac} + \mathfrak{P}$ (dubbed ``exp-jac''),
$\mathfrak{f}^{\rm obj} = \mathfrak{f}_{\rm msh} + \mathfrak{P}$ (dubbed ``exp-mesh'') --- we refer to  section \ref{sec:rectangular_domain} for the definitions. 
We also investigated other objective functions
including the ones associated with the linear elasticity strain energy
\eqref{eq:isotropic_linear_strain} and with the neohookean strain energy 
\eqref{eq:neohookean_strain}: 
we provide a representative test for these objective functions in  Appendix \ref{sec:further_tests}.

Figure \ref{fig:vis_mesh_three}(a) shows the reference and deformed points and the mesh $\mathcal{T}_{\rm hf}$ that is used to assess performance. We resort to the Matlab function \texttt{fmincon} based on the interior-point method; in all our tests, we consider the initial condition $\Phi = \texttt{id}$.
Figures \ref{fig:vis_mesh_three}(b) and (c) show the deformed mesh obtained using the exp-jac objective and the exp-mesh objective, respectively. We consider $\delta=10^{-6}$ in \eqref{eq:morozov_registration}, and we set 
 $\epsilon=0.1$ for the exp-jac objective,  and $\kappa_{\rm msh}=10$
for the exp-mesh objective. Note that the two deformed meshes significantly differ from each other in the proximity of the origin.

\begin{figure}[h!]
\centering
 \subfloat[ ] 
{  \includegraphics[width=0.33\textwidth]
 {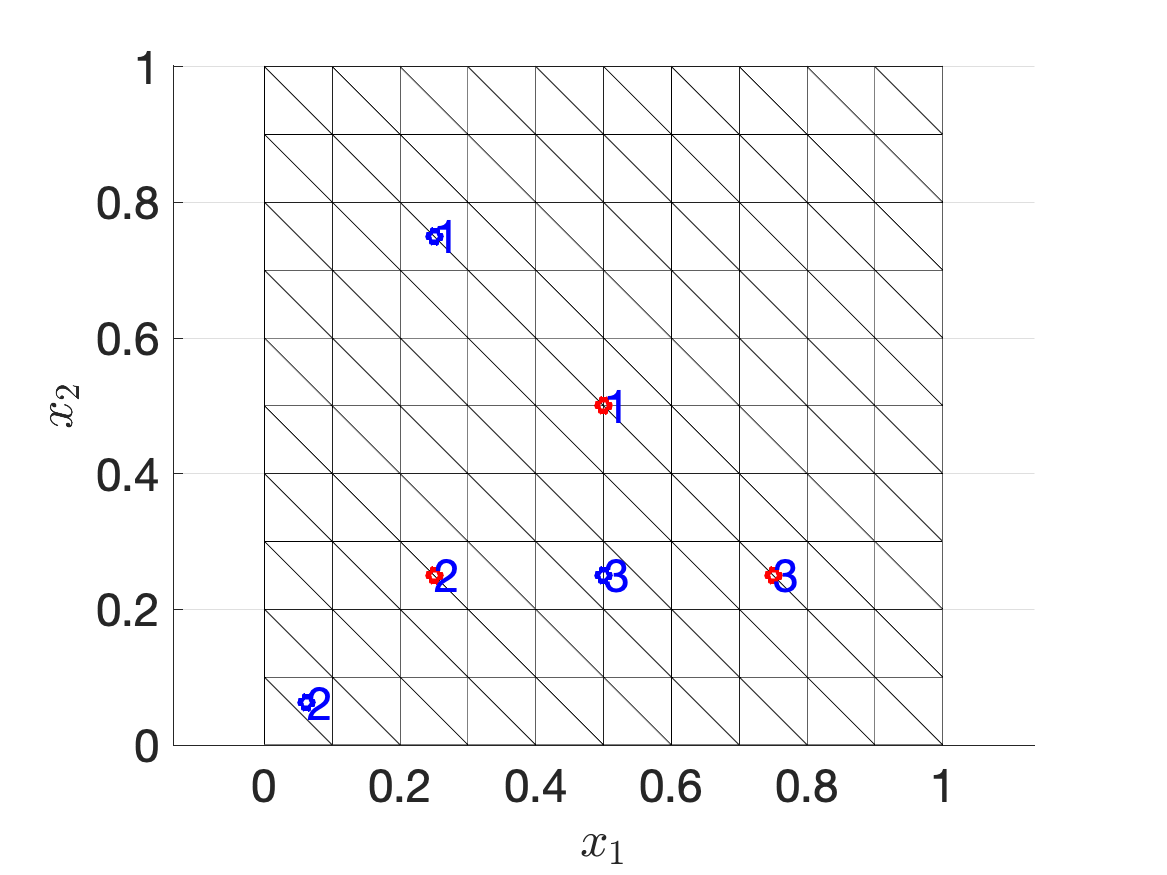}}
   ~~
 \subfloat[exp-mesh] 
{  \includegraphics[width=0.33\textwidth]
 {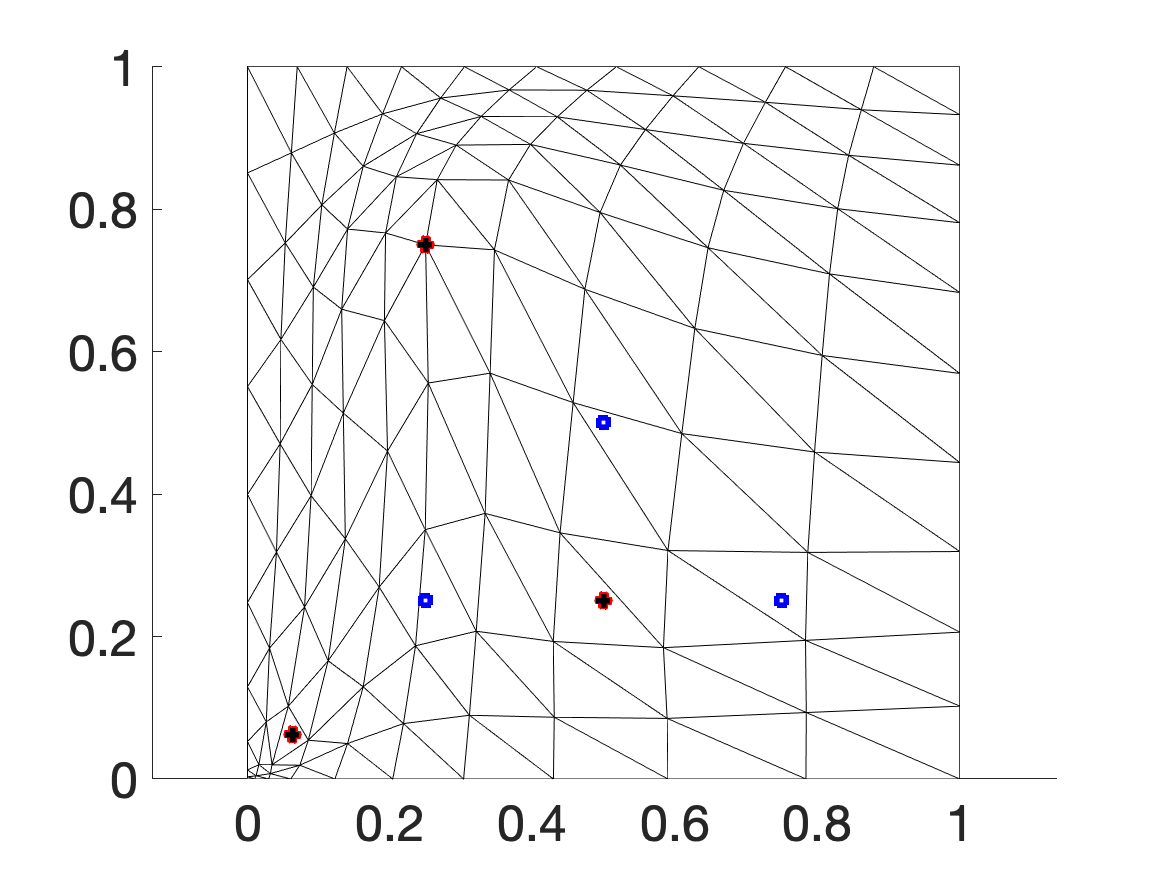}}
~~
 \subfloat[exp-jac ] 
{  \includegraphics[width=0.33\textwidth]
 {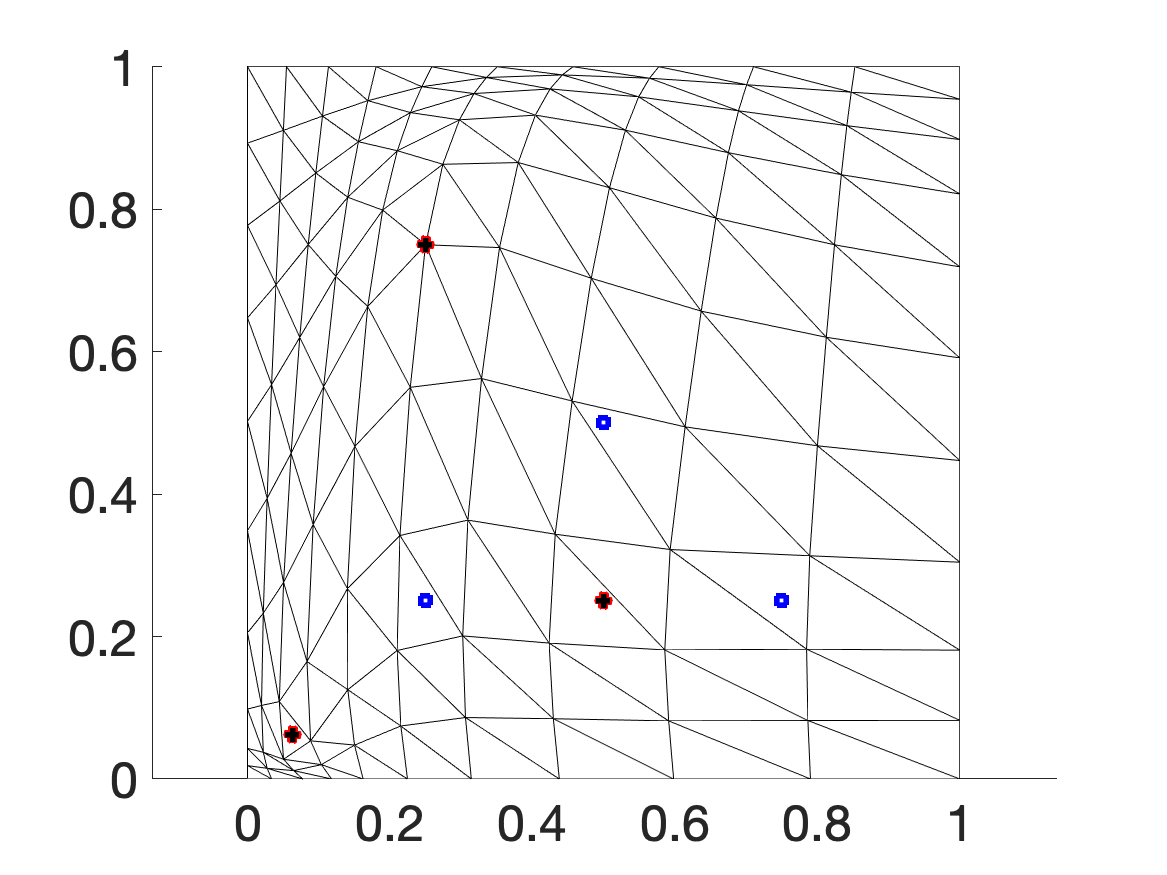}}
  
 \subfloat[exp-mesh] 
{  \includegraphics[width=0.33\textwidth]
 {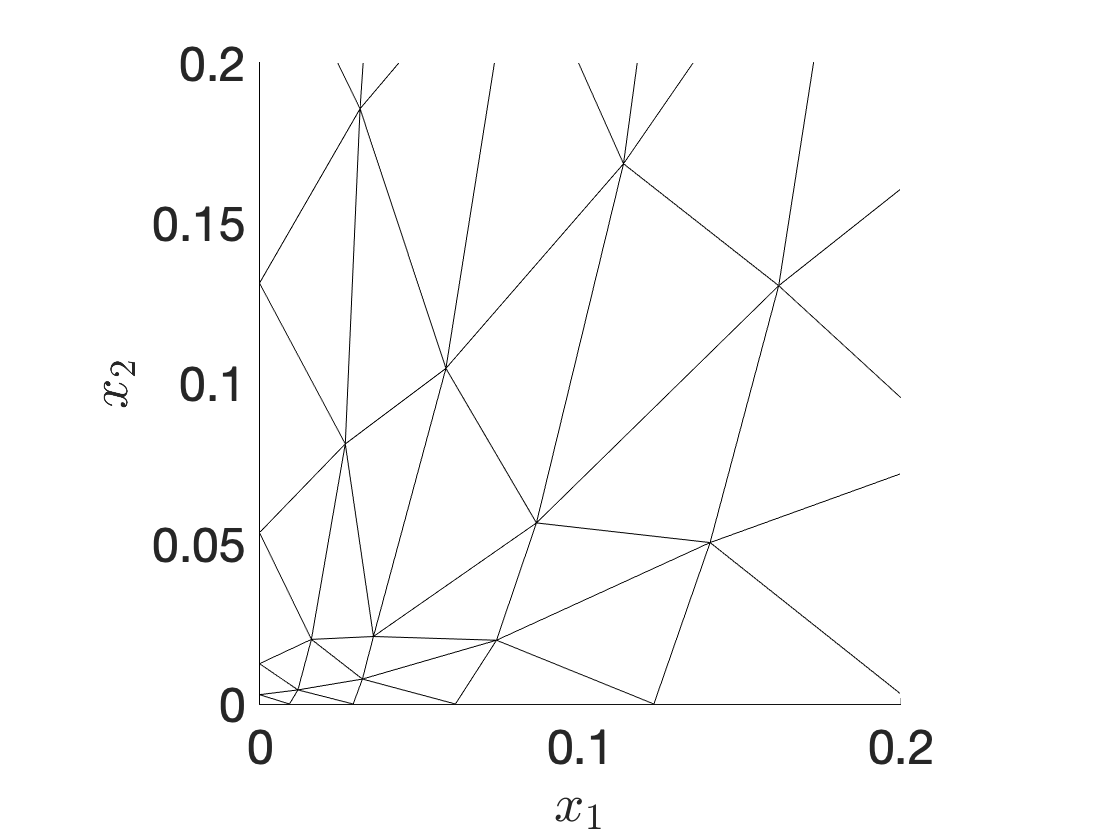}}
~~
 \subfloat[exp-jac ] 
{  \includegraphics[width=0.33\textwidth]
 {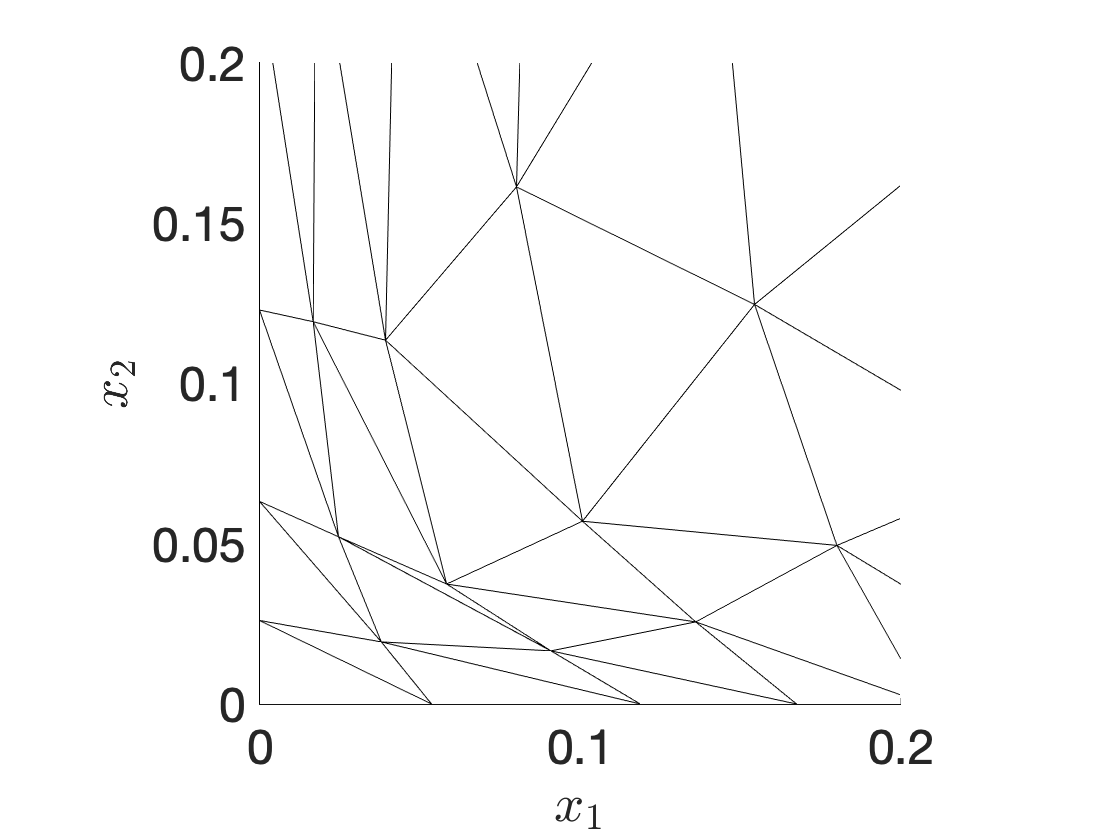}}  
  
 \caption{
 three-point registration problem; visualization of the reference mesh and of the deformed mesh for two different registration algorithms ($\delta=10^{-6}$, $\epsilon=0.1$, $\kappa_{\rm msh}=10$).
Figures (d)-(e) show a zoom of the deformed meshes in the proximity of the origin.
}
 \label{fig:vis_mesh_three}
  \end{figure}

\begin{figure}[h!]
\centering
\subfloat[$H^2$]{
\begin{tikzpicture}[scale=0.6]
\begin{loglogaxis}[
xmode=linear,
ymode=log,
grid=both,
minor grid style={gray!25},
major grid style={gray!25},
xlabel = {\LARGE {$t$} },
  ylabel = {\LARGE {$q_{\rm min}$}},
  line width=1.2pt,
  mark size=3.0pt,
xmin=0,   xmax=1,
 ymin=0.01,   ymax=1,
 ]

\addplot[line width=1.pt,color=red,mark=square]  table {data/three/H2/qmin_delta1m6.dat}; 
\label{three:delta1m6}
      
\addplot[line width=1.pt,color=blue,mark=triangle*] table {data/three/H2/qmin_delta1m4.dat};
\label{three:delta1m4}
       
\addplot[line width=1.pt,color=black,mark=pentagon] table {data/three/H2/qmin_delta1m2.dat};
\label{three:delta1m2}
   
\end{loglogaxis}
\end{tikzpicture}

\bgroup
\sbox0{\ref{data}}%
\pgfmathparse{\ht0/1ex}%
\xdef\refsize{\pgfmathresult ex}%
\egroup
}
\subfloat[exp-jac, $\epsilon=0.1$]{
\begin{tikzpicture}[scale=0.6]
\begin{loglogaxis}[
xmode=linear,
ymode=log,
grid=both,
minor grid style={gray!25},
major grid style={gray!25},
xlabel = {\LARGE {$t$} },
  line width=1.2pt,
  mark size=3.0pt,
xmin=0,   xmax=1,
 ymin=0.01,   ymax=1,
 ]

\addplot[line width=1.pt,color=red,mark=square]  table {data/three/expjac_1/qmin_delta1m6.dat}; 
       
\addplot[line width=1.pt,color=blue,mark=triangle*] table {data/three/expjac_1/qmin_delta1m4.dat};

\addplot[line width=1.pt,color=black,mark=pentagon] table {data/three/expjac_1/qmin_delta1m2.dat};

\end{loglogaxis}
\end{tikzpicture}
}
\subfloat[exp-mesh, $\kappa=10$]{
\begin{tikzpicture}[scale=0.6]
\begin{loglogaxis}[
xmode=linear,
ymode=log,
grid=both,
minor grid style={gray!25},
major grid style={gray!25},
xlabel = {\LARGE {$t$} },
  line width=1.2pt,
  mark size=3.0pt,
xmin=0,   xmax=1,
 ymin=0.01,   ymax=1,
 ]

\addplot[line width=1.pt,color=red,mark=square]  table {data/three/expmesh/qmin_delta1m6.dat}; 
       
\addplot[line width=1.pt,color=blue,mark=triangle*] table {data/three/expmesh/qmin_delta1m4.dat};

\addplot[line width=1.pt,color=black,mark=pentagon] table {data/three/expmesh/qmin_delta1m2.dat};
   
\end{loglogaxis}
\end{tikzpicture}
}

\subfloat[$H^2$]{
\begin{tikzpicture}[scale=0.6]
\begin{loglogaxis}[
xmode=linear,
ymode=log,
grid=both,
minor grid style={gray!25},
major grid style={gray!25},
xlabel = {\LARGE {$t$} },
  ylabel = {\LARGE {$J_{\rm min}$}},
  line width=1.2pt,
  mark size=3.0pt,
xmin=0,   xmax=1,
 ymin=0.001,   ymax=1,
 ]

\addplot[line width=1.pt,color=red,mark=square]  table {data/three/H2/Jmin_delta1m6.dat}; 
\label{three:delta1m6}
      
\addplot[line width=1.pt,color=blue,mark=triangle*] table {data/three/H2/Jmin_delta1m4.dat};
\label{three:delta1m4}
       
\addplot[line width=1.pt,color=black,mark=pentagon] table {data/three/H2/Jmin_delta1m2.dat};
   
\end{loglogaxis}
\end{tikzpicture}
}
\subfloat[exp-jac, $\epsilon=0.1$]{
\begin{tikzpicture}[scale=0.6]
\begin{loglogaxis}[
xmode=linear,
ymode=log,
grid=both,
minor grid style={gray!25},
major grid style={gray!25},
xlabel = {\LARGE {$t$} },
  line width=1.2pt,
  mark size=3.0pt,
xmin=0,   xmax=1,
 ymin=0.001,   ymax=1,
 ]

\addplot[line width=1.pt,color=red,mark=square]  table {data/three/expjac_1/Jmin_delta1m6.dat}; 
       
\addplot[line width=1.pt,color=blue,mark=triangle*] table {data/three/expjac_1/Jmin_delta1m4.dat};

\addplot[line width=1.pt,color=black,mark=pentagon] table {data/three/expjac_1/Jmin_delta1m2.dat};

\end{loglogaxis}
\end{tikzpicture}
}
\subfloat[exp-mesh, $\kappa=10$]{
\begin{tikzpicture}[scale=0.6]
\begin{loglogaxis}[
xmode=linear,
ymode=log,
grid=both,
minor grid style={gray!25},
major grid style={gray!25},
xlabel = {\LARGE {$t$} },
  line width=1.2pt,
  mark size=3.0pt,
xmin=0,   xmax=1,
 ymin=0.001,   ymax=1,
 ]

\addplot[line width=1.pt,color=red,mark=square]  table {data/three/expmesh/Jmin_delta1m6.dat}; 
       
\addplot[line width=1.pt,color=blue,mark=triangle*] table {data/three/expmesh/Jmin_delta1m4.dat};

\addplot[line width=1.pt,color=black,mark=pentagon] table {data/three/expmesh/Jmin_delta1m2.dat};
   
\end{loglogaxis}
\end{tikzpicture}
}
\caption[Caption in ToC]{three-point registration problem.
Performance of the Morozov-regularized statement \eqref{eq:morozov_registration} for the  objective functions H2, exp-jac and exp-mesh, for three choices of $\delta$ and for the deformation points $\{y_i(t) = (1-t) x_i + t y_i \}_{i=1}^3$ for $t\in [0,1]$. 
$\delta=10^{-6}$
\tikzref{three:delta1m6};
$\delta=10^{-4}$
\tikzref{three:delta1m4},
$\delta=10^{-2}$
\tikzref{three:delta1m2}.}
\label{fig:three_qmin_Jmin}
\end{figure}
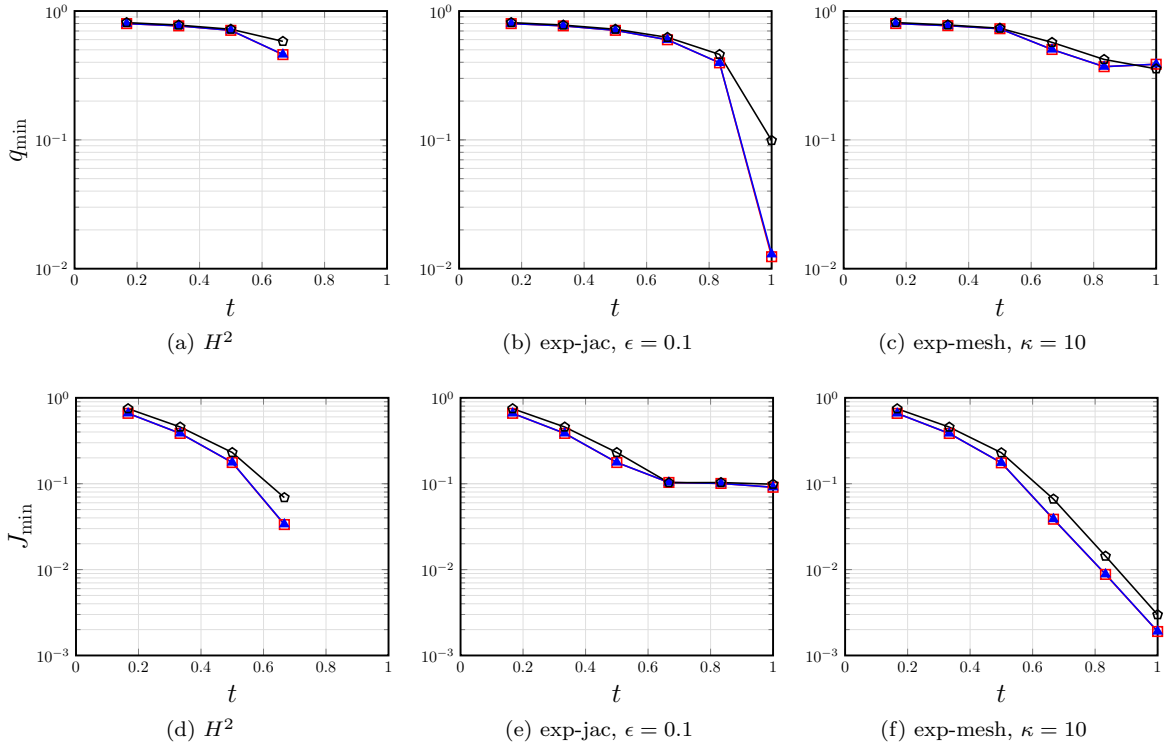

 Figure \ref{fig:three_qmin_Jmin} shows the performance of the Morozov-regularized statement \eqref{eq:morozov_registration} for the  objective functions H2, exp-jac and exp-mesh, for three choices of $\delta$ and for the deformation points $\{y_i(t) = (1-t) x_i + t y_i \}_{i=1}^3$ for $t\in [0,1]$. The quadratic objective function $H^2$ fails to deliver a proper mesh for $t>0.6$ --- i.e., the deformed mesh contains inverted elements. On the other hand, the two nonlinear approaches are able to deliver bijective maps that are also bijective with respect to the target mesh $\mathcal{T}_{\rm hf}$.

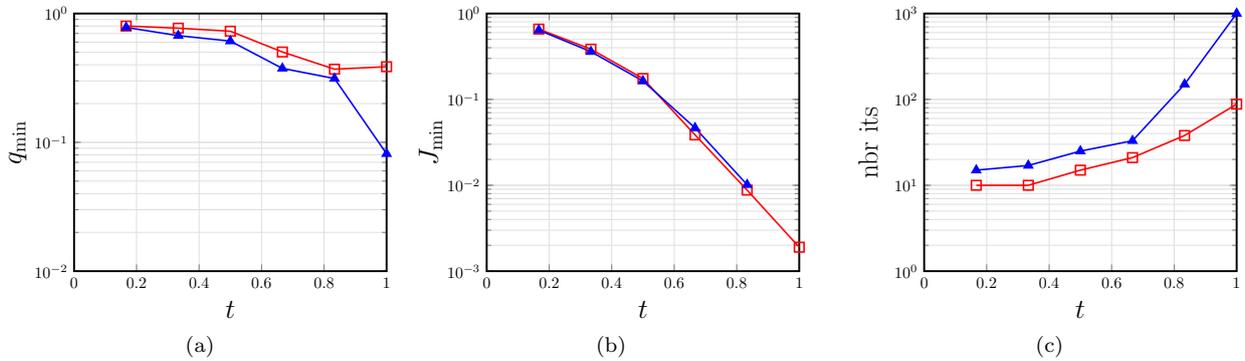
\begin{figure}[h!]
\centering
\subfloat[]{
\begin{tikzpicture}[scale=0.6]
\begin{loglogaxis}[
xmode=linear,
ymode=log,
grid=both,
minor grid style={gray!25},
major grid style={gray!25},
xlabel = {\LARGE {$t$} },
  ylabel = {\LARGE {$q_{\rm min}$}},
  line width=1.2pt,
  mark size=3.0pt,
xmin=0,   xmax=1,
 ymin=0.01,   ymax=1,
 ]

\addplot[line width=1.pt,color=red,mark=square]  table {data/three/pot_vs_full/qmin_full.dat}; 
      
\addplot[line width=1.pt,color=blue,mark=triangle*] table {data/three/pot_vs_full/qmin_pot.dat};
         
\end{loglogaxis}
\end{tikzpicture}
}
\subfloat[]{
\begin{tikzpicture}[scale=0.6]
\begin{loglogaxis}[
xmode=linear,
ymode=log,
grid=both,
minor grid style={gray!25},
major grid style={gray!25},
xlabel = {\LARGE {$t$} },
 ylabel = {\LARGE {$J_{\rm min}$}},
  line width=1.2pt,
  mark size=3.0pt,
xmin=0,   xmax=1,
 ymin=0.001,   ymax=1,
 ]

\addplot[line width=1.pt,color=red,mark=square]  table {data/three/pot_vs_full/Jmin_full.dat}; 
       
\addplot[line width=1.pt,color=blue,mark=triangle*] table {data/three/pot_vs_full/Jmin_pot.dat};

\end{loglogaxis}
\end{tikzpicture}
}
~~
\subfloat[]{
\begin{tikzpicture}[scale=0.6]
\begin{loglogaxis}[
xmode=linear,
ymode=log,
grid=both,
minor grid style={gray!25},
major grid style={gray!25},
xlabel = {\LARGE {$t$} },
 ylabel = {\LARGE {nbr its}},
  line width=1.2pt,
  mark size=3.0pt,
xmin=0,   xmax=1,
 ymin=1,   ymax=1000,
 ]

\addplot[line width=1.pt,color=red,mark=square]  table {data/three/pot_vs_full/nbrit_full.dat}; 
       
\addplot[line width=1.pt,color=blue,mark=triangle*] table {data/three/pot_vs_full/nbrit_pot.dat};

\end{loglogaxis}
\end{tikzpicture}
}
 
\caption[Caption in ToC]{three-point registration problem.
Comparison of performance for ``full'' space and ``potential'' space.
Objective function exp-mesh, $\kappa_{\rm msh}=10$, $\delta=10^{-6}$ (
full basis
\tikzref{three:delta1m6};
potential 
\tikzref{three:delta1m4}).
}
\label{fig:three_potvsfull}
\end{figure}

Figure \ref{fig:three_potvsfull} compares the performance of registration based on the ``full'' polynomial space (cf. \eqref{eq:polynomial_space}) with the performance of registration based on the ``potential space''
$\mathcal{W}_{\rm hf}^{\rm pot} = \{\texttt{id} + \nabla \phi : \phi \in \mathbb{Q}_{n_{\rm lp}'}, \partial_n \phi |_{\partial \Omega} = 0 \}$. We recall that the condition 
${\rm det} (\nabla \Phi) >0$ for $\Phi = \texttt{id} +\nabla \phi \in \mathcal{W}_{\rm hf}^{\rm pot}$ is equivalent to the convexity of 
$\phi$. The interest for this particular choice of the search space is that there exist several effective convexification procedures, 
which have been recently developed for optimal transport problems
(see,  e.g., \cite{jacobs2020fast}),
 that might pave  the way for effective registration strategies.
To ensure the fairness of the test, we choose two different polynomial degrees in order to have trial spaces of comparable dimension --- $M=1150$ for the ``full'' space and $M=1155$ for the ``potential'' space.
We observe that the potential approximation significantly deteriorates the performance of registration for large deformations.

\begin{figure}[h!]
\centering
\subfloat[]{
\begin{tikzpicture}[scale=0.6]
\begin{loglogaxis}[
xmode=linear,
ymode=log,
grid=both,
minor grid style={gray!25},
major grid style={gray!25},
xlabel = {\LARGE {$n_{\rm lp}$} },
  ylabel = {\LARGE {$q_{\rm min}$}},
  line width=1.2pt,
  mark size=3.0pt,
 ymin=0.01,   ymax=1,
 ]

\addplot[line width=1.pt,color=red,mark=square]  table {data/three/Nconv/qmin.dat}; 
      
         
\end{loglogaxis}
\end{tikzpicture}
}
\subfloat[]{
\begin{tikzpicture}[scale=0.6]
\begin{loglogaxis}[
xmode=linear,
ymode=log,
grid=both,
minor grid style={gray!25},
major grid style={gray!25},
xlabel = {\LARGE {$n_{\rm lp}$} },
 ylabel = {\LARGE {$J_{\rm min}$}},
  line width=1.2pt,
  mark size=3.0pt,
 ymin=0.001,   ymax=1,
 ]

\addplot[line width=1.pt,color=red,mark=square]  table {data/three/Nconv/Jmin.dat}; 
       

\end{loglogaxis}
\end{tikzpicture}
}
 
\caption[Caption in ToC]{three-point registration problem; sensitivity with respect to the polynomial order $n_{\rm lp}$ for $t=1$.  Objective function: exp-mesh; $\delta=10^{-6}$, $\kappa_{\rm msh}=10$.
}
\label{fig:three_Nconv}
\end{figure}
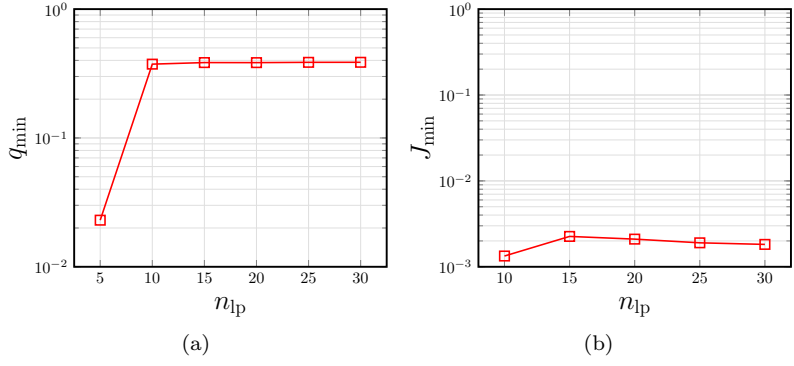

Figure \ref{fig:three_Nconv} shows the sensitivity of the algorithm with respect to the choice of the polynomial order $n_{\rm lp}$ in \eqref{eq:polynomial_space}. We here consider the exp-mesh objective function with $\kappa_{\rm msh}=10$. We observe that for this particular test case results are nearly independent of $n_{\rm lp}$ for $n_{\rm lp} \gtrsim 10$.

\subsection{Large deformation mesh morphing}
\label{sec:square_square_deformation}
\begin{figure}[h!]
\centering
 \subfloat[$\theta=0^o$] 
{  \includegraphics[width=0.33\textwidth]
 {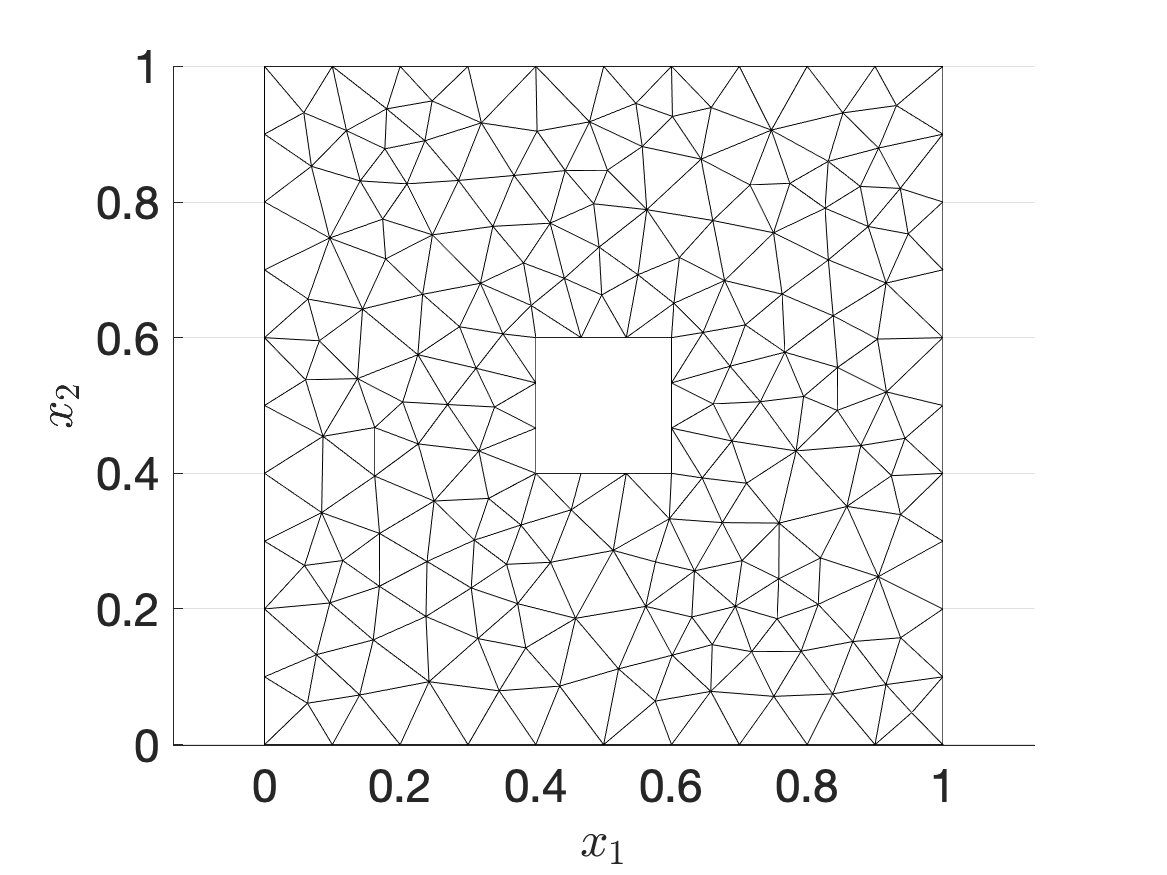}}
   ~~
 \subfloat[$\theta=60^o$] 
{  \includegraphics[width=0.33\textwidth]
 {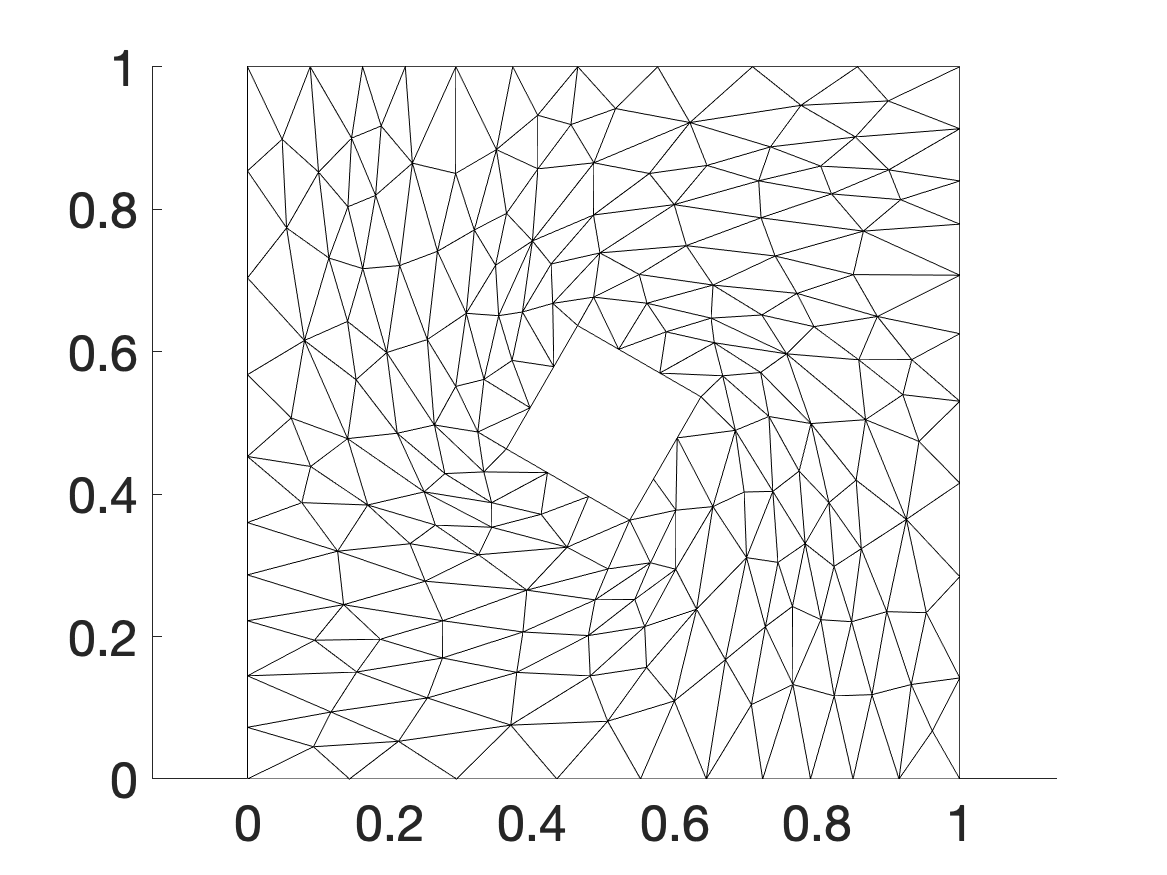}}
~~
 \subfloat[$\theta=120^o$] 
{  \includegraphics[width=0.33\textwidth]
 {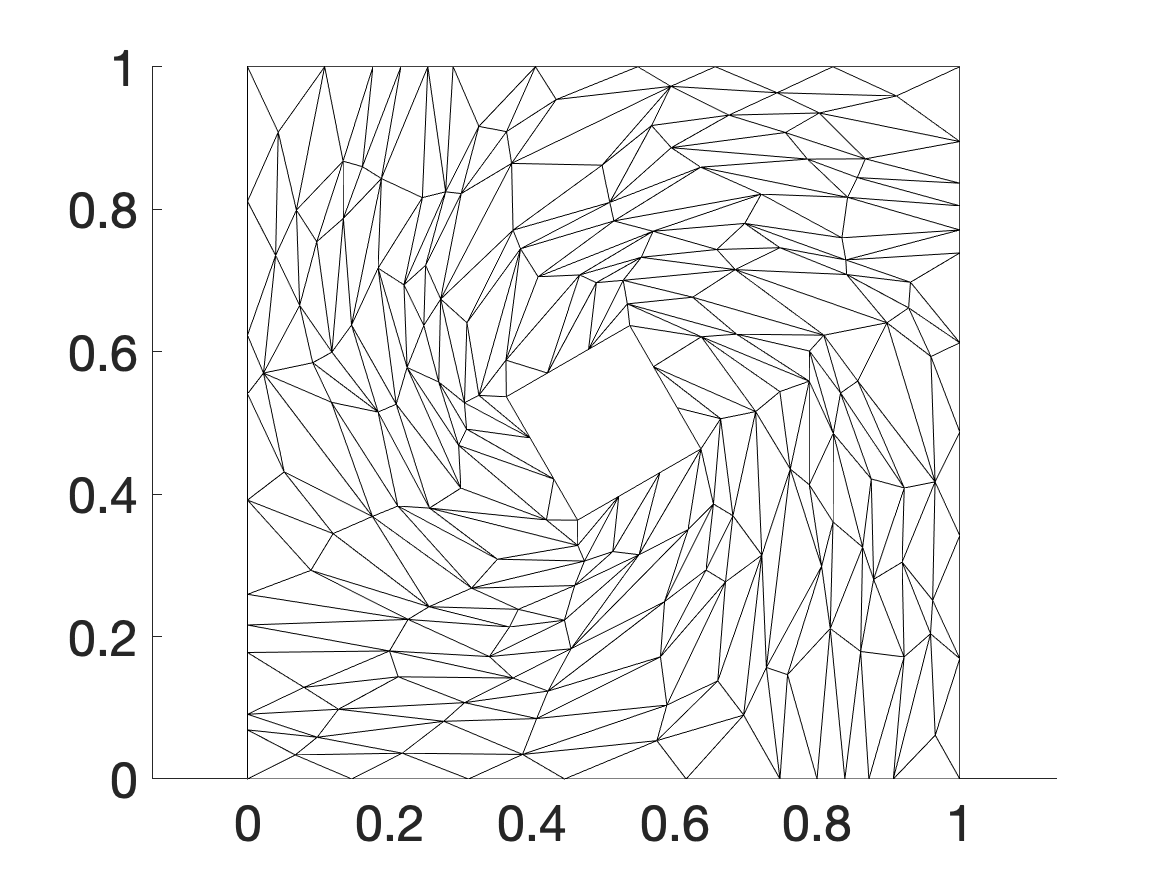}}
 
 \caption{large deformation mesh morphing.
 (a) reference configuration.
 (b)-(c) deformed meshes obtained using 
 \eqref{eq:morozov_registration} with exp-mesh objective, $\delta=10^{-6}$, $\kappa_{\rm msh}=10$, $n_{\rm lp}=25$.
}
 \label{fig:vis_mesh_square_square}
  \end{figure}

We consider the domain $\Omega = (0,1)^2 \setminus (0.4,0.6)^2$ in   Figure \ref{fig:vis_mesh_square_square}(a) and we consider a P2 FE triangulation $\mathcal{T}_{\rm hf}$ with linear elements of $\Omega$ that is also depicted in    Figure \ref{fig:vis_mesh_square_square}(a).
We fix the outer boundary and we rotate the inner boundary about the center $\bar{x}=[0.5,0.5]$ by an angle $\theta\in (0,120^o]$. The same test case was previously considered in    \cite{froehle2015nonlinear} to demonstrate the quality of an elasticity-based mesh deformation method.

We consider the registration statement
\eqref{eq:morozov_registration} with $\delta=10^{-6}$,   $\{  x_i\}_{i=1}^N$ equal to the points of the mesh on the interior boundary $\partial \Omega_{\rm in}$ and 
$\{ y_i = \mathbf{Rot}(\theta) (x_i-\bar{x})\}_{i=1}^N$; we consider the exp-jac objective with $\epsilon=0.05$ and the exp-mesh objective with $\kappa_{\rm msh}=10$. 
We consider  polynomials of degree  $n_{\rm lp}=25$.
In order to solve the optimization problem, we split the interval 
$(0,120^o]$ into $N_{\theta}=15$ equispace sub-intervals $(\theta_{k-1},\theta_{k})$ with $k=1,\ldots,N_{\theta}$; then, we consider a continuation strategy in which we initialize the optimizer for $\theta = \theta_k$ using the previously-computed solution for $\theta = \theta_{k-1}$.
 
\begin{figure}[h!]
\centering
\subfloat[]{
\begin{tikzpicture}[scale=0.6]
\begin{loglogaxis}[
xmode=linear,
ymode=log,
grid=both,
minor grid style={gray!25},
major grid style={gray!25},
xlabel = {\LARGE {$\theta$} },
  ylabel = {\LARGE {$q_{\rm min}$}},
  line width=1.2pt,
  mark size=3.0pt,
xmin=0,   xmax=120,
 ymin=0.01,   ymax=1,
 ]

\addplot[line width=1.pt,color=red,mark=square]  table {data/square_square/qmin_expmesh.dat}; 
      
\addplot[line width=1.pt,color=blue,mark=triangle*] table {data/square_square/qmin_expjac.dat};
         
\end{loglogaxis}
\end{tikzpicture}

}
\subfloat[]{
\begin{tikzpicture}[scale=0.6]
\begin{loglogaxis}[
xmode=linear,
ymode=log,
grid=both,
minor grid style={gray!25},
major grid style={gray!25},
xlabel = {\LARGE {$\theta$} },
  ylabel = {\LARGE {$J_{\rm min}$}},
  line width=1.2pt,
  mark size=3.0pt,
xmin=0,   xmax=120,
ymin=0.001,   ymax=1,
 ]

\addplot[line width=1.pt,color=red,mark=square]  table {data/square_square/Jmin_expmesh.dat}; 
      
\addplot[line width=1.pt,color=blue,mark=triangle*] table {data/square_square/Jmin_expjac.dat};
         
\end{loglogaxis}
\end{tikzpicture}
}
 \subfloat[]{
\begin{tikzpicture}[scale=0.6]
\begin{loglogaxis}[
xmode=linear,
ymode=log,
grid=both,
minor grid style={gray!25},
major grid style={gray!25},
xlabel = {\LARGE {$\theta$} },
ylabel = {\LARGE {nbr its}},
  line width=1.2pt,
  mark size=3.0pt,
xmin=0,   xmax=120,
ymin=1,   ymax=1000,
 ]

\addplot[line width=1.pt,color=red,mark=square]  table {data/square_square/nbrits_expmesh.dat}; 
      
\addplot[line width=1.pt,color=blue,mark=triangle*] table {data/square_square/nbrits_expjac.dat};
         
\end{loglogaxis}
\end{tikzpicture}
}
\caption[Caption in ToC]{ Large deformation mesh morphing.
Behavior of the minimum radius ratio $q_{\rm min}$, the minimum Jacobian $J_{\rm min}$, and the number of iterations for exp-jac and exp-mesh objective function and for several values of $\theta$.
($\delta=10^{-6}$, $\epsilon=0.05$, $\kappa_{\rm msh}=10$, $n_{\rm lp}=25$).
exp-mesh
\tikzref{three:delta1m6};
exp-jac
\tikzref{three:delta1m4}.
}
\label{fig:largedef_meshmorphing_thetasens}
\end{figure}
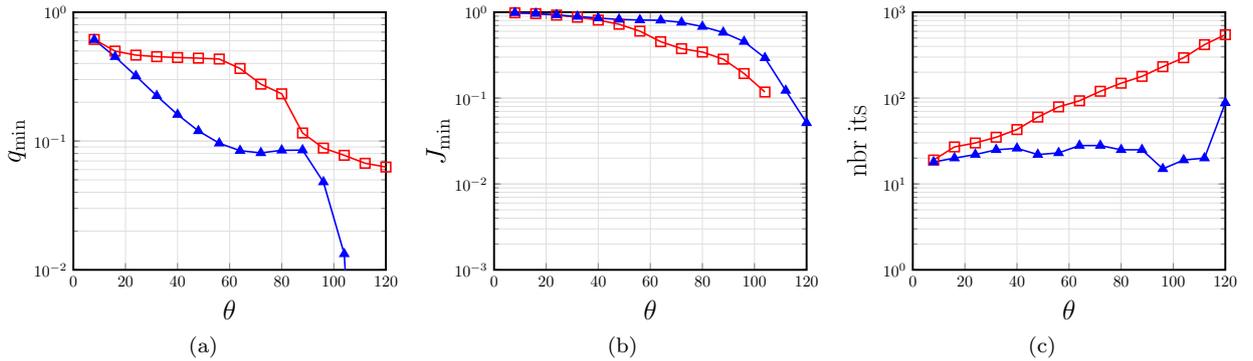

 Figures  \ref{fig:vis_mesh_square_square}(b) and (c) show the deformed meshes for $\theta=60^o$ and
 $\theta=120^o$ for the exp-mesh objective.
 Figures \ref{fig:largedef_meshmorphing_thetasens}(a)-(b)-(c) compare the behavior of the minimum radius $q_{\rm min}$, the minimum Jacobian determinant $J_{\rm min}$, and the number of iterations of the interior point method for the two choices of the objective function and for several values of $\theta$.
 
Some comments are in order. Both approaches succeed to deliver valid deformations with respect to their target --- that is,  the bijectivity of $\Phi$ for the exp-jac objective and discrete bijectivity for the exp-mesh objective.
Note, however, that for $\theta \gtrsim 100^o$ the solution to \eqref{eq:morozov_registration}  with 
exp-jac objective fails to deliver a valid deformed mesh; similarly, the solution to \eqref{eq:morozov_registration}  with 
exp-mesh objective is not bijective over $\Omega$. Recalling Remark \ref{remark:MtD_vs_DtM}, we conclude that the choice of the objective function should be driven by the choice of the strategy to deal with geometry variations --- namely, discretize-then-map or map-then-discretize.
 
\subsection{Geometry reduction of a parametric two-dimensional domain}
\label{sec:TT_sisc2020_revised}

\begin{figure}[h!]
\centering
 \subfloat[ ] 
{  \includegraphics[width=0.33\textwidth]
 {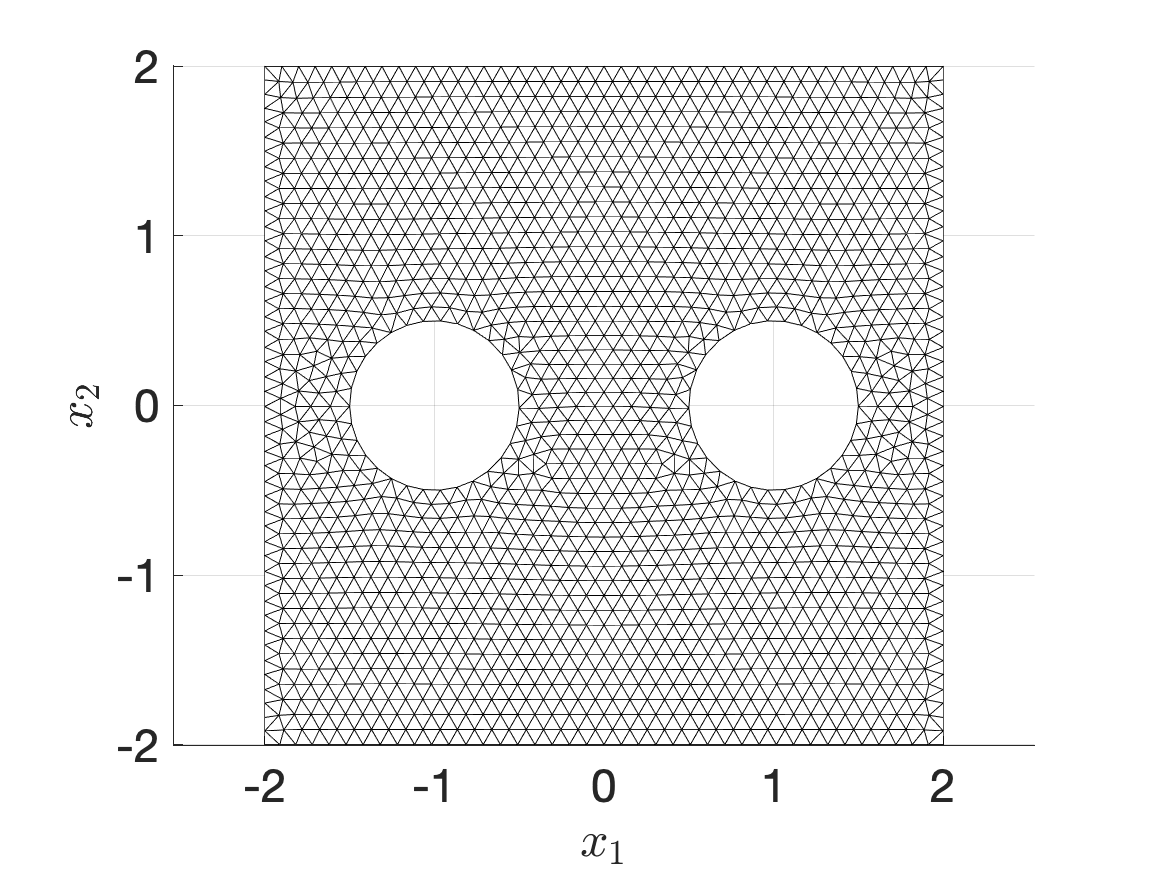}}
   ~~
 \subfloat[exp-mesh] 
{  \includegraphics[width=0.33\textwidth]
 {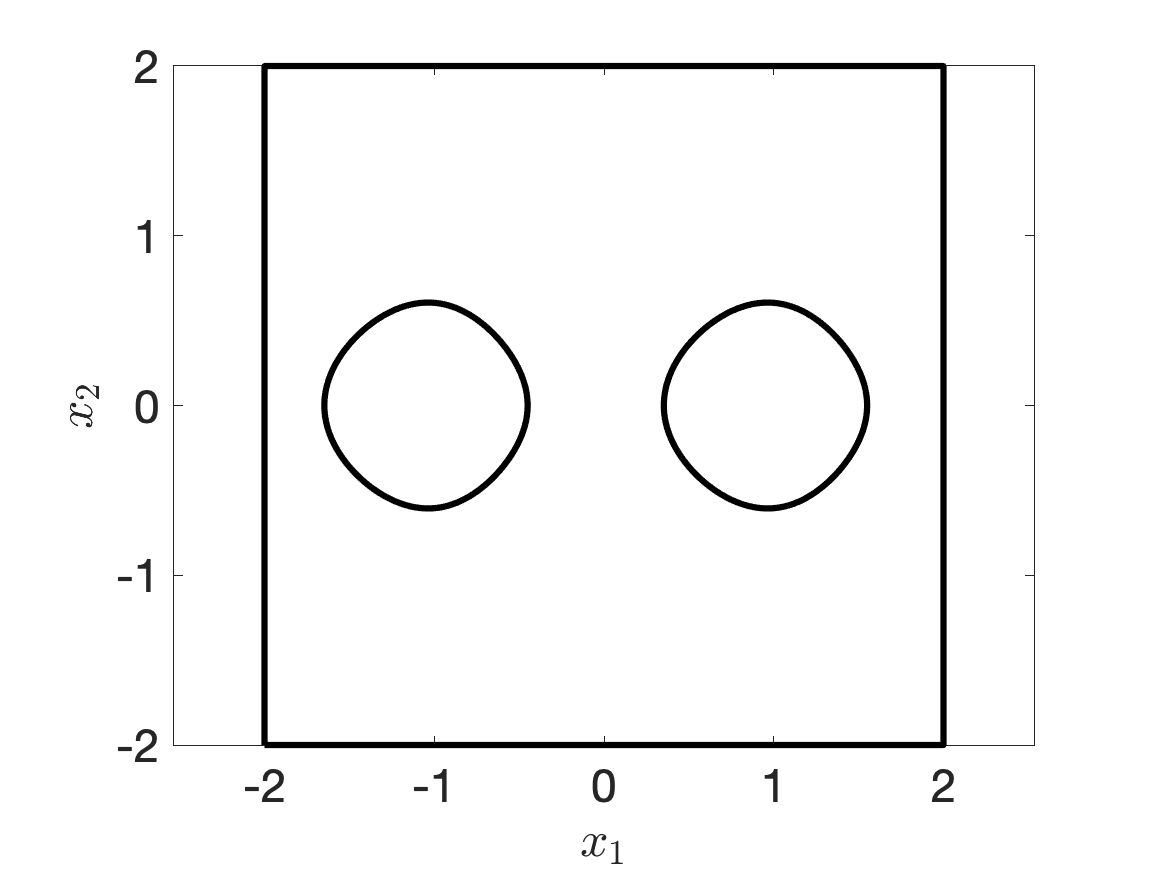}}
~~
 \subfloat[exp-jac ] 
{  \includegraphics[width=0.33\textwidth]
 {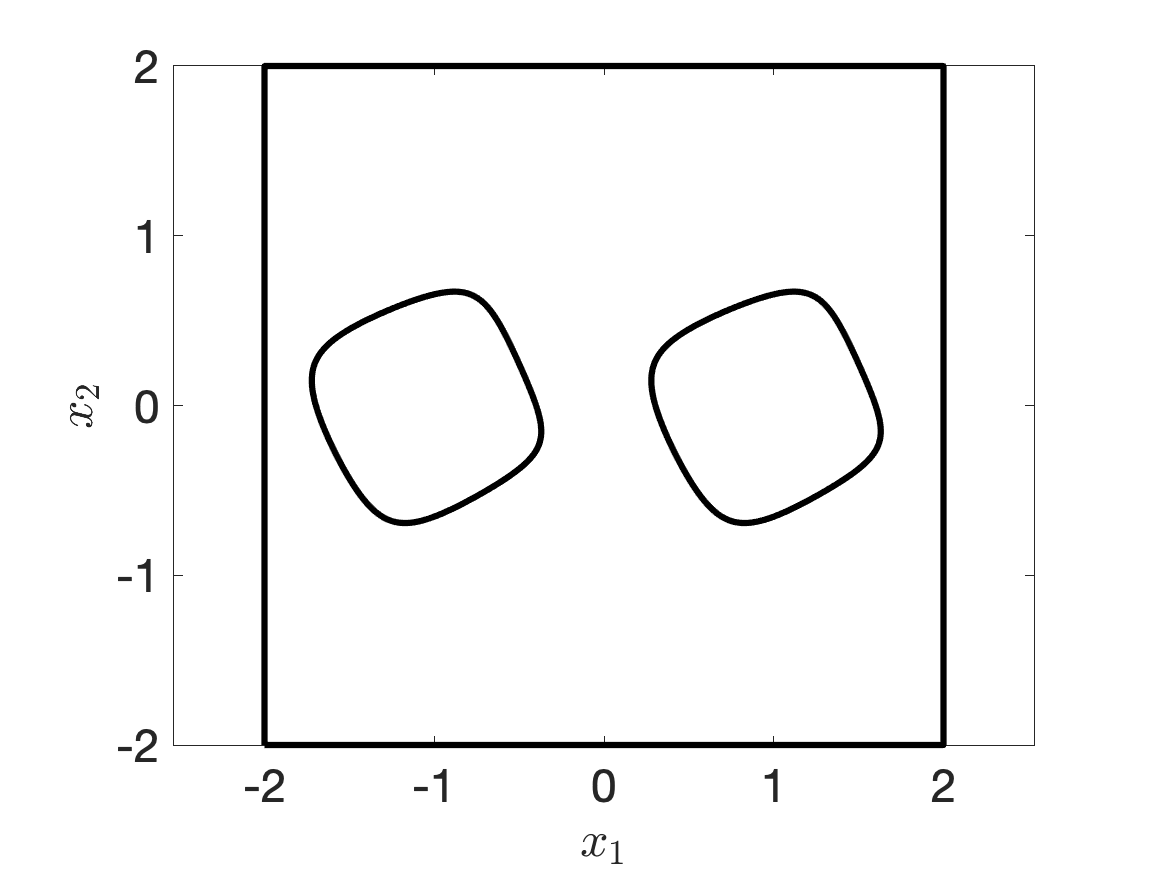}}
  
\caption{geometry reduction of a parametric two-dimensional domain. 
(a) reference configuration. 
(b) domain $\Omega_{\mu}$ 
for $\mu=[\nu_{\rm min},\nu_{\rm min}]$, with $(\nu_{\rm min})_j =\min_{\nu\in \mathcal{P}_{\rm in}} (\nu)_j$, $j=1,2,3$.
(c) domain $\Omega_{\mu}$ for 
$\mu=[\nu_{\rm max}^{(1)},\nu_{\rm max}^{(2)}]$, with $(\nu_{\rm max})_j =\max_{\nu\in \mathcal{P}_{\rm in}} (\nu)_j$, $j=1,2,3$.
}
\label{fig:geored_vis1}
 \end{figure}

We consider a variant of the parametric geometry reduction problem considered in
\cite{taddei2020registration}. We introduce $\Omega_{\rm box}=(-2,2)^2$ and we define the reference domain
$\Omega = \Omega_{\rm box} \setminus \bigcup_{i=1}^2 \Omega_{\rm in}^{(i)}$ where $\Omega_{\rm in}^{(1)}, \Omega_{\rm in}^{(2)}$ are two circles of radius $r=1/2$  centered in $\bar{x}_1 = (-1,0)$ and  
$\bar{x}_2 = (1,0)$, respectively. 
We parameterize the inner circles using the functions
$\gamma_{\rm ref,in}^{(i)}:[0,2\pi) \to \mathbb{R}^2$ with $\gamma_{\rm ref,in}^{(i)}(t) = \bar{x}_i + \frac{1}{2}\left[\cos(t) ,\sin(t) \right]$.
We further introduce the parameter domain $\mathcal{P}_{\rm in} = [0.1,0.4]^2\times [0,\pi/4]$ and the curves
\begin{equation}
\label{eq:geored_pbdefinition}
\boldsymbol{\gamma}_{\rm in, \nu^{(i)}}^{(i)}(t) = 
\bar{x}_i + \frac{1}{2}
\left[
\begin{array}{l}
\displaystyle{ \cos(t) \left(  1 + \nu_1^{(i)} \, \left( \cos (t + \nu_3^{(i)}) \right)^2 \, + \, 2 \cdot 10^{-3} \left( (2\pi - t) \, t \right)^2 \right) }
\\[3mm]
\displaystyle{ \sin(t) \left(  1 + \nu_2^{(i)} \, \left( \sin (t + \nu_3^{(i)}) \right)^2 \, + \, 2 \cdot 10^{-3} \left( (2\pi - t) \, t \right)^2 \right)  }
\\
\end{array}
\right],
\end{equation}
with $\bar{x}_1=[-1,0]$, $\bar{x}_2=[1,0]$ and $r=1/2$. Then, we define the vector of parameters
$\mu=[\nu^{(1)}, \nu^{(2)}]$ in the parameter region 
$\mathcal{P} = \mathcal{P}_{\rm in}\times \mathcal{P}_{\rm in}$ and we introduce the family of parameterized domains
\begin{equation}
\label{eq:two_hole_domains}
\Omega_{\mu} = \Omega_{\rm box} \setminus
\left( \Omega_{\rm in,\nu^{(1)} }^{(1)} \cup 
\Omega_{\rm in,\nu^{(2)} }^{(2)}   \right),
\;\; {\rm with} \;\;
\partial \Omega_{\rm in,\nu^{(i)}}^{(i)}
=
\gamma_{\rm in, \nu^{(i)}}^{(i)}([0,2\pi)), \;\; i=1,2.
\end{equation}
 Figure \ref{fig:geored_vis1}(a) shows the reference configuration and the mesh $\mathcal{T}_{\rm hf}$ employed for the numerical investigations;
Figures \ref{fig:geored_vis1}(b)  and (c) show two elements of the family $\{ \Omega_{\mu} : \mu\in \mathcal{P} \}$.

We consider $n_{\rm lp}=20$ and we enforce $(\Phi-\texttt{id})\cdot \mathbf{n} |_{\partial \Omega_{\rm box}}= 0$ to ensure bijectivity over $\Omega_{\rm box}$ --- note that the resulting   affine space $\mathcal{W}_{\rm hf}$ is of dimension $M=720$.
We generate training and test sets $\mathcal{P}_{\rm train}, \mathcal{P}_{\rm test} \subset \mathcal{P}$ of cardinality $n_{\rm train}=100$ and 
$n_{\rm test}=20$ using independent realizations of an  uniform random variable over $\mathcal{P}$.
We apply registration with exp-mesh objective function ($\kappa_{\rm msh}=10$) and we describe the internal reference and parameterized domains 
using $N_{\rm v}=100$ points,
\begin{equation}
\label{eq:sorted_points_geored}
\begin{array}{l}
\displaystyle{
\left\{
x_{i:=k+(j-1)N_{\rm v}} = \gamma_{\rm in,ref}^{(j)}(t_k) : \, k=1,\ldots,N_{\rm v}, j=1,2
\right\}, 
\quad
{\rm with} \;\; t_k = \frac{k-1}{N_{\rm v}}}\\[3mm]
\displaystyle{
\left\{
y_{i:=k+(j-1)N_{\rm v}}(\mu) = \gamma_{\rm in,\nu^{(j)}}^{(j)}(t_k) : \, k=1,\ldots,N_{\rm v}, j=1,2
\right\}
}.\\[3mm]
\end{array}
\end{equation}
Note that the reference and target points are sorted \emph{a priori} and that $N=2 N_{\rm v}=200$.  In all our tests, we initialize the registration algorithm with $\Phi=\texttt{id}$.

We  here design a GRR algorithm that takes as input the points $\{x_i \}_{i=1}^N$  and $\{y_i^{\rm raw} \}_{i=1}^Q$ and the mesh $\mathcal{T}_{\rm hf}$ of $\Omega$ and returns a deformed mesh $\Phi(\mathcal{T}_{\rm hf})$,
\begin{equation}
\label{eq:abstract_morphing_algorithm}
\Phi(  \mathcal{T}_{\rm hf}  )
\, = \,
{\rm registration} \left(
\mathcal{T}_{\rm hf}, \;
\{x_i \}_{i=1}^N, \; 
\{y_i^{\rm raw} \}_{i=1}^Q
\right).
\end{equation} 
We do not exploit the knowledge of the parameterization.
We first consider the simplified case of ``sorted'' points (cf. \eqref{eq:sorted_points_geored}) in which we feed the algorithm \eqref{eq:abstract_morphing_algorithm} with 
$\{y_i^{\rm raw}=y_i  \}_{i=1}^{Q=N}$; then, we consider the more challenging problem of unsorted and sparse points in which we feed  \eqref{eq:abstract_morphing_algorithm} with 
$$
\{y_i^{\rm raw}=y_{\texttt{I}_i}  \}_{i=1}^Q,
\quad
Q = 0.8 N, \;\;
\texttt{I} \subset \{1,\ldots,N\}
$$
where $\texttt{I} $ is obtained by selecting the first $Q$ points of a random permutation of $\{1,\ldots,N\}$. We envision that the unsorted problem is an adequate proxy of the geometry registration problems encountered in practical applications.
 Table \ref{tab:details_optimization_geored}
 summarizes the choice of the hyper-parameters considered in the numerical experiments.
 Here, the geometric error is practically estimated as
\begin{equation}
\label{eq:geo_error_geored}
E_{\rm geo}(\Phi) = \max_{i=1,\ldots,N} \min_{j=1,\ldots, 5N} \| \Phi(x_i) - y_j^{\rm test}  \|_2,
\end{equation} 
 where $\{ y_j^{\rm test} \}_{j=1}^{5 N}$ are
defined  as in  \eqref{eq:sorted_points_geored} with $N_{\rm v}=500$. 
 
\begin{table}[h!]
\begin{tabular}{|l|l|l|}
\hline
Method & Hyper-parameters  & Optimization method  \\
\hline
Tykhonov 
\eqref{eq:mesh_morphing_optimization_tikhonov}
&  exp-mesh objective, $\kappa_{\rm msh}=10$,
$\xi \in \{10^{-4}, 10^{-5}\}$ & quasi-Newton \\
\hline
Morozov 
\eqref{eq:mesh_morphing_optimization_morozov}

&  exp-mesh objective,
$\delta \in \{10^{-3}, 10^{-4}\}$
 &
interior point (linear constraints)\\
\hline
Inverted 
\eqref{eq:inverted_formulation}
 & exp-mesh  constraint, 
$\xi \in \{10^{-4}, 10^{-5}\}$, $\delta_{\rm con}=1$ & interior point (nonlinear constraints)\\
\hline
Coherent point drift
&
$\beta=1$, $\lambda=1$, $w=0$ 
&
 Algorithm \ref{alg:CPD},  Appendix \ref{sec:CPD_appendix}
\\
\hline
\end{tabular}
\caption{geometry reduction of a parametric two-dimensional domain. Details of the optimization method.}
\label{tab:details_optimization_geored}
\end{table}
 
Table \ref{tab:training_geored} and Table \ref{tab:test_geored} 
summarize the average and worst performance of the registration methods of   section \ref{sec:geometry_reduction} on the training and the  test set, respectively. 
In more detail, in 
 Table \ref{tab:training_geored}, 
we consider the statements
\eqref{eq:mesh_morphing_optimization_tikhonov},
\eqref{eq:mesh_morphing_optimization_morozov}, and
\eqref{eq:inverted_formulation}: 
Tykhonov regularization enables faster predictions, particularly in the worst-case cost;
Morozov regularization enables a sharp control of the geometric error --- which is in the order of $\sqrt{2} \delta$ for all numerical experiments --- at the price of significantly larger worst-case costs;
the inverted formulation leads to results that are close to Tykhonov regularization.
In Table \ref{tab:test_geored},  we compare performance of the same registration algorithms over the test set based on the full space ($M=720$) and on the reduced space obtained using POD with tolerance 
$tol_{\rm pod}=10^{-5}$ (cf. 
\eqref{eq:POD_cardinality_selection}):
we set $\xi=10^{-5}$ in 
\eqref{eq:mesh_morphing_optimization_tikhonov} and 
\eqref{eq:inverted_formulation}, and 
$\delta=10^{-4}$ in 
\eqref{eq:mesh_morphing_optimization_morozov}.
 we observe that POD enables speedups in the order of $\mathcal{O}(100)$ without any significant deterioration of performance. 
Even more, we observe that dimensionality reduction slightly improves reconstruction performance for certain configurations: this empirical finding suggests that the reduction of the number of unknowns simplifies the optimization task and might prevent convergence to suboptimal local minima.

\begin{table}[h!]
\begin{tabular}{|l|c|c|c|c|c|c|c|c|}
\hline
\textbf{Training (sorted)} & 
\multicolumn{2}{|c|}{$q_{\rm min}$ }
 &  
\multicolumn{2}{|c|}{ geo error} 
&
\multicolumn{2}{|c|}{nbr its} 
&
\multicolumn{2}{|c|}{cost [s]} \\
 & avg & min &  avg & max   &  avg  & max &  avg  & max   \\
\hline
Tykhonov ($\xi=10^{-4}$)
& $0.44$ &  $0.32$   &   $9.53\cdot 10^{-5}$ & $1.98 \cdot 10^{-4}$ &  $287.5$ &  $401$  & $31.85$ &  $44.04$\\
\hline
Tykhonov ($\xi=10^{-5}$)
 & $0.42$&  $0.27$  
& $8.18 \cdot 10^{-5}$    &  $1.17 \cdot 10^{-4}$   &  $274.9$ &  $383$  & $30.62$ & $42.39$ \\
\hline
Morozov     ($\delta=10^{-3}$)
&  $0.45$ &  $0.34$    &  $1.40 \cdot 10^{-3}$  & $1.42 \cdot 10^{-3}$  &  $78.6$  &  $176$  & $26.58$ &  $68.18$ \\
\hline
Morozov     ($\delta=10^{-4}$)
&  $0.45$ &  $0.34$    &  $1.36 \cdot 10^{-4}$  & $1.42 \cdot 10^{-4}$  &  $75.7$  &  $211$  & $30.06$ &  $108.04$ \\
\hline
Inverted   ($\xi=10^{-4}$)
&  $0.41$ &  $0.24$    &  
$3.37 \cdot 10^{-4}$  &
$6.20 \cdot 10^{-4}$  &  $199.6$    &  $285$    & $33.03$ &  $49.38$ \\
\hline
Inverted   ($\xi=10^{-5}$)
&  $0.41$ &  $0.24$    
  &  $8.43 \cdot 10^{-5}$   &  $1.29 \cdot 10^{-4}$  &  $204.14$  &  $349$    & $32.92$ &    $55.89$\\
\hline
\end{tabular}
\medskip
 
\caption{geometry reduction of parametric two-dimensional domain. Registration performance on the  training set $\mathcal{P}_{\rm train}$ for sorted data, $M=720$.}
\label{tab:training_geored}
\end{table}

\begin{table}[h!]
\begin{tabular}{|l|c|c|c|c|c|c|c|c|}
\hline
\textbf{Test (sorted)} & 
\multicolumn{2}{|c|}{$q_{\rm min}$ }
 &  
\multicolumn{2}{|c|}{ geo error} 
&
\multicolumn{2}{|c|}{nbr its} 
&
\multicolumn{2}{|c|}{cost [s]} \\
 & avg & min &  avg & max   &  avg  & max &  avg  & max   \\
\hline
Tykhonov, full ($M=720$)
& $0.41$ &  $0.30$  
& $0.81 \cdot 10^{-4}$   & $1.01\cdot 10^{-4}$   &  $270.0$ &  $378$  &  $30.10$ &  $41.91$ \\
\hline
Tykhonov, reduced  ($M=20$)
& $0.42$ &  $0.31$   
& $0.73\cdot 10^{-4}$  & $1.01\cdot 10^{-4}$   &  $59.3$ &  $91$  &  $0.23$ &  $0.31$ \\
\hline
Morozov, full  ($M=720$)
& $0.44$ &  $0.34$    & $1.37\cdot 10^{-4}$ & $1.40 \cdot 10^{-4}$  &  $68.3$ &  $86$  &  $26.31$ &  $73.12$ \\
\hline
Morozov, reduced ($M=21$)
& $0.44$ &  $0.35$    & $1.35\cdot 10^{-4}$ & $1.42 \cdot 10^{-4}$  &  $47.1$ &  $61$  &  $0.51$ &  $0.63$ \\
\hline
Inverted,  full   ($M=720$)
& $0.41$ &  $0.27$    & $0.79\cdot 10^{-4}$ & $1.13 \cdot 10^{-4}$  &  $217.0$ &  $320$  &  $34.68$ &  $51.71$ \\
\hline
Inverted, reduced  ($M=21$)
& $0.40$ &  $0.26$    & $0.69\cdot 10^{-4}$ & $1.20 \cdot 10^{-4}$  &  $77.5$ &  $117$  &  $0.35$ &  $0.51$ \\
\hline
\end{tabular}
\medskip 

\caption{geometry reduction of a parametric two-dimensional domain. Registration performance on the  test set $\mathcal{P}_{\rm test}$ for sorted data.}
\label{tab:test_geored}
\end{table}

\begin{figure}[h]
\centering

 \subfloat[Tykhonov, full, sorted] 
{  \includegraphics[width=0.32\textwidth]
 {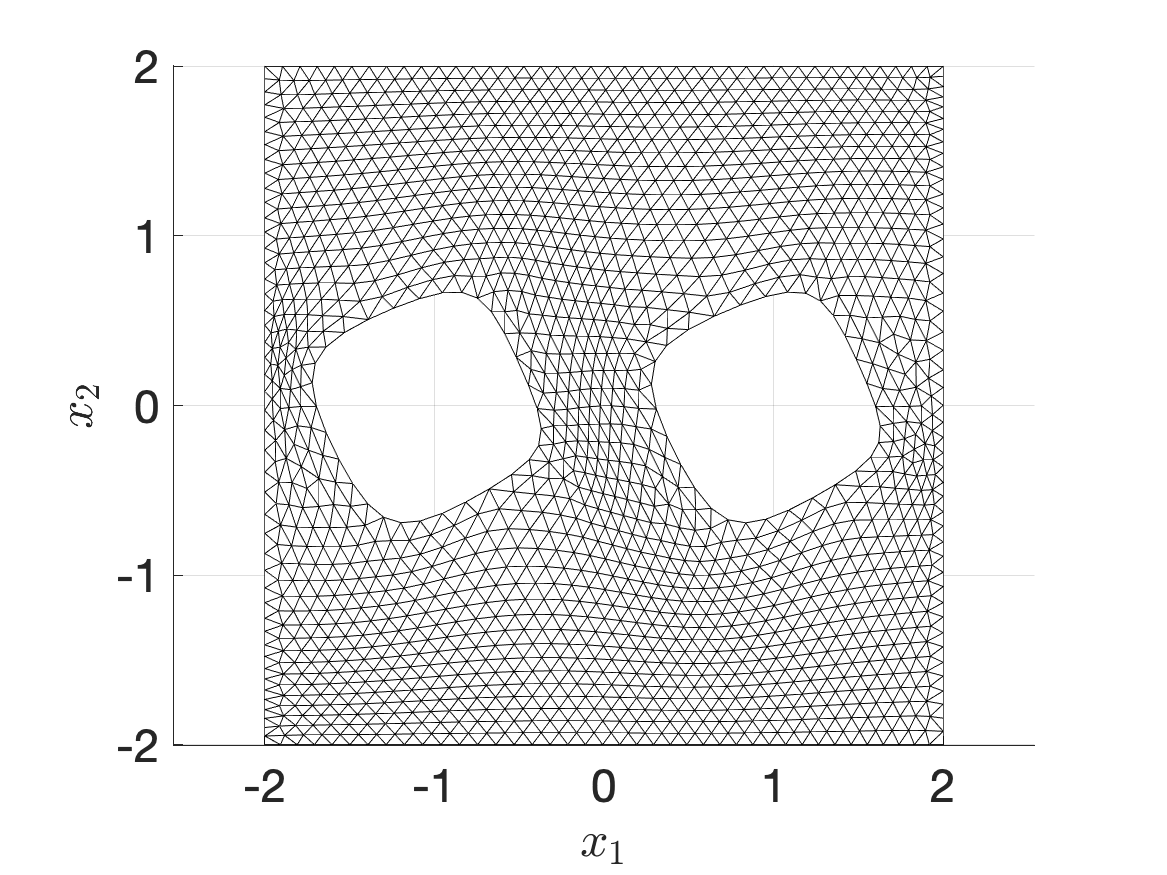}}
~~
 \subfloat[Morozov, full, sorted] 
{  \includegraphics[width=0.32\textwidth]
 {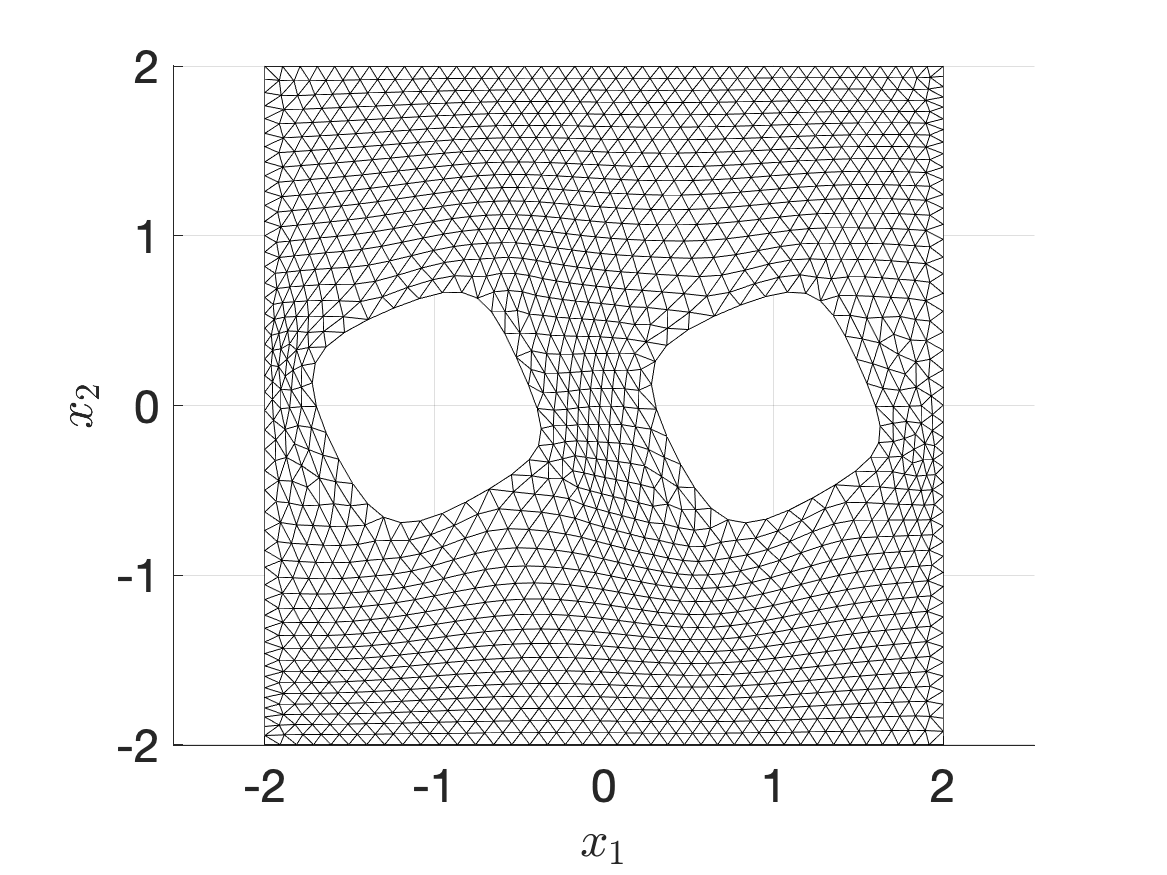}}
~~
 \subfloat[Tykhonov, reduced, unsorted] 
{\includegraphics[width=0.32\textwidth]
 {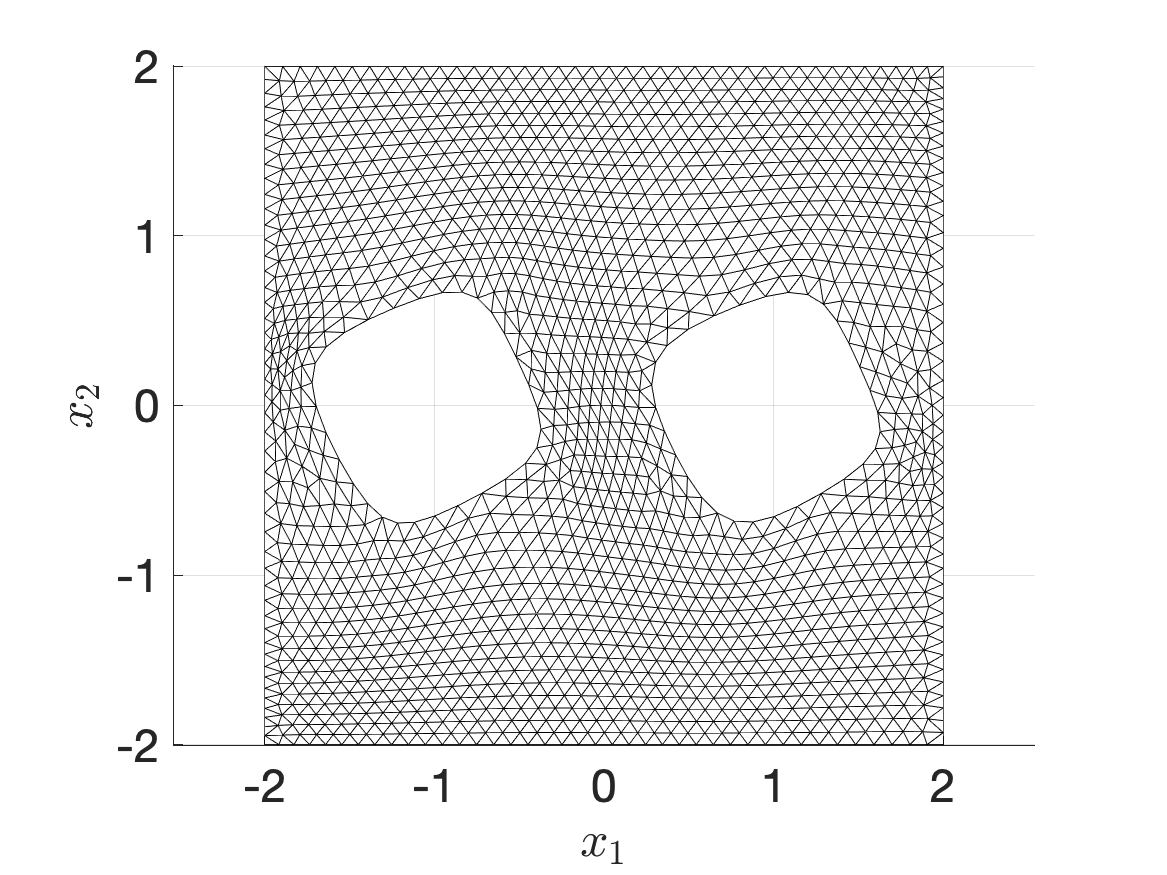}}
 
 \caption{geometry reduction of a parametric two-dimensional domain. Deformed mesh obtained using three different methods 
 for $\mu=[\nu_{\rm max},\nu_{\rm max}]$.
}
 \label{fig:vis_geored}
  \end{figure}

Table \ref{tab:test_geored_unsorted} shows performance for the unsorted problem;
in this test,  we consider unsorted data both for training and assessment, and we resort to the full CPD space.
  We observe that both Tykhonov and the inverted formulation lead to performance that are comparable to the performance obtained using CPD in terms of geometric error. Note also that the number of POD modes required to achieve the tolerance $tol_{\rm pod}=10^{-5}$  in 
\eqref{eq:POD_cardinality_selection} is significantly larger (see also  Figure \ref{fig:vis_geored}) for the unsorted dataset: this empirical finding suggests  that the application of CPD to the unsorted dataset with missing points  might introduce numerical noise that  ultimately hinders the compressibility of the manifold associated with the mapping $\Phi$.

\begin{table}[h!]

\begin{tabular}{|l|c|c|c|c|c|c|c|c|}
\hline
\textbf{Test (unsorted)} & 
\multicolumn{2}{|c|}{$q_{\rm min}$ }
 &  
\multicolumn{2}{|c|}{ geo error} 
&
\multicolumn{2}{|c|}{nbr its} 
&
\multicolumn{2}{|c|}{cost [s]} \\
 & avg & min &  avg & max   &  avg  & max &  avg  & max   \\
\hline
CPD 
& ---  &  ---   &  $6.68\cdot 10^{-3}$ & $8.61 \cdot 10^{-3}$  &  $53.0$ & $97$  & $0.03$ &  $0.05$ \\
\hline
Tykhonov, full ($M=720$)
& $0.43$  &    $0.29$    & $6.93\cdot 10^{-3}$ & $8.61 \cdot 10^{-3}$    &  $334.5$ &  $401$  & $37.39$ &  $44.75$ \\
\hline
Tykhonov, reduced ($M=48$)
&  $0.43$  &   $0.29$   & $6.93\cdot 10^{-3}$ & $8.67 \cdot 10^{-3}$    &  $181.1$ &  $237$  & $1.19$ &  $1.52$ \\
\hline
Inverted, full ($M=720$)
&  $0.40$   &  $0.26$   & $7.03\cdot 10^{-3}$ & $8.86 \cdot 10^{-3}$    &  $208.0$ &  $245$  & $33.91$ &  $41.41$ \\
\hline
Inverted, reduced ($M=40$)
&  $0.40$   &  $0.26$   & $7.32\cdot 10^{-3}$ & $9.03 \cdot 10^{-3}$    &  $130.3$ &  $143$  & $0.85$ &  $1.02$ \\
\hline
\end{tabular}
\caption{geometry reduction of a parametric two-dimensional domain. Registration performance on the  test set $\mathcal{P}_{\rm test}$ for unsorted data, with $80\%$ training data 
($\xi=10^{-4}$).
CPD is based on the full space  \eqref{eq:CPD_affine_space}.}
\label{tab:test_geored_unsorted}
\end{table}

\begin{figure}[h]
\centering
\subfloat[Tykhonov]{
\begin{tikzpicture}[scale=0.7]
\begin{loglogaxis}[
xmode=linear,
ymode=log,
grid=both,
minor grid style={gray!25},
major grid style={gray!25},
xlabel = {\LARGE {$k$} },
ylabel = {\LARGE {$1 - {\rm EC}_k(\boldsymbol{\lambda})$}},
legend entries = {\LARGE {sorted}, \LARGE {unsorted}},
  line width=1.2pt,
  mark size=3.0pt,
xmin=1,   xmax=100,
 ymin=0.0000000001,   ymax=1,
 ]

\addplot[line width=1.pt,color=red,mark=square]  table {data/param_geo/lambdas_sorted.dat}; 
       
\addplot[line width=1.pt,color=blue,mark=triangle*] table {data/param_geo/lambdas_unsorted.dat};
   
\end{loglogaxis}
\end{tikzpicture}
}
~~
\subfloat[Inverted]{
\begin{tikzpicture}[scale=0.7]
\begin{loglogaxis}[
xmode=linear,
ymode=log,
grid=both,
minor grid style={gray!25},
major grid style={gray!25},
xlabel = {\LARGE {$k$} },
ylabel = {\LARGE {$1 - {\rm EC}_k(\boldsymbol{\lambda})$}},
legend entries = {\LARGE {sorted}, \LARGE {unsorted}},
  line width=1.2pt,
  mark size=3.0pt,
xmin=1,   xmax=100,
 ymin=0.0000000001,   ymax=1,
 ]

\addplot[line width=1.pt,color=red,mark=square]  table {data/param_geo/lambdas_sorted_inverted.dat}; 
       
\addplot[line width=1.pt,color=blue,mark=triangle*] table {data/param_geo/lambdas_unsorted_inverted.dat};
   
\end{loglogaxis}
\end{tikzpicture}
}

 \caption{geometry reduction of a parametric two-dimensional domain.
POD eigenvalues for sorted and unsorted data for 
Tykhonov  and inverted formulations
 with $\xi=10^{-4}$; CPD   is based on the full space  \eqref{eq:CPD_affine_space}.}
 \label{fig:pod_eigenvalues}
  \end{figure}
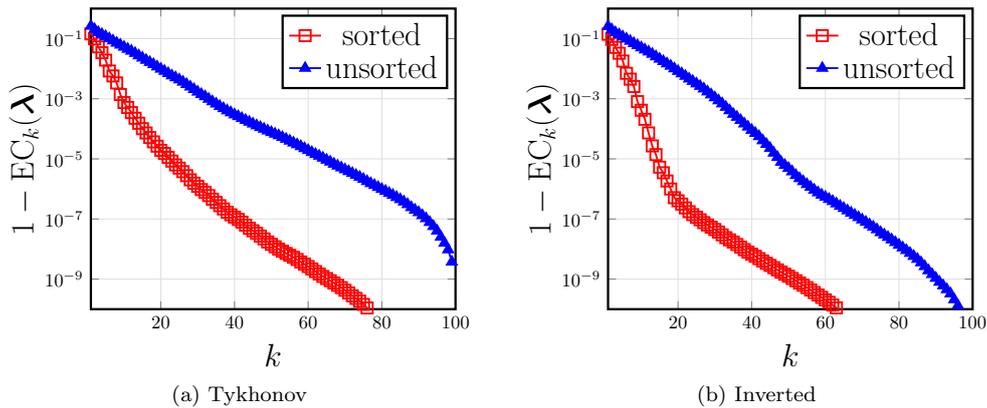

Figure \ref{fig:vis_geored}  shows the deformed meshes obtained for an out-of-sample parameter using 
Tykhonov regularization with full space and sorted data ($\xi=10^{-5}$),
Morozov regularization with full space and sorted data ($\delta=10^{-4}$),
Tykhonov regularization (in combination with CPD) with  reduced space and unsorted data, respectively  ($\xi=10^{-5}$, $M=48$).
For this test case, we obtain that the geometric error 
\eqref{eq:geo_error_geored} is given by 
$7.8\cdot 10^{-5}, 1.4\cdot 10^{-4}, 5.9\cdot 10^{-3}$, respectively; the geometric error of CPD is 
$5.9\cdot 10^{-3}$.

Table \ref{tab:test_geored_sortedunsorted}
investigates performance of the registration method based on sorted training on a  unsorted test dataset: we envision that this scenario is of interest for real-time applications for which online data are considerably more noisy than data used for training. 
We consider the Tykhonov formulation with $\xi=10^{-4}$; we perform POD with tolerance $tol_{\rm pod}=10^{-5}$ to obtain a reduced space of cardinality $M=22$; then, we compare performance obtained using  CPD based on the full space
\eqref{eq:CPD_affine_space} with the performance obtained using   CPD with the reduced space
\eqref{eq:CPD_reduced_space}.
We observe that dimensionality reduction in the CPD algorithm has a beneficial effect on the accuracy of the registration procedure and also reduces the cost of CPD.

\begin{table}[h!]
\begin{tabular}{|l|c|c|c|c|c|c|c|c|}
\hline
\textbf{Full-space CPD} & 
\multicolumn{2}{|c|}{$q_{\rm min}$ }
 &  
\multicolumn{2}{|c|}{ geo error} 
&
\multicolumn{2}{|c|}{nbr its} 
&
\multicolumn{2}{|c|}{cost [s]} \\
 & avg & min &  avg & max   &  avg  & max &  avg  & max   \\
\hline
CPD
&   &    
& $6.68 \cdot 10^{-3}$   & $8.61\cdot 10^{-3}$   &  $53.0$ &  $97$  &  $0.04$ &  $0.08$ \\
\hline
Tykhonov, reduced  ($M=22$)
& $0.46$ &  $0.30$   
& $2.95\cdot 10^{-2}$  & $5.91\cdot 10^{-2}$   &  $95.3$ &  $119$  &  $0.33$ &  $0.39$ \\
\hline 
\end{tabular}
\medskip 

\begin{tabular}{|l|c|c|c|c|c|c|c|c|}
\hline
\textbf{Reduced-space CPD} & 
\multicolumn{2}{|c|}{$q_{\rm min}$ }
 &  
\multicolumn{2}{|c|}{ geo error} 
&
\multicolumn{2}{|c|}{nbr its} 
&
\multicolumn{2}{|c|}{cost [s]} \\
 & avg & min &  avg & max   &  avg  & max &  avg  & max   \\
\hline
CPD
&   &    
& $6.09 \cdot 10^{-3}$   & $1.13\cdot 10^{-2}$   &  $26.4$ &  $39$  &  $0.01$ &  $0.02$ \\
\hline
Tykhonov, reduced  ($M=22$)
& $0.41$ &  $0.34$   
& $3.29\cdot 10^{-3}$  & $1.14\cdot 10^{-2}$   &  $68.4$ &  $116$  &  $0.24$ &  $0.38$ \\
\hline 
\end{tabular}
\medskip 

\caption{geometry reduction of parametric two-dimensional domain. Registration performance on the  test set $\mathcal{P}_{\rm test}$ for unsorted data; training performed on sorted data, 
$\xi=10^{-4}$, $tol_{\rm pod}=10^{-5}$.}
\label{tab:test_geored_sortedunsorted}
\end{table}


\section{Summary and discussion}
\label{sec:conclusions}
We presented an optimization-based approach to the problem of geometry registration and reduction (GRR). GRR is of paramount importance in pMOR to deal with PDEs in parametric geometries.
We presented a thorough mathematical analysis that offers the theoretical foundations for the methodology; we further presented thorough numerical investigations for three two-dimensional model problems.

The analysis of section \ref{sec:analysis} rigorously justifies the proposed approach for smooth $C^2$ domains; it also illustrates the issues that the approach faces when dealing with slender bodies and/or domains with corners. In particular, ensuring appropriate approximation --- via interpolation --- of the target shape at corners and cusps is key to  recover near-optimal bounds for the geometry error.

The numerical results of section \ref{sec:numerics} provide a thorough overview of the performance of   optimization-based registration methods, for several choices of the objective function and of the parameters.
First, Tykhonov regularization appears to be superior in terms of computational performance at the price of a less straightforward choice of the free parameter $\xi$ (compared to Morozov regularization).
In addition, the combination with a representative PSR procedure did not lead to any unstable behavior and/or lack of convergence.
Second, we observe that potential approximations of the form $\Phi=\texttt{id} +\nabla \phi$ are highly suboptimal
(cf. Figure \ref{fig:three_potvsfull}) in our setting:
this is in striking contrast with the related problem of optimal transportation for which the solution is guaranteed to be the gradient of a convex function \cite{brenier1987polar}.
Third, the results of section \ref{sec:TT_sisc2020_revised}  show that geometry reduction based on POD might greatly reduce the cost of registration. In this respect, the development of specialized hyper-reduction (see, e.g.,   \cite{farhat2021computational,yano2021model}) techniques might lead to much more significant computational gains for large-scale problems.
Fourth, we observed that the   deformation is highly sensitive to the choice of the objective function (see, e.g., Figure \ref{fig:vis_mesh_three} and  Figure \ref{fig:largedef_meshmorphing_thetasens}):
for projection-based pMOR applications, 
this choice should thus be driven by the strategy to deal with geometry variations --- namely, discretize-then-map or map-then-discretize.
Fifth, we observe that the approach can also cope with unsorted data when trained on sorted data (cf.
Table \ref{tab:test_geored_sortedunsorted}): towards this end, the application of the reduced-space CPD seems crucial to achieve accurate reconstructions.

\section*{Acknowledgement}
The author thanks Prof. Angelo Iollo (University of Bordeaux, Inria) for fruitful discussions; he also thanks Dr. Pierre Mounoud (University of Bordeaux) for extensive discussions on the proofs of section \ref{sec:analysis}.

\appendix

\section{Coherent point drift}
\label{sec:CPD_appendix}

In this Appendix, we review the coherent point drift (CPD) algorithm first introduced in
\cite{myronenko2010point}.
CPD is a probabilistic method based on a Gaussian mixture model (GMM): the PSR problem is formulated as a maximum likelihood estimation problem with a motion coherence constraint and is solved using an expectation-maximization (EM) procedure.
First, we briefly present EM algorithms for maximum likelihood estimation; then, we present the probabilistic model of CPD and we discuss the application of EM; finally, we comment on dimensionality reduction.
In  Appendix \ref{sec:EM_intro} and Appendix \ref{sec:prob_model_CPD}, we denote by $\mathbf{Y}=[y_1,\ldots,y_Q]^{\star}$ the observed data to shorten notation.

\subsection{Expectation maximization procedures}
\label{sec:EM_intro}
The EM algorithm is a prominent method in machine learning to find maximum likelihood solutions for models with latent variables; EM is used when the maximization of the likelihood of the observed variables is difficult but can be made easier by enlarging the sample with latent (unobserved) data.
We refer to \cite{bishop2006pattern} Chapter 9, 
 and to 
 \cite[Chapter 8]{hastie2009elements}    for a thorough review of the methodology and for a rigorous theoretical justification.

We denote by $\boldsymbol{\Theta}$ the parameter values of the model, by $\mathbf{Y}=[y_1,\ldots,y_Q]^{\star}   \in \mathbb{R}^{Q\times d}$ the observed data and by
$\mathbf{Z}=[z_1,\ldots,z_Q]^{\star} \in \{1,\ldots,N\}^Q$ the latent variables  for some $N\in \mathbb{N}$ --- we here consider the case of a single  discrete latent variable; however, the method can also deal with  continuous latent variables.
We denote by $p(\mathbf{Y},\mathbf{Z} | \Theta)$ the joint likelihood; recalling the Bayes' theorem, we find the expression for the posterior of $\mathbf{Z}$:
\begin{equation}
\label{eq:posterior_Z}
{\rm Pr}(\mathbf{Z}| \mathbf{Y}, \boldsymbol{\Theta})
=
\frac{ {\rm Pr}(\mathbf{Y},\mathbf{Z} |\boldsymbol{\Theta})   }{{\rm Pr}(\mathbf{Z}| \Theta)}.
\end{equation}

Assuming that the samples are independent identically distributed, the EM procedure  can be stated as in Algorithm \ref{alg:EM}. Note that, since the logarithm acts directly on the joint distribution, the M-step maximization might be tractable. The choice of the expectation in the M-step is rigorously justified by interpreting EM as an alternative  minimization algorithm (see,
\cite[Chapter 9.4]{bishop2006pattern}  
\cite[Chapter 8.5.3]{hastie2009elements}). Several authors have also considered to partition the parameters $\boldsymbol{\Theta}$ into groups and then break down the M-step into multiple steps each of which  involves optimizing one of the subsets with the remainders held fixed 
\cite{meng1993maximum}.

\begin{framed}
\captionof{algorithm}{Expectation-Maximization algorithm. \label{alg:EM}}
\begin{algorithmic}[1]
\State
Set $\boldsymbol{\Theta}^{(k)} = \boldsymbol{\Theta}_0$.
\medskip
			
\For {$k=1,2,   \ldots, \texttt{until  convergence}$ }
\medskip		
 
 \State
 \textbf{E-step:} Evaluate 
 $(\mathbf{P})_{i,j} =   {\rm Pr}(Z=z_j | \mathbf{Y}=y_i, \boldsymbol{\Theta}^{(k)} )$  using \eqref{eq:posterior_Z}.
 \medskip
 
  \State
 \textbf{M-step:} Maximize the expectation of the complete data log-likelihood 
 
 $\displaystyle{\boldsymbol{\Theta} \mapsto 
 \mathfrak{Q}(\boldsymbol{\Theta}; \mathbf{P}):=
 \sum_{i=1}^Q \sum_{j=1}^K
 (\mathbf{P})_{i,j} \log \left(
 {\rm Pr}(Z=z_i, Y = y_j |  \boldsymbol{\Theta})
 \right)}
 $ to find $\boldsymbol{\Theta}^{(k+1)}$.

 \medskip
 
 \EndFor
\end{algorithmic}
\end{framed}

\subsection{Probabilistic model of CPD}
\label{sec:prob_model_CPD}
We   introduce the mixture model
\begin{equation}
\label{eq:CPD_probabilistic_model}
Y = \sum_{j=1}^N \mathbbm{1}(Z=j) Y^{(j)} +
 \mathbbm{1}(Z=N+1) Y^{(N+1)},
\end{equation}
where $Y^{(j)} \sim \mathcal{N}(\mu_j, \sigma^2 \mathbbm{1} )$ for $j=1,\ldots,N$,
$Y^{(N+1)} \sim {\rm Uniform}(\mathcal{D}),$ 
and $Z\sim {\rm Multinomial}(\{1,\ldots,N+1\})$ are independent random variables. Given $w\in [0,1)$, we set
$$
{\rm Pr}(Z=j )
= 
\left\{
\begin{array}{ll}
\displaystyle{\frac{1-w}{N}} & j=1,\ldots,N \\
\displaystyle{w} & j=N+1 \\
\end{array}
\right.
\quad
|\mathcal{D}|=Q.
$$
Exploiting the previous hypotheses, we find
$$
{\rm Pr}(Y = y, Z=j )
=
{\rm Pr}(Y^{(j)} = y)  \,  {\rm Pr}(Z=j )
=
\left\{
\begin{array}{ll}
\displaystyle{
\frac{1-w}{N} \; \frac{1}{(2\pi \sigma^2)^{d/2}} \;
{\rm exp} \left(
-\frac{\|y-\mu_j\|_2^2}{2\sigma^2}
\right)
} &
{\rm if} \; 
j\in \{1,\ldots,N\} \\[3mm]
\displaystyle{\frac{w}{Q}
} &
{\rm if} \; 
j = N+1.\\
\end{array}
\right.
$$
We consider the set of tunable parameters
$\boldsymbol{\Theta} = \{ \mu_1,\ldots, \mu_N, \sigma^2 \}$, while we fix $w, \mathcal{D}$  \emph{a priori}. We  observe that \eqref{eq:CPD_probabilistic_model} reads as a  GMM with a  perturbance given by the uniform random variable $Y^{(N+1)}$: the latter is intended to account for noise and outliers in the dataset. We   do not explicitly construct the domain $\mathcal{D}$: we simply assume that all the observed datapoints belong to 
$\mathcal{D}$.

Given $j\in \{1,\ldots,N+1\}$ and $y\in \mathcal{D}$, using \eqref{eq:posterior_Z}, we find
\begin{equation} 
\label{eq:CPD_expressionP}
{\rm Pr}\left(
Z=j | Y=y
\right)
=
\frac{{\rm Pr}(Y^{(j)} = y, Z=j )}{  \sum_{k=1}^{N+1}  {\rm Pr}(Z=k ) {\rm Pr}(Y^{(k)} = y)   
}
=
\frac{ 
{\rm Pr}(Y^{(j)} = y)   \, 
{\rm Pr}(Z=j )    }{  \sum_{k=1}^{N+1}  {\rm Pr}(Z=k ) {\rm Pr}(Y^{(k)} = y)   
}
=
\frac{  {\rm exp}\left( -\frac{\| y - \mu_j \|_2^2}{2\sigma^2}  \right)       }{    \sum_{k=1}^N 
 {\rm exp}\left( -\frac{\| y - \mu_k \|_2^2}{2\sigma^2}  \right)
+ 
c },
\end{equation}
with $c=\frac{w}{1-w}  \frac{N}{Q}  (2\pi \sigma^2)^{d/2}$, for $j=1,\ldots,N$.
We further introduce the matrix 
$\mathbf{P}\in \mathbb{R}^{Q\times N}$ such that
$( \mathbf{P}  )_{i,j} = {\rm Pr}\left(
Z=j | Y=y_i \right)$
 for $i=1,\ldots,Q$ and  $j=1,\ldots,N$ and we define
$Q_{\rm p} = \sum_{i,j} ( \mathbf{P}  )_{i,j} $;
then,   
assuming that $y_1,\ldots,y_Q\in \mathcal{D}$, 
we obtain
\begin{equation}
\label{eq:loglikelihood_expression}
\begin{array}{rl}
-\mathfrak{Q}(\boldsymbol{\Theta}; \mathbf{P})
=
&
\displaystyle{
- \sum_{i=1}^Q
\left(
\left(
\sum_{j=1}^N
(\mathbf{P})_{i,j}
\log {\rm Pr} \left(
Y=y_i, Z = j | \boldsymbol{\Theta}
\right)
\right)
+
\left(1 -   \sum_{k=1}^N
(\mathbf{P})_{i,k} \right)
\log (1/|\mathcal{D}|)
\right)
}
\\[3mm]
=
&
\displaystyle{
 \sum_{i=1}^Q
\sum_{j=1}^N
(\mathbf{P})_{i,j}
\left(
\frac{1}{2\sigma^2}
\| y_i - \mu_j  \|_2^2
+\frac{d}{2} \log(\sigma^2) +
 \log\left(
\frac{1-w}{N (2\pi)^{d/2}}
\right)
\right)
+C_1;
}
\\[3mm]
=
&
\displaystyle{
\left(
\frac{1}{2\sigma^2}
 \sum_{i=1}^Q
\sum_{j=1}^N
(\mathbf{P})_{i,j}
\| y_i - \mu_j  \|_2^2
\right)
+\frac{d Q_{\rm p}}{2} \log(\sigma^2) 
+C_2,
}
\\
\end{array}
\end{equation}
where $C_1,C_2$ are constants that 
 are independent of $\boldsymbol{\Theta}$.

So far, we have not forced the GMM centroids to move coherently; towards this end, we propose the model
$\mu_j = x_j + v(x_j)$ where $\{x_j\}_{j=1}^N$ is the reference point cloud that we wish to deform and $v:\mathbb{R}^d\to \mathbb{R}^d$ is a displacement field that is assumed to belong to the native Reproducing Kernel Hilbert space (RKHS) $\mathcal{H}_{\phi}$ 
associated with the radial basis function (RBF)
$\phi: r \mapsto  {\exp}\left( -\frac{r^2}{2 \beta^2} \right)$. Then, we introduce the objective function 
$\mathfrak{E}$
for the M-step by adding  the regularization $\frac{\lambda}{2} \|v\|_{\mathcal{H}_{\phi}}^2$ to \eqref{eq:loglikelihood_expression}:
\begin{equation}
\label{eq:regularized_loglikelihood}
\mathfrak{E}(v, \sigma^2| \mathbf{P})
\, =\,
\frac{1}{2\sigma^2}
 \sum_{i=1}^Q
\sum_{j=1}^N
(\mathbf{P})_{i,j}
\| y_i - x_j - v(x_j)  \|_2^2
\;
+ \; \frac{d Q_{\rm p}}{2} \log(\sigma^2) 
\; + \; 
\frac{\lambda}{2} \|v\|_{\mathcal{H}_{\phi}}^2.
\end{equation}
Exploiting the representation theorem for RKHS
 (e.g.,  \cite{wendland2004scattered}), we find that minimizers of \eqref{eq:regularized_loglikelihood} belong to the space $\mathcal{U}_{\rm hf}^{\rm cpd}$ introduced in  \eqref{eq:CPD_affine_space}; we also observe that in the Bayesian setting the regularization 
 $\frac{\lambda}{2} \|v\|_{\mathcal{H}_{\phi}}^2$ might be associated to the logarithm of a  prior on the displacement field.

We observe that the probabilistic model introduced in this section  depends on three free parameters:
$w,\beta,\lambda$.  The parameter
$w\in [0,1]$ is designed to account for noise and outliers in the datasets: it should be set to zero for noiseless data.
The parameter $\beta$ reflects the strength of interaction between points: the value of $\beta$ should thus depend on the characteristic length-scale of the displacement field we wish to approximate. The value of $\lambda$ is associated with  the Tikhonov regularization of the $E$-step and thus reflects the trade-off between data fitting and smoothness regularization.
Finally, the algorithm also depends on the choice of the RBF $\phi$: following \cite{myronenko2010point}, we here consider the Gaussian kernel; however,we observe that several other choices are possible and there is an extensive literature on the comparison between different RBFs
\cite{rocha2009selection}.

\subsection{Coherent point drift procedure}
We adapt the EM algorithm \ref{alg:EM} to the mixture model \eqref{eq:CPD_probabilistic_model}: we present the detailed procedure in Algorithm \ref{alg:CPD}.  To facilitate the M-step, we first solve for the displacement $v$ and then for $\sigma^2$.

\begin{framed}
\captionof{algorithm}{Coherent point drift \cite{myronenko2010point}. \label{alg:CPD}}

\begin{flushleft}
\emph{Inputs:} $\mathbf{X} = [x_1,\ldots,x_N]^{\star}  \in \mathbb{R}^{N\times d}$,
$\mathbf{Y}^{\rm raw}  = [y_1^{\rm raw} ,\ldots,y_Q^{\rm raw} ]^{\star}  \in \mathbb{R}^{Q\times d}$. 
\smallskip

\emph{Hyper-parameters:} $w\in [0,1), \beta,\lambda>0$.
\smallskip

\emph{Outputs:} 
$\mathbf{Y} = [x_1+v(x_1),\ldots,x_N+v(x_N)]^{\star}  \in \mathbb{R}^{N\times d}$.
\end{flushleft}                      

\begin{algorithmic}[1]
\State
Set $v^{(0)} = 0$ and
$\sigma^{2,(0)}
=\frac{1}{d QN } \sum_{i=1}^Q \sum_{j=1}^N \| x_j - y_i^{\rm raw}   \|_2^2
$.
\medskip
			
\For {$k=1,2,   \ldots, \texttt{until  convergence}$ }
\medskip		
 
 \State
 \textbf{E-step:} Find  
 $\mathbf{P}^{(k)}\in \mathbb{R}^{Q\times N}$ such that
$( \mathbf{P}^{(k)}\  )_{i,j} = {\rm Pr}\left(
Z=j | Y=y_i^{\rm raw}  \right)$ using  
\eqref{eq:CPD_expressionP}  with $\mu_j=x_j$

$+ v^{(k-1)} (x_j)$ and $\sigma^2=\sigma^{2,(k-1)}$.
 \medskip
 
  \State
 \textbf{M-step (I):} find
$
v^{(k)}  = {\rm arg} \min_{v\in \mathcal{H}_{\phi}} 
\mathfrak{E}(v, \sigma^{2,(k-1)} )| \mathbf{P}^{(k)})
$ 
 (cf. \eqref{eq:regularized_loglikelihood}).
  \medskip
 
  \State
 \textbf{M-step (II):} find
$
\sigma^{2,(k)}   = {\rm arg} \min_{\sigma^2\in  \mathbb{R}_+} 
\mathfrak{E}(v^{(k)},  \sigma^2| \mathbf{P}^{(k)})
$.
 \medskip
 
 \EndFor
\end{algorithmic}

\end{framed}

Exploiting the representation theorem for $v$, we find that 
$v^{(k)}(\cdot)= \sum_{j=1}^N w_j^{(k)} \phi(\|\cdot - x_i \|_2)$ and
$ \mathbf{W}^{(k)} = [w_1^{(k)}, \ldots,w_N^{(k)}  ]^{\star}  \in \mathbb{R}^{N\times d}$ solves the linear system (we omit dependence on the iteration count $k$ at the right-hand-side to shorten notation)
\begin{equation}
\label{eq:Mstep_part1}
\mathbf{W}^{(k)}  = \left(\mathbf{G} + \lambda \sigma^2   \left( \mathbf{d}(  \mathbf{P} \mathbf{1}_Q     \right)^{-1}  \right) ^{-1}
 \left(  
\left( \mathbf{d}(  \mathbf{P} \mathbf{1}_Q ) \right)^{-1}   
 \mathbf{P}  \mathbf{Y}^{\rm raw}  -  \mathbf{X}
 \right),  
\end{equation}
where 
$\mathbf{G} \in \mathbb{R}^{N\times N}$ satisfies
$\left( \mathbf{G} \right)_{i,j}
={\rm exp} \left(
-
\frac{1}{2\beta^2} \|x_i - x_j \|_2^2
\right)$, $i,j=1,\ldots,N$,
$\mathbf{1}_Q$ is the $Q$-dimensional vector with entries all equal to one, and
$\mathbf{d} \left(  \mathbf{P} \mathbf{1}_Q     \right) = {\rm diag} \left(
( \mathbf{P} \mathbf{1}_Q)_1,\ldots,
( \mathbf{P} \mathbf{1}_Q)_N 
\right) \in \mathbb{R}^{N\times N}$.
By tedious but straightforward calculations, we also find the explicit expression for $\sigma^{2,(k)}$
\begin{equation}
\label{eq:Mstep_part2}
\sigma^{2,(k)} = 
\frac{1}{Q_{\rm p} d}
\left(
{\rm Tr} \left(
( \mathbf{Y}^{\rm raw} )^{\star}
\mathbf{d}(  \mathbf{P}^{\star} \mathbf{1}_N )
\mathbf{Y}^{\rm raw} 
\right)
\, -\,
2 \, 
{\rm Tr} \left(
\left( \mathbf{P} \mathbf{Y}^{\rm raw}  \right)^{\star}
\mathbf{Y}
\right)
\, + \,
{\rm Tr} \left(
 \mathbf{Y}^{\star} 
\mathbf{d}(  \mathbf{P}^{\star} \mathbf{1}_N )
\mathbf{Y}
\right)
\right).
\end{equation}
with $Q_{\rm p} = \sum_{i,j}  \left( \mathbf{P}^{(k)}  \right)_{i,j}$ and 
$\mathbf{Y}   = \mathbf{X} + \mathbf{G} \mathbf{W}^{(k)} $.
 
Since the seminal work   \cite{myronenko2010point}, several authors have discussed how to effectively implement Algorithm \ref{alg:CPD} to enable its application to large datasets and also to prevent unstable behaviors.
In this work, we resort to an eigenvalue decomposition of the matrix $\mathbf{G}$ and to the Woodbury identity to invert the system in \eqref{eq:Mstep_part1} (cf.  
\cite[section 6]{myronenko2010point}). Furthermore, before solving the system, we remove all lines $i$ for which  $\sum_{j=1}^Q (\mathbf{P})_{i,j} < 10^{-12}$ and we set $\mathbf{W}(i,:)=0$. In all our numerical experiments, we consider the termination condition
\begin{equation}
\label{eq:CDP_termination}
\| \mathbf{W}^{(k)} -\mathbf{W}^{(k-1)}   \|_2 < 10^{-4} \;\; 
\texttt{OR} \;\;
\max_{j=1,\ldots,Q} \,
\min_{i=1,\ldots,N}
\| \mathbf{Y}^{\rm raw}(j,:) -\mathbf{Y}(i,:)   \|_2 < 10^{-5}.
\end{equation}

\subsection{Low-dimensional representation of the displacement field $v$}
As discussed in section \ref{sec:dimensionality_reduction}, and also shown in the numerical experiment of section \ref{sec:TT_sisc2020_revised}, we might wish to restrict computations over a low-dimensional subspace
$\mathcal{U}_M^{\rm cpd}$.
Since we are ultimately interested in bounded-domain registration, the natural choice is $\mathcal{U}_M^{\rm cpd} = \mathcal{U}_M$. Since in general $\mathcal{U}_M$ is not contained in $\mathcal{H}_{\phi}$, we propose to replace the regularization $\|\cdot\|_{\mathcal{H}_{\phi}}^2$ 
in \eqref{eq:regularized_loglikelihood}
with $\|\cdot \|_{H^2(\Omega_{\rm box})}^2$. In conclusion, we replace Line 4 of Algorithm \ref{alg:CPD}  with 
\begin{equation}
\label{eq:reducedCPD}
v^{(k+1)}
=
 {\rm arg} \min_{v\in \mathcal{U}_M} 
\;
\frac{1}{2\sigma^2}
 \sum_{i=1}^Q
\sum_{j=1}^N
(\mathbf{P})_{i,j}
\| y_i - x_j - v(x_j)  \|_2^2
\;
+ \;
\frac{\lambda}{2} \|v\|_{H^2(\Omega_{\rm box})}^2.
\end{equation}
Provided that $\{\psi_m\}_{m=1}^M$ is an orthonormal basis of $\mathcal{U}_M$, by tedious but straightforward calculations, we find that 
$v^{(k+1)}(\cdot) = \sum_{m=1}^M  \left( \mathbf{A}^{-1} \mathbf{b} \right)_m \psi_m(\cdot)$ where 
$$
\begin{array}{l}
\displaystyle{
\left( \mathbf{A} \right)_{m,m'}
=
\lambda\sigma^{2} \,  \delta_{m,m'} \, + \, 
\sum_{j=1}^N \sum_{\ell=1}^d 
\left( \psi_m(x_j)  \right)_{\ell}
\left( \mathbf{P} \mathbf{1}_Q \right)_{j}
\left( \psi_{m'}(x_j)  \right)_{\ell}
}
\\[3mm]
\displaystyle{
\left( \mathbf{b} \right)_{m}
=
\sum_{j=1}^N \sum_{\ell=1}^d 
\left( \psi_m(x_j)  \right)_{\ell}
\left( \mathbf{P} \mathbf{Y}^{\rm raw}
- 
\mathbf{d}(  \mathbf{P} \mathbf{1}_Q  ) \mathbf{X}
 \right)_{j,\ell},
}
\\
\end{array}
$$
for $m,m'=1,\ldots,M$.

\section{Proofs}
\label{sec:proofs}

\subsection{Proposition  \ref{th:quasi_hausdorff_bound}}
We first observe that if $\{ x_i \}_{i=1}^N$ is an $\epsilon$-cover of $\partial \Omega$ then $\{ \Phi(x_i) \}_{i=1}^N$ is 
$K \epsilon$-cover of $\partial \Phi(\Omega)$.
To prove this statement, consider $\tilde{x} \in \partial \Phi(\Omega)$; since $\Phi$ is bijective, $\partial \Phi(\Omega) = \Phi(\partial \Omega)$ and thus  there exists $x\in \partial \Omega$ such that
$\tilde{x} = \Phi(x)$. We denote by $i^{\star}\in \{1,\ldots,N\}$ the index that satisfies $\| x - x_{i^{\star}}\|_2 = {\rm dist}\left(
x, \{ x_i \}_{i=1}^N \right)$: since $\{ x_i \}_{i=1}^N$ is an $\epsilon$-cover of $\partial \Omega$,  we have 
$\| x - x_{i^{\star}}\|_2 \leq \epsilon$. Then, we find
$$
{\rm dist}\left( \tilde{x}, \{ \Phi(x_i) \}_{i=1}^N \right)
=
{\rm dist}\left(\Phi(x), \{ \Phi(x_i) \}_{i=1}^N \right)
\leq
\| \Phi(x) - \Phi(x_{i^{\star}})  \|_2
\leq K \epsilon.
$$

Exploiting the previous estimate, we obtain that for all 
$\tilde{x} = \Phi(x) \in  \partial \Phi(\Omega)$ we have
$$
\begin{array}{rl}
\displaystyle{
{\rm dist}\left( \tilde{x},  \partial V   \right)
=
} &
\displaystyle{
\inf_{ y \in   \partial V }
\|\Phi(x) - y \|_2
\leq
\min_{ i=1,\ldots,N}
\|\Phi(x) - y_i \|_2
}
\\[3mm]
\displaystyle{
\leq
} &
\displaystyle{
\min_{ i = 1,\ldots,N  }
\left(
\|\Phi(x_i) - y_i \|_2
+
\|\Phi(x_i) - \Phi(x) \|_2
\right)
}
\\[3mm]
\leq
&
\displaystyle{
\left(
\max_{ i = 1,\ldots,N  }
\|\Phi(x_i) - y_i \|_2
\right)
+
\left(
\min_{ i = 1,\ldots,N  }
\|\Phi(x_i) - \Phi(x) \|_2
\right)
\leq
\max_{ i = 1,\ldots,N  }
\|\Phi(x_i) - y_i \|_2
+ K\epsilon.
\quad
\qed
}
\\
\end{array}
$$

\subsection{Lemma \ref{th:neigh_prel}}
Proof of the second statement is an immediate consequence of the Weyl's tube formula \cite{gray2003tubes}, which states that if  the curves (or surfaces for $d=3$) $\{x\pm \delta \mathbf{n}(x) \,: \, x\in \partial U \}$ do not have self-intersections, we have
$$
{\rm Neigh}_{\delta}(\partial U) = 2\delta | \partial U|_{(d-1)}
+
\frac{4\pi}{3} \chi(\partial U) \delta^3,
$$
where $\chi(\partial U) $ is the so-called Euler characteristics. The quantity
$\chi(\partial U)$ is invariant under isomorphisms and is equal to $0$ for the unit circle and to $2$ for the unit sphere.

We now prove the first statement. In more detail, we show that  any $y\in {\rm Neigh}_{\delta}^{\rm t}( \partial U  )$ belongs to $ {\rm Neigh}_{\delta}( \partial U  )$ and then that
any $y\in {\rm Neigh}_{\delta}( \partial U  )$ belongs to $ {\rm Neigh}_{\delta}^{\rm t}( \partial U  )$.

Let $y= x + t \mathbf{n}(x)\in {\rm Neigh}_{\delta}^{\rm t}( \partial U  ) $. Clearly, we have
$
{\rm dist}(y, \partial U)
\leq
\| y -  x \|_2
=
|t| < \delta$, 
which implies that $y  \in 
{\rm Neigh}_{\delta}( \partial U)$.

  Let $y \in   {\rm Neigh}_{\delta}( \partial U ) $.
We define   $f:\partial U\to \mathbb{R}$  such that
  $f(x) =  \| x - y  \|_2^2$; we denote by $\nabla_{ \partial U}  f := \nabla f - \mathbf{n} ( \mathbf{n} \cdot \nabla f)$   the surface gradient of $f$.
   Since $ \partial U$ is a closed set in $\mathbb{R}^d$, there exists a minimizer $\widetilde{x} \in \partial U$ that  
 $f(\widetilde{x} ) = \inf_{x\in \partial U} f(x)$; furthermore, 
$$
 \| \widetilde{x}  - y  \|_2^2 = \inf_{x\in \partial U} f(x) =\left(  {\rm dist}(y,\partial U ) \right)^2 < \delta^2
 \quad \Rightarrow
 \|  \widetilde{x} - y  \|_2 < \delta.
$$ 
Since $ \partial U$ is a closed hyper-surface  without boundary,
 $\widetilde{x}$ must be a stationary point of  $f$ over $\partial U$, that is $\nabla_{\partial U} f( \widetilde{x}   ) = 0$. Exploiting the definition of the surface gradient, we find that
$\widetilde{x}   - y =   \widetilde{\mathbf{n}}
 (  \widetilde{\mathbf{n}} \cdot (x^{\star}  - y)) $
with  $ \widetilde{\mathbf{n}} = \mathbf{n}(\widetilde{x} )$, which implies that 
$\widetilde{x}  - y$ is parallel to $ \widetilde{\mathbf{n}}$ and then that
$$
y = \widetilde{x}  - \left(  \theta  \|  \widetilde{x}  - y \|_2 \right)  \widetilde{\mathbf{n}} ,
$$
for either $\theta=1$ or $\theta=-1$. Recalling \eqref{eq:tubular_neigh}, we find that $y\in  {\rm Neigh}_{\delta}^{\rm t}( \partial U)$.  \qed
  
\subsection{Lemma \ref{th:equivalent_condition_deltareg}}
The no-intersection hypothesis implies that 
${\rm Neigh}_{\delta}(\partial U)$ is diffeomorphic to a spherical shell. Therefore, there exists a mapping $\Phi$ that transforms $U$ into the unit ball and 
${\rm Neigh}_{\delta}(\partial U)$ into the shell 
$\mathcal{B}_{1+\delta}(0) \setminus \mathcal{B}_{1-\delta}(0)$. We conclude that it suffices to prove the result for $U$ equal to the unit ball.

By contradiction, let 
$y \in (V \setminus  {\rm Neigh}_{\delta}(\partial U) ) \cap \mathcal{B}_{1-\delta}(0)$ and consider $p = x - \delta \mathbf{n}(x) \in U \setminus \overline{V}$; note that
$p\in \partial \mathcal{B}_{1-\delta}(0)$. Clearly, the segment
$\Gamma = \{t p + (1-t) y: t\in (0,1) \}$  is contained in
$ \mathcal{B}_{1-\delta}(0)$ and --- since $y\in V$ and $p\notin V$ ---   intersects $\partial V$ at some point
$z\in \partial V \cap \Gamma$. 
We thus proved the existence of a point $z\in \partial V$ that does not belong to  $ {\rm Neigh}_{\delta}(\partial U) )$, which is a contradiction.
The case $y \in (V \setminus  {\rm Neigh}_{\delta}(\partial U) ) \setminus \mathcal{B}_{1+\delta}(0)$ can be treated similarly. \qed

\subsection{Lemma \ref{th:lemma_smooth}.}

The aim of this section is to further discuss Definition
 \ref{def:delta_regularity} and to prove a refined version of Lemma \ref{th:lemma_smooth}.
We focus on the three-dimensional case: the same argument applies   to two-dimensional domains  with minor modifications.
 
We first recall the following elementary result.

\begin{lemma}
\label{th:mean_value_theorem_vector}
Let $f:[a,b]\to \mathbb{R}^k$ be differentiable.
Then, $\|f(b) -f(a) \|_2 \leq (b-a) \max_{\zeta\in [a,b]} \|f(\zeta)\|_2$.
\end{lemma}
\begin{proof}
It suffices to apply the mean value theorem to the scalar function
$\phi(\zeta) = (f(b) - f(a))\cdot f(\zeta)$.
\end{proof}

Next Lemma introduces relevant quantities for the subsequent result.  

\begin{lemma}
\label{th:weierstrass_kindof}
Let $U$ be a three-dimensional bounded domain of class $C^2$. Then, there exist positive constants $L, L_{\rm inv}, r$  such that for any $x\in \partial U$ there exists a bijective mapping $\Phi_x: A \to    U \cap \mathcal{B}_r(x)$ 
with $A = \mathcal{B}_r(0) \cap \{x':  (x')_2<0\}$
with the following properties:
(i) $\Phi_x$ has  Lipschitz constant $L$,
(ii) the inverse of $\Phi_x$ has  Lipschitz constant $L_{\rm inv}$,
(iii) $s\in  \mathcal{C}_r(0)  \mapsto \Phi_x([s,0])$ is a regular parameterization of  $\partial U \cap \mathcal{B}_r(x)$ where
$\mathcal{C}_r(0)$ denotes the two-dimensional ball of radius $r$ centered in the origin. Furthermore, 
there exists $\kappa_{\rm max}>0$  such that  the  principal curvatures are uniformly bounded   in absolute value  by $\kappa_{\rm max}$.
\end{lemma}
\begin{proof}
If $U$ is of class $C^2$, it is in particular a bi-Lipschitz domain. This implies that for any $x\in \partial U$ there exists $R=R(x)$ and a bijection 
 $\Phi_x: A \to    U \cap \mathcal{B}_R(x)$ 
with $A = \mathcal{B}_R(0) \cap \{x':  (x')_2<0\}$
with the following properties:
(i) $\Phi_x$ has  Lipschitz constant $L$,
(ii) the inverse of $\Phi_x$ has  Lipschitz constant $L_{\rm inv}$,
(iii) $s\in  \mathcal{C}_R(0)  \mapsto \Phi_x([s,0])$ is a regular parameterization of  $\partial U \cap \mathcal{B}_r(x)$.
We denote by $R_{\rm max}:\partial U \to \mathbb{R}_+$ the maximum radius for which the above conditions hold: it is easy to verify that 
$R_{\rm max}$ is a continuous function; therefore, since 
$\partial U$ is compact, 
exploiting the Weierstrass theorem, we find that 
$R_{\rm max}$ attains a minimum at some 
$x\in \partial U$; we thus find 
$r:=R_{\rm max}(x)>0$.

If $U$ is of class $C^2$, the principal curvatures are continuous over $\partial U$ and thus uniformly bounded.
\end{proof}

Then, Lemma \ref{th:lemma_smooth_refined} provides a refined version of Lemma  \ref{th:lemma_smooth}. Note that the bound on $\delta_0$ depends on the global property 
$r$  of $\partial U$, and on the local properties
$\kappa_{\rm max}, L,L_{\rm inv}$.
The proof involves the introduction of the so-called shape operator 
(see, e.g.,  \cite[Chapter 5]{o2006elementary}).
Given $x\in \partial U$, we denote by $T_x$ the tangent space at $x$ and we define the operator $S[x]$ such that
$$
S[x] \mathbf{t} = - \lim_{\epsilon \to 0} \; \frac{\mathbf{n}(x+\epsilon \mathbf{t}) - \mathbf{n}(x)     }{\epsilon},
\quad
\forall \; \mathbf{t} \in T_x.
$$
It is possible to prove that $S[x]$ is a linear symmetric operator from $T_x$ to $T_x$ whose eigenvalues are the principal curvatures of $\partial U$ at $x$ $\kappa_1(x),\kappa_2(x)$. Therefore, we find that the dual norm of $S[x]$ satisfies 
\begin{equation}
\label{eq:onemore_formula}
\| S[x]   \|_{\star}
\, :=\,
\sup_{\mathbf{t} \in T_x }
\frac{\| S[x] \mathbf{t}  \|_2}{\| \mathbf{t}  \|_2} 
 =
 \max\left\{
 |\kappa_1(x)   |,
 \;
 |\kappa_2(x) |
  \right\}.
 \end{equation}

\begin{lemma}
\label{th:lemma_smooth_refined}
Let $L,L_{\rm inv}, r,\kappa_{\rm max}>0$ be the quantities introduced in Lemma  \ref{th:weierstrass_kindof}.
Then, we find that $U$ is $\delta$-regular for $\delta< \min \{ \frac{r}{2}, \frac{1}{L L_{\rm inv} \kappa_{\rm max}} \}$.
\end{lemma}

\begin{proof}
Let $x,y\in \partial U$ satisfy
\begin{equation}
\label{eq:contradiction_3D}
x -  \delta \mathbf{n}(x) = y -  \delta \mathbf{n}(y).
\end{equation}
Below, we prove that if $\delta< r/2$ then \eqref{eq:contradiction_3D} implies that $\delta \geq  \frac{1}{L L_{\rm inv} \kappa_{\rm max}}$. 
Recalling Lemma \ref{th:equivalent_condition_deltareg}, 
the latter implies that $U$ is guaranteed to be $\delta$-regular if  $\delta<r/2$ and $\delta< \frac{1}{L L_{\rm inv} \kappa_{\rm max}}$.

If $\delta<r/2$, we have $\| x - y \|_2 \leq \delta \| \mathbf{n}(x)  - \mathbf{n}(y)  \|_2 \leq 2 \delta< r$; therefore, given
$A = \mathcal{B}_r(0) \cap \{x':  (x')_3<0\}$, 
there exists  a mapping $\Phi: A \to  U \cap \mathcal{B}_r(x)$
that satisfies the conditions of Lemma \ref{th:weierstrass_kindof}.
 Then, if we denote by $t,s\in \mathbb{R}^2$  the points satisfying $x=\Phi([t,0]), y=\Phi([s,0])$ such that $\| s \|_2,  \|t\|_2 < r$, we find
$$
\|t - s\|_2 
\overset{\rm (i)}{\leq}
 L_{\rm inv} \|x-y\|_2
\overset{\rm (ii)}{=}
\left( 
 L_{\rm inv}
\| \mathbf{n}(x) - \mathbf{n}(y)  \|_2    \right) \delta.
$$
Note that
 in (i) we used the fact that $\Phi^{-1}$ is Lipschitz,  while in (ii) we used \eqref{eq:contradiction_3D}.
 
Let $\gamma: \xi \mapsto  [(1-\xi) t + \xi s,0] $ be a path from $t$ to $s$ in $A$ and define 
 $\mathbf{N}:\xi \mapsto \mathbf{n}( \Phi(\gamma(\xi)) )$. Clearly, we have
 $\mathbf{N}(0) = \mathbf{n}(x)$,  $\mathbf{N}(1) = \mathbf{n}(y)$, and we also have that
$\nabla \Phi(  \gamma(\xi)   ) \gamma'(\xi))$ 
 belongs the tangent space $T_x$ with 
 $x=  \Phi(  \gamma(\xi)   )$. 
 . Then, applying Lemma \ref{th:mean_value_theorem_vector}, the chain rule, and the definition in Eq. \eqref{eq:onemore_formula}, 
  we obtain
 $$
 \| \mathbf{n}(x) - \mathbf{n}(y)  \|_2
 \leq
 \max_{\xi \in [0,1]} \| \mathbf{N}'(\xi)  \|
 =
 \max_{\xi\in [0,1]}
 \| S\left[
 \Phi(\gamma(\xi))
\right]    \|_{\star} \|  \nabla \Phi( \gamma(\xi))   \|_2 \| t - s \|_2 
 \leq\kappa_{\rm max} L \|t - s \|_2.
 $$
 By combining the latter two equations, we obtain
 $$
 \|t - s\|_2 
\leq 
L 
 L_{\rm inv} 
 \kappa_{\rm max}   \|t - s \|_2 \delta.
 $$
 By dividing  both sides by  $ \|t - s\|_2$, we obtain
 $\delta\geq (  L 
 L_{\rm inv} 
 \kappa_{\rm max}     )^{-1}$.
\end{proof}

\subsection{Lemma \ref{th:geo_lemma}}
We first remark that 
\begin{equation}
\label{eq:easy_estimate_geo}
{\rm dist}_{\rm H}(U,V) \leq
{\rm dist}_{\rm H}(\partial U, \partial V) 
\end{equation}
We have indeed that
$
\sup_{{x} \in V} {\rm dist}({x} , U)
=
\sup_{{x} \in V\setminus U} {\rm dist}({x}, U)
=
\sup_{{x} \in \partial V\setminus U} {\rm dist}({x}, U)
=
\sup_{{x} \in \partial V\setminus U} {\rm dist}({x}, \partial U)
\leq
{\rm dist}_{\rm H}(\partial U, \partial V) .
$
Similarly, we verify that  $\sup_{{x} \in U} {\rm dist}({x}, V) \leq {\rm dist}_{\rm H}(\partial U, \partial V)$.

Let ${x} \in \partial V$. By contradiction, assume that
${\rm dist} ({x}, \partial U) > \delta$, that is
$\mathcal{B}_{\delta}({x}) \cap \partial U = \emptyset$, which implies 
$\partial U \subset  {\rm Neigh}_{\delta}(\partial V) \setminus \mathcal{B}_{\delta}({x})$. Since $V$ is $\delta$-regular, we must have that $U \subset {\rm Neigh}_{\delta}(\partial V)$ and 
thus 
$| U  |_{(d)} \leq   {\rm Neigh}_{\delta}(\partial V)$, which contradicts Hypothesis (ii). 

In conclusion, we find
$$
{\rm dist}_{\rm H}(U,V)
\overset{\eqref{eq:easy_estimate_geo} }{\leq}
{\rm dist}_{\rm H}(\partial U, \partial V)
=\max \Big\{
\underbrace{\max_{{x} \in \partial V} {\rm dist}({x}, \partial U)}_{\leq\delta}
, 
\underbrace{
\max_{{x} \in \partial U} {\rm dist}({x}, \partial V)}_{\leq \delta}
\Big\}
\leq \delta. 
\qquad
\qed
$$

\section{Further investigations for the three-point registration problem}
\label{sec:further_tests}

We assess performance of the elasticity-based objectives introduced in Remark \ref{remark:link_elasticity} for the model problem in section \ref{sec:three_point_registration}. We consider $t=0.7$, $\delta=10^{-6}$, $n_{\rm lp}=25$; we set $\lambda_1=\frac{\nu}{(1-2\nu)(1+\nu)}$ and 
$\lambda_2=\frac{E}{1+\nu}$ with $E=1$ and $\nu=0.3$ in \eqref{eq:isotropic_linear_strain} and \eqref{eq:neohookean_strain}. 
Note that the choice of $\lambda_1,\lambda_2$ corresponds to the plane strain assumption; note also that, since \eqref{eq:isotropic_linear_strain} and \eqref{eq:neohookean_strain} are linear with respect to the Young's modulus $E$, results are independent of its choice.

\begin{figure}[h!]
\centering
 \subfloat[linear elasticity strain] 
{  \includegraphics[width=0.33\textwidth]
 {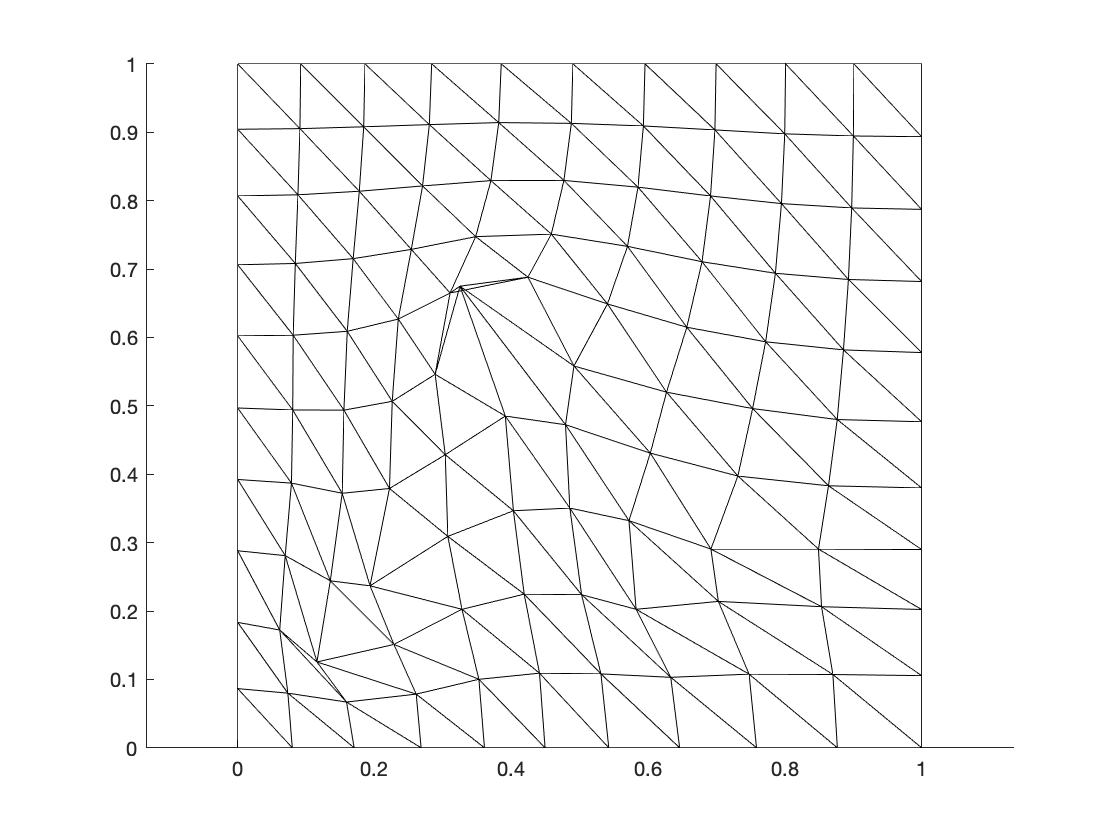}}
   ~~
 \subfloat[neohookean strain] 
{  \includegraphics[width=0.33\textwidth]
 {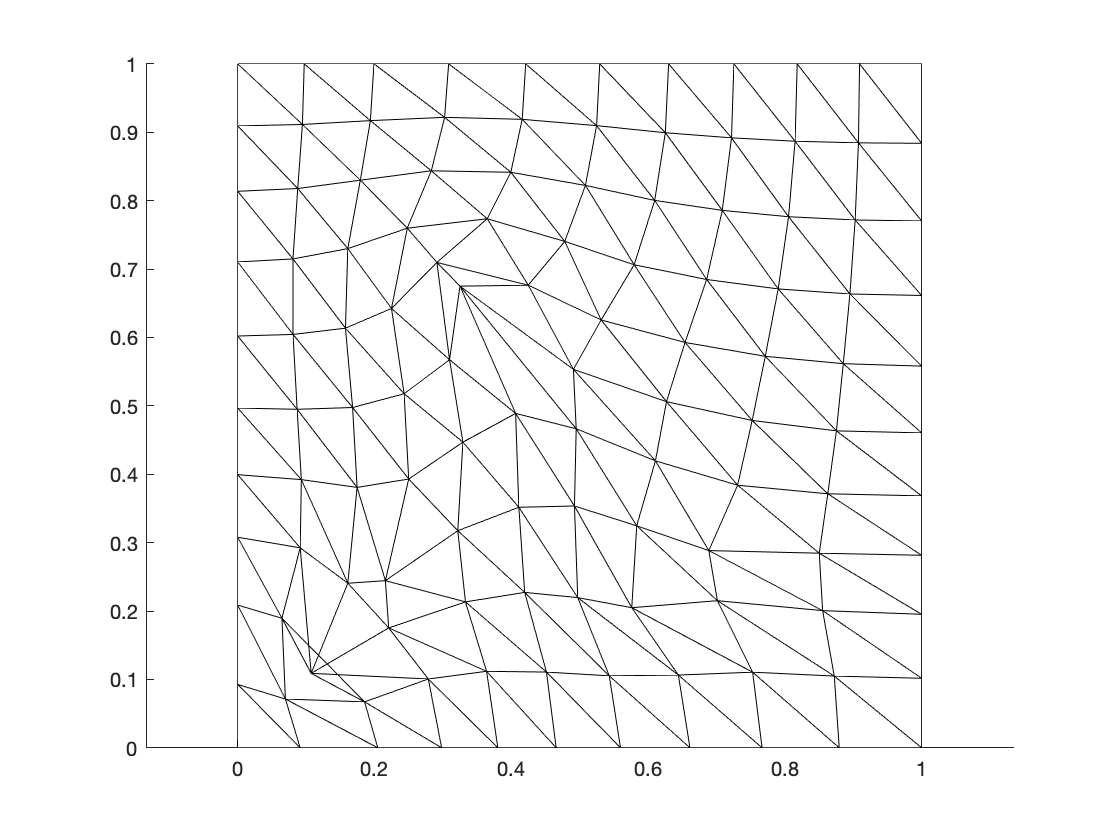}}
~~
 \subfloat[H2] 
{  \includegraphics[width=0.33\textwidth]
 {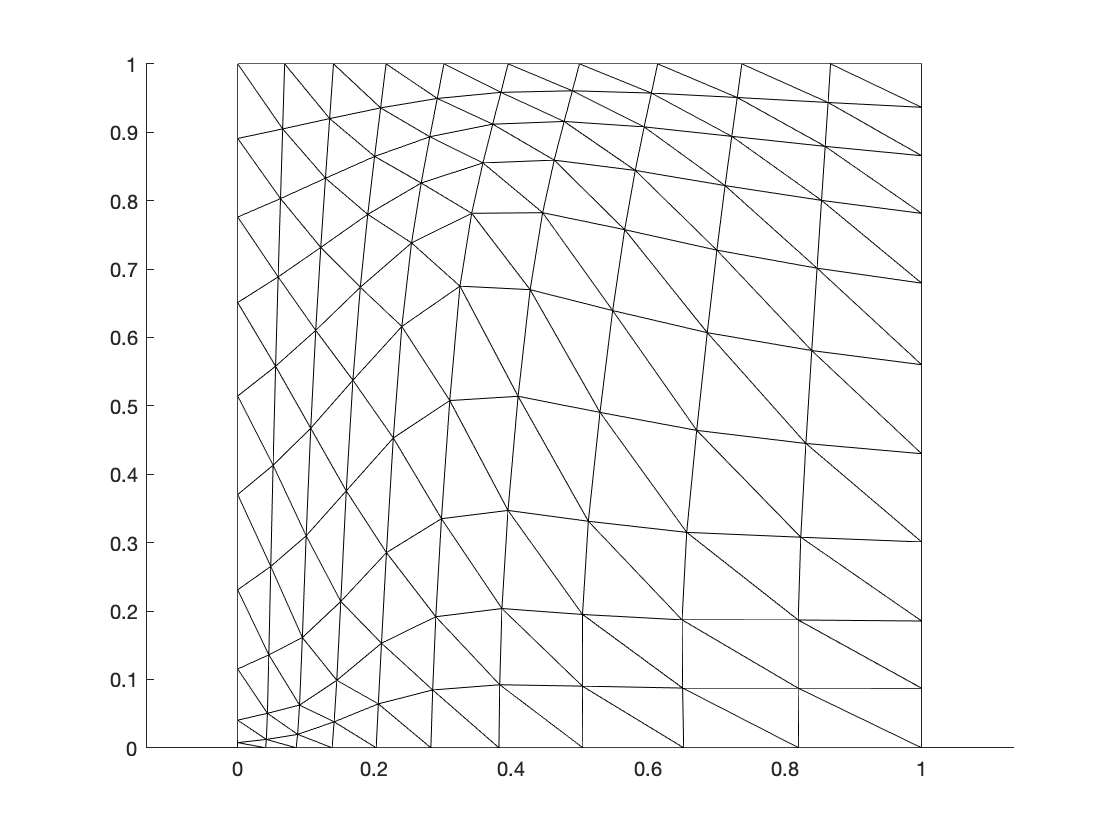}}
 
 \caption{
 three-point registration; visualization of the deformed meshes for three choices of the objective function in \eqref{eq:morozov_registration}.
}
 \label{fig:elasticity_bad}
  \end{figure}

Figure \ref{fig:elasticity_bad} shows the results of \eqref{eq:morozov_registration} with 
(a) linear-elasticity objective,
(b)  H2 objective, and
(c) neo-hookean objective.
We observe that both the  linear-elasticity model and the 
neo-hookean model fail to deliver a proper deformed mesh for this test case, while the H2 objective does.
The linear-elasticity model  also leads to a non-bijective mapping --- the minimum Jacobian determinant $J_{\rm min}$ is equal to $0.012$ for H2, to $0.47$ for the neohookean model and to $-0.09$ for linear elasticity.
We also remark that the introduction of higher-order derivatives in the objective function has the effect of smoothening the mapping $\Phi$: we conjecture  that this feature of the approach might be important in the pMOR framework to improve the compressibility of the solution manifold.

\end{document}